%% file: paper20rev.tex
\magnification=\magstep 1
\baselineskip=15pt     
\parskip=3pt plus1pt minus.5pt
\overfullrule=0pt
\font\hd=cmbx10 scaled\magstep1
\input defs.tex

\def\tr{{\mathop{\rm Tr}}}
\input epsf.tex
\input amssym.tex

\centerline{\hd Topological Mirror Symmetry.}
\medskip
\centerline{\it Mark Gross\footnote{*}{Supported in part by NSF grant 
DMS-9700761}}
\medskip
\centerline{September 2nd, 1999, {\it revised} November 2nd, 1999}
\medskip
\centerline{Mathematics Institute}
\centerline{University of Warwick}
\centerline{Coventry, CV4 7AL}
\centerline{mgross@maths.warwick.ac.uk}
\bigskip
\bigskip
{\hd \S 0. Introduction.}

The Strominger-Yau-Zaslow conjecture proposes that mirror symmetry
can be explained by the existence, in a mirror pair of
Calabi-Yau manifolds, of dual special Lagrangian $T^n$-fibrations.
(See [19,8,6,7] for further clarification of this conjecture.)
Recently, Zharkov in [21] proved that non-singular Calabi-Yau hypersurfaces
in toric varieties have topological $T^n$-fibrations, and Ruan
in [17] has shown further that the quintic threefold
has a Lagrangian torus fibration.
He in addition announced in that paper an extension of these results
to arbitrary Calabi-Yau threefolds which are hypersurfaces in toric
varieties. This gives initial evidence for the SYZ conjecture.

In order for the SYZ conjecture to be useful in constructing mirror 
manifolds, one should be able to produce the mirror manifold by
dualizing and compactifying torus fibrations with no a priori 
knowledge of the mirror manifold.
Typically, one expects to start with a fibration
$f:X\rightarrow B$ whose general fibre is a torus, and discriminant
locus $\Delta\subseteq B$ of $f$ with $B_0=B-\Delta$. One then has a diagram
$$\matrix{X_0&\hookrightarrow&X\cr
\mapdown{f_0}&&\mapdown{f}\cr
B_0&{\smash{\mathop{\hookrightarrow}\limits^i}}&B\cr}$$ 
with $f_0$ being a $T^n$-fibre bundle.
The dual of $f_0$ can be described as the torus bundle
$\check f_0:\check X_0=R^1f_{0*}({\bf R}/\boldz)\rightarrow B_0$. A dual
of $f$ is then a compactification of $\check f_0$ to $\check f:\check X
\rightarrow B$. Thus important questions are: {\it Does there exist a class
of torus fibrations for which duals can be constructed? Do these fibrations
exist on Calabi-Yau manifolds such as the quintic threefold? If they do,
does dualizing produce the mirror?} In this paper we answer
these questions affirmatively in the three-dimensional case. 
This proves that in at least
one non-degenerate case, the SYZ conjecture explains
mirror symmetry on a topological level.

To construct a compactification
in some reasonably unique way, one expects that one needs the cohomology of the
looked-for singular fibres to be determined by the fibration
$\check X_0\rightarrow B_0$.
More specifically, in [6] and [7] a property called $G$-simplicity
was introduced. A $T^n$-fibration $f:X\rightarrow B$
is {\it $G$-simple}, for $G$ an abelian group, if $i_*R^pf_{0*}G
=R^pf_*G$ for all $p$. This essentially says that the cohomology
of the singular fibres is determined by the monodromy of the local
system $R^qf_{0*} G$ on $B_0$. Without this condition, it would appear
to be difficult to construct a dual, and it is certainly difficult
to relate the cohomology of $X$ with that of a dual fibration. 
By contrast, if $f:X\rightarrow B$ and $\check f:\check X\rightarrow B$
are dual ${\bf R}$-simple $T^3$ fibrations, then it was shown in
[6] that $H^{2*}(X,{\bf R})\cong H^3(\check X,{\bf R})$ and
$H^3(X,{\bf R})\cong H^{2*}(\check X,{\bf R})$, the desired relationship
of mirror symmetry.

For the fibrations of [17] and [21],
$\Delta$ is codimension one, and then such $f:X\rightarrow B$ cannot be simple,
making it unlikely one can construct the dual of $f$ in a natural way. 
Ruan discusses
this issue in [17], and points out that it should not be difficult to
modify his fibration topologically to obtain one with a codimension two
locus. He also announces that he can do this in such a way that the fibration
is Lagrangian. This will come a long way towards addressing the above problems.
However, the hypothetical fibration he analyzes in [17] in some depth
with codimension two discriminant locus is also not 
simple, because the fibres are not all irreducible. It was argued in [6]
that for a generic choice of complex structure all fibres of a special 
Lagrangian $T^3$-fibration would be irreducible. Thus it should be possible
to construct a simple torus fibration on the quintic.

In this paper we delve more deeply into
the structure of $T^3$-fibrations with sufficient regularity hypotheses,
similar to the
ones assumed in [7]. We call such fibrations ``well-behaved''; see
Definitions 1.1 and 1.2. These fibrations in particular are $G$-simple
for $G=\boldz$ and $\boldz/n\boldz$, which prove to be
very strong conditions. In \S 1,
we study the local topology of well-behaved fibrations,
and focus on what we term {\it semi-stable} singular fibres,
i.e. singular fibres for which the local monodromy group is unipotent.
We obtain a classification of possible monodromy groups
of such singular fibres, and furthermore give a precise description of the
discriminant locus near such fibres.
In \S 2 we use this information to
show how to construct well-behaved duals of such fibrations.
The main result of \S 2 is

\proclaim Theorem 0.1. If $f:X\rightarrow B$ is a well-behaved
$T^3$-fibration, with only semi-stable fibres, then $f$
has a well-behaved dual $\check f:\check X\rightarrow B$.

This gives a class of $T^3$-fibrations whose duals can now be constructed.

\S 3 is a bit off of the main thrust of the paper. There we give a
construction of topological $T^{n-1}\times {\bf R}$-fibrations on crepant
resolutions of $n$-dimensional toric Gorenstein canonical singularities.
We make a connection with the local mirror symmetry of
[4], and use this construction to motivate some of the choices
made in \S 4, as well as aid in the identification of the mirror quintic.

In \S 4, we study the specific example of the quintic. We prove the main
theorem of the paper:

\proclaim Theorem 0.2. The quintic threefold $X\subseteq\Pfour$ has a
well-behaved $T^3$-fibration with semistable fibres $f:X\rightarrow B$.
If $\check f:\check X\rightarrow B$ is the dual fibration given by
Theorem 0.1, then $\check X$ is diffeomorphic to a specific
non-singular minimal model of the mirror quintic.

This theorem is proved by first constructing a $T^3$-fibrations $f:X\rightarrow
B$ ``by hand,'' making use of Ruan's description of monodromy in [17],
and then showing $X$ is diffeomorphic to the quintic
using Wall's classification theory [20]. This classification theory
is also applied
to show $\check X$ is diffeomorphic to the mirror quintic.

In doing this, we also discover a beautiful geometric explanation
for the phenomenon of topology change in mirror symmetry
discovered by Aspinwall, Greene and Morrison. In [1], 
arguments involving calculations of Yukawa couplings
at large complex structure limits suggested
that birationally equivalent Calabi-Yau threefolds have the same
mirror, but correspond to different large complex structure
limit points in the complex moduli of the mirror. In general it is
conjectured that special Lagrangian fibrations only exist near
large complex structure limit points. Thus one might well expect to find
two distinct $T^3$-fibrations $f_1:X\rightarrow B$ and $f_2:X\rightarrow B$,
with duals $\check f_1:\check X_1\rightarrow B$ and $\check f_2:\check X_2
\rightarrow B$, such that $\check X_1$ and $\check X_2$ are topologically 
distinct manifolds, perhaps related by a flop. In fact, we will see this
explicitly. In defining the $T^3$-fibration in the quintic, there will
be some ambiguity in the construction which gives rise to distinct
fibrations, and whose duals are related by flops.
 
There are several unsatisfactory aspects to the paper.
It will be clear from the constructions of \S 4 that they should generalize
to arbitrary Calabi-Yau
hypersurfaces in toric varieties. However, the proof that we
have constructed the quintic and the mirror quintic is, at this point,
very ad hoc. Applying Wall's classification result means that one
has to identify the cubic intersection forms, and this means much
messing about with explicit representatives of cohomology classes
and orientations. Before generalising to arbitrary toric hypersurfaces,
a more natural proof that the fibrations we have constructed in
fact coincide with the desired hypersurfaces should be found. Another point
missing is that we make no connections with special Lagrangian geometry,
and at this point cannot argue that the fibrations we have
constructed actually resemble hypothetical special Lagrangian 
fibrations.

Finally, I note that shortly after the first version of this paper
was released, Ruan released a paper [18] announcing a proof of
existence of Lagrangian $T^3$-fibrations on the quintic and
mirror quintic which are provably dual to each other via a comparison
of monodromy. The technique used is quite different from
the one given here, and makes use of a perturbation of the
fibrations obtained by Ruan's symplectic flow technique. I expect
however that these fibrations coincide with my own on the topological
level.

I would like to thank D. Joyce, D. Matessi, M. Micallef and P.M.H. Wilson
for useful discussions concerning the SYZ conjecture, and the hospitality
of M. Kobayashi at RIMS where part of this work was inspired by
a talk of I. Zharkov. I also thank
Tenryuji for a key idea.

{\hd \S 1. Torus Fibrations and Semistable Fibres.}

We are interested in studying topological $T^n$-fibrations
$f:X\rightarrow B$, i.e continuous proper maps with
connected fibres between topological manifolds
whose general fibres are tori. However, without any restrictions on
the nature of such a map, there will be no useful theorems one can
prove about such fibrations. We begin by giving the restrictions
on such fibrations we will use in this paper. These restrictions
are motivated by results of [7] on the structure of special Lagrangian
fibrations.

\proclaim Definition 1.1. Let $f:X\rightarrow B$ be a continuous, proper
mapping of topological manifolds with connected fibres, 
$\dim X=2n$, $\dim B=n$, such that
for a dense open set $B_0\subseteq B$, the fibre $X_b$ is an $n$-torus
for all $b\in B_0$. We say $f$ is {\it admissible} if there exists
an open subset $X^{\#}\subseteq X$ with $f^{\#}:X^{\#}\rightarrow B$
the restriction of $f$ to $X^{\#}$, $Crit(f):=X-X^{\#}$,
 satisfying the following properties:
\item{(1)} $Sing(X_b):=Crit(f)\cap X_b$ is a union of a finite number of
locally closed submanifolds of dimension at most $n-2$,
and $\Delta:=f(Crit(f))=B-B_0$.
\item{(2)} For each $b\in B$, there exists an open neighbourhood
$U\subseteq B$ of $b$, a rank $n$ vector bundle $\F$ on
$U$, and an exact sequence
$$\exact{(R^{n-1}_cf^{\#}_*\boldz)|_U}{\F}{f^{\#-1}(U)}.$$
Furthermore, the map $\F\rightarrow f^{\#-1}(U)$ is a local isomorphism.
\item{(3)} $R^nf_*\boldz\cong\boldz$.

The last condition in this definition, besides providing a uniform
orientation on fibres (which would be the case for a special Lagrangian
fibration), also ensures the singular fibres are irreducible in
a certain sense. We don't expect this to be the case for all special
Lagrangian fibrations (for example an elliptic fibration can have reducible
fibres) but it was argued in [6] that gnerically this should
be the case for $n\le 3$.

We will also need some regularity assumptions involving
the discriminant locus $\Delta$, which in dimensions 2 and
3 are encapsulated in

\proclaim Definition 1.2. 
Let $f:X\rightarrow B$ be a $T^n$-fibration,
$n=2$ or $3$. Let $\Delta$ be the discriminant locus of $f$. We say
$f$ is {\it well-behaved} if it is admissible and either
\item{(1)} $n=2$ and $\Delta$ is discrete; or,
\item{(2)} $n=3$ and $\Delta$ can be written as a disjoint union of 
two sets, $\Delta_g$ and $\Delta_d$, where $\Delta_d$ is a discrete
set (the ``dissident'' points)
and $\Delta_g$ is a locally closed 
topological 1-submanifold (the ``generic'' points).
For all $b\in \Delta_g$, there exists an open neighbourhood
$U$ of $b\in B$, a homeomorphism $\alpha:U\rightarrow D\times (0,1)$, $D$ a
two-dimensional disk, a homeomorphism $\alpha':f^{-1}(U)\rightarrow
X'\times (0,1)\times S^1$ with $X'$ a four-manifold, 
$g:X'\rightarrow D$
a well-behaved $T^2$-fibration, such that the diagram
$$\matrix{f^{-1}(U)&\mapright{\alpha'}&X'\times (0,1)\times S^1\cr
\mapdown{f}&&\mapdown{f'}\cr
U&\mapright{\alpha}&D\times (0,1)\cr}$$
is commutative, where $f'$ is the composition of projection
onto $X'\times (0,1)$ and $g\times id$.
Also $\alpha(\Delta_g\cap U)=\{v\}\times (0,1)$ for some point
$v\in D$. Furthermore, for each point $b\in\Delta_d$,
there is an open neighbourhood $U$ of $\Delta_d$ such that 
there is a commutative diagram
$$\matrix{U-\Delta&\hookrightarrow&U-\Delta_d\cr
\mapdown{\cong}&&\mapdown{\cong}\cr
(S^2-\{p_1,\ldots,p_m\})\times (0,1)&\hookrightarrow&
S^2\times (0,1).\cr}$$
In addition, there is a five-manifold $X'$ and a map $g:X'\rightarrow S^2$
along with a commutative diagram
$$\matrix{f^{-1}(U-\Delta_d)&\mapright{\cong}&X'\times (0,1)\cr
\mapdown{f}&&\mapdown{g\times id}\cr
U-\Delta_d&\mapright{\cong}&S^2\times (0,1).\cr}$$

Let me emphasize that at this point we do not know that special Lagrangian
fibrations are either admissible or well-behaved, but the arguments of
[7] \S 1 and 2
show integral special Lagrangian fibrations
are admissible given sufficient regularity assumptions, e.g.
$f$ real analytic. (However, this may be too strong a hypothesis.)
On the other hand, the fibrations we will construct in this paper are
well-behaved, so we do not require any further justification
of these definitions. 

We make the observation that these restrictions are sufficient to prove
simplicity.

\proclaim Theorem 1.3. If $f:X\rightarrow B$ is a well-behaved 
$T^2$-fibration and $G$ is an abelian group, then $f$
is $G$-simple. If $f:X\rightarrow B$ is a well-behaved
$T^3$-fibration, then $f$ is $\boldz$-simple and 
$\boldz/n\boldz$-simple for all $n$. 

Proof. As remarked at the beginning of [7], \S 3, $f_*G=G$ since
the fibres of $f$ are connected, and $R^nf_*G=G$ follows from Definition
1.1, (3).
Thus $i_*R^pf_{0*}G=R^pf_*G$ for $p=0,n$, where $i:B-\Delta\rightarrow B$
is the inclusion. For $p=1$, this follows from
[7], Theorem 3.3, where it can be checked that the conditions of
well-behaved are sufficient for the proof of that theorem in [7] to work.
Thus, in particular, a well-behaved $T^2$-fibration
is $G$-simple.

Furthermore, the proof of [7] Theorem 3.5 also applies to well-behaved
$T^3$-fibrations, and thus such $f$ are $\boldz$-simple.
To show $f$ is $\boldz/n\boldz$-simple, we still need to show
that $i_*R^2f_{0*}\boldz/n\boldz=R^2f_*\boldz/n\boldz$. The
proof of [7], Theorem 3.5, and especially Lemma 3.4,
does not carry over to
the coefficient group $\boldz/n\boldz$. Thus a different proof must be given.

Let $i_0:B-\Delta_d\rightarrow B$ be the inclusion. We first need
to show the equivalent of [7], Lemma 3.4, namely that 
$$R^2f_*\boldz/n\boldz\cong i_{0*}i_0^*R^2f_*\boldz/n\boldz.$$
Let $b_0\in\Delta_d$, and let $U\subseteq B$ be an open neighbourhood
of $b_0$ as given in Definition 1.2. Thus $U-\{b_0\}\cong
S^2\times (0,1)$, $f^{-1}(U-\{b_0\})\cong X'\times (0,1)$, and the
map $f$ is induced by $g:X'\rightarrow S^2$.
Put $X^*=f^{-1}(U-\{b_0\})$. Using the relative cohomology
exact sequence, the fact that $H^i(f^{-1}(U),\boldz)\cong H^i(X_{b_0},\boldz)$,
and $H^i_{X_{b_0}}(f^{-1}(U),\boldz)=H_{6-i}(X_{b_0},\boldz)$, one sees that
$$H^i(X^*,\boldz)=\cases{\boldz&$i=0,5$\cr
\boldz^{b_1}&$i=1,4$\cr
\boldz^{b_2+1}&$i=2,3$,\cr}\leqno{(1.1)}$$
where $b_i=\rank H^i(X_{b_0},\boldz)$.
The only part of this calculation which requires comment is $i=2$ and
$3$. Here we have a commutative diagram
$$\matrix{
&&&&H^2(U-\{p\},\boldz)&\mapright{\alpha'}&H^3_p(U,\boldz)=\boldz\cr
&&&&\mapdown{}&&\mapdown{\beta}\cr
0&\mapright{}&H^2(f^{-1}(U),\boldz)&\mapright{}&
H^2(f^{-1}(U-\{p\}),\boldz)&\mapright{\alpha}&H^3_{X_b}(f^{-1}(U),\boldz)=\boldz
\cr}$$
Since $\alpha'$ and $\beta$ are isomorphisms,
$\alpha$ is surjective.

While we will need to make use of this full calculation later in the paper,
the only part we need now is the conclusion that
$H^*(X^*,\boldz)$ is torsion-free, and so
$H^*(X^*,\boldz/n\boldz)=H^*(X^*,\boldz)\otimes\boldz/n\boldz$.

Now look at the Leray spectral sequence for $f:X^*\rightarrow B^*=U-\{p\}$,
with coefficients in $\boldz$. This looks like
$$\matrix{\boldz&0&\boldz\cr
H^0(B^*,R^2f_*\boldz)&H^1(B^*,R^2f_*\boldz)&H^2(B^*,R^2f_*\boldz)\cr
H^0(B^*,R^1f_*\boldz)&H^1(B^*,R^1f_*\boldz)&H^2(B^*,R^1f_*\boldz)\cr
\boldz&0&\boldz\cr}$$
Admissible fibrations have local sections, so
we can assume that $f:X^*\rightarrow B^*$ has a section. Then standard
arguments as in [6], Lemma 2.4 
show that the only possible non-zero differential is
$d:H^0(B^*,R^2f_*\boldz)\rightarrow H^2(B^*,R^1f_*\boldz)$. However,
in the proof of [7], Lemma 3.4, it was shown that $H^2(X^*,\boldz)
\rightarrow H^0(B^*,R^2f_*\boldz)$ was surjective, and hence $d$ 
is zero.
Thus the spectral sequence degenerates at the $E_2$ term.
In particular, $H^2(B^*,R^if_*\boldz)$ is torsion-free since
$H^{i+2}(X^*,\boldz)$ is torsion-free.

Now suppose that $H^1(B^*,R^2f_*\boldz)$ were not torsion-free.
First note that because of the product structure of the map $X^*\rightarrow
B^*$, $H^i(B^*,R^jf_*\boldz)\cong H^i(S^2,R^jg_*\boldz)$. We can then use
Poincar\'e-Verdier duality to compare $H^1(S^2,R^2g_*\boldz)_{tors}$
with $H^2(S^2,R^1g_*\boldz)_{tors}$.
By $\boldz$-simplicity, $\lhom(R^1g_*\boldz,\boldz)
\cong R^2g_*\boldz$, and thus by the local-global Ext spectral sequence,
$H^1(S^2,R^2g_*\boldz)$ injects into $\ext^1(R^1g_*\boldz,\boldz)$,
which in turn, by Poincar\'e-Verdier duality (as applied in the proof
of [7] Theorem 3.9, see especially (3.2) there), 
is isomorphic to a $\boldz$-module
$M$ fitting in the exact sequence
$$\exact{\ext^1(H^2(S^2,R^1g_*\boldz),\boldz)}{M}{\hom(H^1(S^2,R^1g_*\boldz),
\boldz)}.$$
Thus, if $H^1(S^2,R^2g_*\boldz)$ contains torsion, so does $M$,
and hence so does $H^2(S^2,R^1g_*\boldz)\cong H^2(B^*,R^1f_*\boldz)$, 
contradicting the torsion-freeness
of this $\boldz$-module. Thus $H^1(B^*,R^2f_*\boldz)$ is torsion-free.
It then follows from the exact sequence
$$0\mapright{} R^2f_*\boldz \mapright{\cdot n} R^2f_*\boldz
\mapright{} R^2f_*\boldz/n\boldz\mapright{} 0$$
that $H^0(B^*,R^2f_*\boldz/n\boldz)=H^0(B^*,R^2f_*\boldz)
\otimes \boldz/n\boldz$. Since the restriction map
$$H^2(X_b,\boldz)\cong H^0(U,R^2f_*\boldz)\rightarrow H^0(B^*,
R^2f_*\boldz)$$ is an isomorphism  
by $\boldz$-simplicity, it then follows that the map $$H^2(X_b,
\boldz/n\boldz)\cong H^0(U,R^2f_*\boldz/n\boldz)
\rightarrow H^0(B^*,R^2f_*\boldz/n\boldz)$$ is an isomorphism, 
so that
$i_{0*}i_0^*R^2f_*\boldz/n\boldz=R^2f_*\boldz/n\boldz$.

Finally, to show $\boldz/n\boldz$-simplicity, we have to show that if $i_1:
B-\Delta\hookrightarrow B-\Delta_d$, then the natural map
$$(R^2f_*\boldz/n\boldz)|_{B-\Delta_d}\mapright{\phi} i_{1*}i_1^*
R^2f_*\boldz/n\boldz$$
is an isomorphism. The kernel of $\phi$ is $\HH^0_{\Delta_g}(B-\Delta_d,
R^2f_*\boldz/n\boldz)=0$, since $\HH^i_{\Delta_g}(B-\Delta_d,R^2f_*\boldz)=0$
for $i=0,1$, by $\boldz$-simplicity. The cokernel of $\phi$ is
$$\HH^1_{\Delta_g}(B-\Delta_d,R^2f_*\boldz/n\boldz)=\ker(
\HH^2_{\Delta_g}(B-\Delta_d,R^2f_*\boldz)\mapright{\cdot n}
\HH^2_{\Delta_g}(B-\Delta_d,R^2f_*\boldz)).$$
Thus to show this is zero for any $n$, it is enough to show it for $n$
prime.

To show this, we take a point $b_0\in\Delta_g$, and $\gamma$ a simple
loop around $b_0$ based at a nearby point $b$. Let $T:H^2(X_b,\boldz/n\boldz)
\rightarrow H^2(X_b,\boldz/n\boldz)$ be the induced monodromy transformation.
The map $\phi$ on stalks at $b_0$ is
$$H^2(X_{b_0},\boldz/n\boldz)\rightarrow H^2(X_b,\boldz/n\boldz)^T.$$
The map is injective, so we just need to show
$\dim H^2(X_b,\boldz/n\boldz)^T=\dim H^2(X_{b_0},\boldz/n\boldz)$.
But $\ker T^t-I=H^1(X_{b_0},\boldz/n\boldz)$, and since $X_{b_0}
=X'_{b_0}\times S^1$, $H^1(X_{b_0},\boldz/n\boldz)$ is the same
dimension as $H^2(X_{b_0},\boldz/n\boldz)$ by K\"unneth, so $\dim\ker(T-I)=
\dim\ker(T^t-I)=\dim H^1(X_{b_0},\boldz/n\boldz)=\dim H^2(X_{b_0},\boldz/
n\boldz)$, as desired.
$\bullet$

We now show how simplicity enables one to obtain strong restrictions on 
monodromy about singular fibres in torus fibrations.

\proclaim Theorem 1.4. Let $f:X\rightarrow B$ be a well-behaved 
$T^2$-fibration, $b_0\in\Delta\subseteq B$, $b$ a nearby point,
and $T:H_1(X_b,\boldz)\rightarrow H_1(X_b,\boldz)$ the monodromy
transformation about a simple loop around $b_0$ based at $b$. Then
there is a basis $e_1,e_2$  of $H_1(X_b,\boldz)$
in which $T=\pmatrix{1&1\cr 0&1}$, $\pmatrix{0&-1\cr 1&1\cr}$ or
$\pmatrix{0&-1\cr 1&3\cr}$. In the first case, the singular fibre is of
type $I_1$ (by which we mean the one-point compactification of 
$S^1\times{\bf R}$),
while in the second and third cases the singular fibre is
topologically an $S^2$.

Proof. By simplicity, $H_1(X_b,\boldz)^T=H_1(X_{b_0},\boldz)$,
and since $X_{b_0}$ is a singular fibre, $\rank H_1(X_{b_0},\boldz)=0$ or $1$.
Note that since $Sing(X_{b_0})$ is connected (as can be shown by the same
argument as in Lemma 2.5 of [7]), $Sing(X_{b_0})$ consists of one point,
and so $X_{b_0}$ is just the one point compactification of $X_{b_0}^{\#}$.
Since $X_{b_0}^{\#}\cong {\bf R}^2/H_1(X_{b_0}^{\#},\boldz)$, 
this compactification
is either an $I_1$ fibre or $S^2$, topologically.

{\it Case 1.} $\rank H_1(X_{b_0},\boldz)=1$. Then $\rank\ker T-I=1$,
and choosing a basis $e_1,e_2$ with $e_1$ invariant
yields $T=\pmatrix{1&a\cr 0&1\cr}$ (using also the fact that
$T\in SL_2(\boldz)$). Now by $\boldz/n\boldz$-simplicity,
$H_1(X_b,\boldz/n\boldz)^{T_n}=H_1(X_{b_0},\boldz/n\boldz)=\boldz/n\boldz$,
where $T_n\in SL_2(\boldz/n\boldz)$ is $T\mod n$. Thus $\rank\ker T_n-I
=1$ for all $n\ge 2$, and in particular $a=\pm 1$. Changing the sign
of $e_2$ if necessary, we obtain the first matrix. In this case the
compactification is an $I_1$ type fibre.

{\it Case 2.} $H_1(X_{b_0},\boldz)=0$. In this case, there are no
invariant vectors modulo $n$, so in particular $T$ does not have
the eigenvalue 1 modulo any $n$. Now the characteristic polynomial
of $T$ is $f(x)=x^2-\tr(T)x+1$, and we require
$f(1)\not\equiv 0\mod n$ for any $n\ge 2$. Thus
$f(1)=2-\tr(T)=\pm 1$, or $\tr(T)=1$ or $3$. 

First assume $\tr(T)=1$. Then $\boldz[T]=\boldz[(1+\sqrt{-3})/2]$,
which is a PID. Thus $H_1(X_b,\boldz)$ is a free $\boldz[T]$
module, and hence is $H_1(X_b,\boldz)=\boldz[T]e_1$ for some $e_1\in
H_1(X_b,\boldz)$. Let $e_2=Te_1$. Then $e_1,e_2$ form a basis
for $H_1(X_b,\boldz)$, and in this basis the matrix for $T$
is $\pmatrix{0&-1\cr 1&1\cr}$ as desired.

Next, suppose $\tr(T)=3$. Then $\boldz[T]=\boldz[(1+\sqrt{5})/2]$,
the number ring of ${\bf Q}(\sqrt{5})$. This is also a PID, so the
same argument allows us to find a basis in which 
$T=\pmatrix{0&-1\cr 1&3\cr}$. In both these cases, $X_{b_0}$ is
topologically a sphere.
$\bullet$

{\it Remark 1.5.} It does not seem possible to rule out the last
matrix in Theorem 1.4
using only topological means, as such a topological fibration appears
to exist. Of course, it does not occur in the elliptic curve
case, or equivalently the special Lagrangian case.
\bigskip

We now move on to the three-dimensional case.

\proclaim Definition 1.6. Let $f:X\rightarrow B$ be a well-behaved
$T^3$-fibration, $b_0\in\Delta$ a point. The {\it monodromy
group} of $f$ at $b_0$ is the image of $\pi_1(U-\Delta,b)
\rightarrow SL(H_1(X_b,\boldz))$, where $U$ is an open neighbourhood
of $b_0$ as given in Definition 1.2, and $b\in U$ a base point.

\proclaim Theorem 1.7. Let $f:X\rightarrow B$ be a well-behaved
$T^3$-fibration. If $b_0\in\Delta_g$, then $X_{b_0}$ is homeomorphic to
$I_1\times S^1$ or $S^2\times S^1$. 
There is a basis of $H_1(X_b,\boldz)$
for $b$ near $b_0$ in which the monodromy group is generated by 
$\pmatrix{1&1&0\cr 0&1&0\cr 0&0&1\cr}$,
$\pmatrix{0&-1&0\cr 1&1&0\cr 0&0&1\cr}$, or
$\pmatrix{0&-1&0\cr 1&3&0\cr 0&0&1\cr}$.

Proof. This follows immediately from the definition of well-behaved
and Theorem 1.4.
$\bullet$

\proclaim Definition 1.8. Let $f:X\rightarrow B$ be a well-behaved 
$T^3$-fibration.  Let $b_0\in B$ be a point in the discriminant
locus $\Delta$ of $f$. We say $X_{b_0}$ is {\it semistable} if the monodromy
group at $b_0$ is unipotent.

This terminology is in analogy with terminology for elliptic curves.

\proclaim Definition 1.9. 
Let $f:X\rightarrow B$ be a well-behaved $T^3$-fibration,
$b_0\in B$. We say the fibre $X_{b_0}$ is of type $(b_1,b_2)$
if $b_i=\rank H^i(X_{b_0},\boldz)$, $i=1,2$.

\proclaim Theorem 1.10. Let $f:X\rightarrow B$ be a well-behaved 
$T^3$-fibration, and let $b_0\in \Delta$ with $X_{b_0}$ be a semistable fibre with
monodromy group $G$. Then either 
\item{(1)} $X_{b_0}$ is of type $(2,2)$ and $G$ is conjugate
to 
$$\left\{\pmatrix{1&0&a\cr 0&1&0\cr 0&0&1\cr}\bigg | a\in\boldz\right\}.$$
\item{(2)} $X_{b_0}$ is of type $(2,1)$ and $G$ is conjugate
to 
$$\left\{\pmatrix{1&a&b\cr 0&1&0\cr 0&0&1\cr}\bigg | a,b\in\boldz\right\}.$$
\item{(3)} $X_{b_0}$ is of type $(1,2)$ and $G$ is conjugate
to 
$$\left\{\pmatrix{1&0&a\cr 0&1&b\cr 0&0&1\cr}\bigg | a,b\in\boldz\right\}.$$
\item{(4)} $X_{b_0}$ is of type $(1,1)$ and $G$ is conjugate
to 
$$\left\{\pmatrix{1&a&b\cr 0&1&c\cr 0&0&1\cr}\bigg | a,b,c\in\boldz\right\}.$$

Proof. Let $b$ be a base-point near $b_0$.
By the theorem of 17.5 of [12], $H_1(X_b,{\bf C})$
contains an element invariant under $G$, and thus $b_2(X_{b_0})\ge 1$.
In particular, there is a $G$-invariant $e_1\in H_1(X_b,\boldz)$.
Applying the theorem again to
$H_1(X_b,{\bf C})/{\bf C}e_1$, one finds a two dimensional
subspace $V\subseteq H_1(X_b,{\bf C})$ containing
$e_1$, and invariant under $G$. Thus $V$ yields an invariant
element of $\bigwedge^2 H_1(X_b,{\bf C})$. Since $H_2(X_b,{\bf C})\cong
\bigwedge^2 H_1(X_b,{\bf C})$, 
$H_2(X_b,{\bf C})^G=H^1(X_b,{\bf C})^G
\not=0$, so $b_1(X_{b_0})\not=0$. Furthermore,
some element of $H_2(X_b,\boldz)^G$ is of the form $e_1\wedge e_2$
for some $e_2$ such that $e_1$, $e_2$ span a primitive sublattice of
$H_1(X_b,\boldz)$. Choosing $e_3$ so that $e_1,e_2,e_3$ span $H_1(X_b,\boldz)$
we see that in this basis $G$ is a subgroup of
$$G'=
\left\{\pmatrix{1&a&b\cr 0&1&c\cr 0&0&1\cr}\bigg | a,b,c\in\boldz\right\}.$$
In particular, $b_1(X_{b_0}),b_2(X_{b_0})>0$, so only the types listed
can occur. (It was shown in [7], \S 3, that $b_i(X_{b_0})=3$ for $i=1$ or
$2$ implies $X_{b_0}$ is non-singular.)

For $X_{b_0}$ of type $(2,2)$, we can choose $e_1$ and $e_2$
to be a basis of the invariant cycles of $H_1(X_b,\boldz)$
so that $G$ is in fact contained in
$$H=
\left\{\pmatrix{1&0&a\cr 0&1&b\cr 0&0&1\cr}\bigg | a,b\in\boldz\right\}
\cong\boldz^2.$$
$G$ is then a subgroup of $H$. Since $G$ acts on $H_2(X_b,\boldz)
\cong \dual{H_1(X_b,\boldz)}$ via the transpose inverse matrices,
we see the fact that $H_2(X_b,\boldz)^G\cong \boldz^2$ implies that
$G\cong\boldz\subseteq\boldz^2$, i.e. is generated by some matrix
$T=\pmatrix{1&0&a\cr 0&1&b\cr 0&0&1\cr}$. If $n={\rm gcd}(a,b)\not=1$,
then $T\equiv I\mod n$, violating $\boldz/n\boldz$-simplicity.
Thus $a$ and $b$ are relatively prime, and a change of basis
gives the desired form for $G$.

For type $(1,2)$, $e_1$ and $e_2$ are again chosen to be invariant
as above, but now $G\subseteq H$ must be a sublattice
of rank 2 to ensure $H_2(X_b,\boldz)^G\cong \boldz$. Let $T_1=
\pmatrix{1&0&a\cr 0&1&b\cr 0&0&1\cr}$, $T_2=
\pmatrix{1&0&c\cr 0&1&d\cr 0&0&1\cr}$ be generators of this
sublattice. Let $e_1^*,e_2^*,e_3^*$ be the dual basis to $e_1,e_2,e_3$.
From $\boldz/n\boldz$-simplicity it follows that
$H_2(X_b,\boldz/n\boldz)^G=(\boldz/n\boldz)\cdot e_3^*$ for all $n$,
hence modulo $n$, the kernel of the matrix $\pmatrix{a&b\cr c&d\cr}$
is always zero. Thus $\det\pmatrix{a&b\cr c&d\cr}\not\equiv 0\mod n$
for any $n$, so $\det\pmatrix{a&b\cr c&d\cr}=\pm 1$. Thus $T_1$ and
$T_2$ generate $H$, as desired.

The argument for type $(2,1)$ is identical, interchanging
the roles of $H_1(X_b,\boldz)$ and $H_2(X_b,\boldz)$.

Finally, if $X_{b_0}$ 
is type $(1,1)$, then $G$ must be generated by transformations
$T_1,\ldots,T_n$ where $T_i$ is the monodromy transformation about a curve in 
$\Delta_g$, and hence, being unipotent, is conjugate to $\pmatrix{1&1&0\cr
0&1&0\cr 0&0&1\cr}$, by Theorem 1.7.
If $T_i$ is written in the basis $e_1,e_2,e_3$ so that $G\subseteq G'$,
then $T_i=\pmatrix{1&a_i&b_i\cr 0&1&c_i\cr 0&0&1\cr}$
for some $a_i,b_i,c_i\in\boldz$. However the requirement that
$\rank\ker T_i-I=2$ guarantees that either $a_i=0$ or $c_i=0$ for each $i$.
If $a_i=0$ for all $i$, then we must in fact be in the $(1,2)$
case; similarly, if $c_i=0$ for all $i$, then we must be in the $(2,1)$
case.

Now the group $G'$ is in fact a central extension
$$0\mapright{}\boldz\mapright{} G'\mapright{\psi}\boldz^2\mapright{}0,$$
where $\psi\pmatrix{1&a&b\cr 0&1&c\cr 0&0&1\cr}=(a,c)$ (in fact $G'$
is the Heisenberg group over $\boldz$). What we have now shown is that in the
$(1,1)$ case $\psi(G)$ must be a rank 2 sublattice generated by
$(a_0,0)$ and $(0,c_0)$ for some $a_0,c_0\in\boldz$. Furthermore,
$a_0=\pm 1$ and $c_0=\pm 1$, as can be seen via $\boldz/n\boldz$
simplicity. Thus $G$ contains elements of the form $\pmatrix{1&1&b_1\cr
0&1&0\cr 0&0&1\cr}$ and $\pmatrix{1&0&b_2\cr 0&1&1\cr 0&0&1\cr}$.
But these two elements generate $G'$, so $G=G'$ as desired. $\bullet$

Next we give a more detailed structure result. If $b_0\in\Delta_d$ is
a dissident point, and $U\subseteq B$ the open neighbourhood of
$b_0$ promised by Definition 1.2, we can write 
$U-\Delta=(S^2-\{p_1,\ldots,p_n\})\times (0,1)$. Take
as generators of $\pi_1(U-\Delta,b)$ $n$ loops $\gamma_1,\ldots,\gamma_n$
in $S^2\times \{1/2\}$ around the points $p_1,\ldots,p_n$ oriented so that
$\gamma_1\cdots\gamma_n=1$ in $\pi_1(U-\Delta,b)$. Let $T_1,
\ldots,T_n$ be the corresponding monodromy transformations in $SL
(H_1(X_b,\boldz))$. Using this notation, we have

\proclaim Theorem 1.11. Let $f:X\rightarrow B$ be a well-behaved
$T^3$-fibration, $b_0\in\Delta_d$, with $X_{b_0}$ a semi-stable fibre,
 $U,T_1,\ldots,T_n$ as above.
Then
\item{(1)} If $X_{b_0}$ is of type $(2,2)$, then $n=2$ and
in a suitable basis, $T_1=\pmatrix{1&0&1\cr 0&1&0\cr 0&0&1\cr}$
and $T_2=T_1^{-1}$.
\item{(2)} If $X_{b_0}$ is of type $(1,2)$ or $(2,1)$, then $n=3$
and in a suitable basis,
$T_1=\pmatrix{1&0&1\cr 0&1&0\cr 0&0&1\cr}$,
$T_2=\pmatrix{1&0&0\cr 0&1&1\cr 0&0&1\cr}$ and $T_3=\pmatrix{1&0&-1\cr
0&1&-1\cr 0&0&1\cr}$ in the $(1,2)$ case, and the transpose of these
matrices in the $(2,1)$ case. 
\item{(3)}
If $X_{b_0}$ is a fibre of type
$(1,1)$, then $n=4$. Furthermore, in a suitable
basis, possibly after relabelling the $p_i$'s,
$T_1=\pmatrix{1&1&0\cr 0&1&0\cr 0&0&1\cr}$,
$T_2=\pmatrix{1&0&0\cr 0&1&1\cr 0&0&1\cr}$,
$T_3=\pmatrix{1&-1&a\cr 0&1&0\cr 0&0&1\cr}$,
and $T_4=\pmatrix{1&0&-a-1\cr 0&1&-1\cr 0&0&1\cr}$.

Proof. Let $U\subseteq B$ be the open neighbourhood of $b_0\in\Delta_d$ as
above, with $f^{-1}(U-\{b_0\})\cong X'\times (0,1)$.
Then $H^i(X',\boldz)=H^i(f^{-1}(U-\{b_0\}),\boldz)$, 
which was computed in the proof of Theorem 1.3, with the result given
in (1.1).

We next study the Leray spectral sequence for $g:X'\rightarrow S^2$,
$H^i(S^2,R^jg_*{\bf Q})\Rightarrow H^n(X',{\bf Q})$. This was
shown in the proof of Theorem 1.3 to degenerate at the $E_2$-term
even with coefficients in $\boldz$. In addition, I claim the $E_2$
term looks like
$$\matrix{{\bf Q}&0&{\bf Q}\cr
{\bf Q}^{b_2}&0&{\bf Q}^{b_1}\cr
{\bf Q}^{b_1}&0&{\bf Q}^{b_2}\cr
{\bf Q}&0&{\bf Q}\cr}$$
First, note that $H^0(S^2,R^ig_*{\bf Q})=H^0(f^{-1}(U),R^if_*{\bf Q})=
H^i(X_b,{\bf Q})$, giving the first column. The third column then
follows from Poincar\'e-Verdier duality $\ext^0(R^ig_*{\bf Q},{\bf Q})
\cong \dual{H^2(S^2,R^ig_*{\bf Q})}$. From the knowledge of $H^i(X',{\bf Q})$ as
given in (1.1), and the degeneration of the $E_2$ term,
it follows that the middle column is 
indeed zero. 

Now let $\F=R^2g_*{\bf Q}$ in the $(1,2)$ or $(2,2)$
cases and $\F=R^1g_*{\bf Q}$ in
the $(2,1)$ case, so that $\F$ has ${\bf Q}^2$ as a subsheaf.
We in fact have an exact sequence
$$\exact{{\bf Q}^2}{\F}{{\bf Q}_V},$$
where $V=S^2-\{p_1,\ldots,p_n\}$, and
${\bf Q}_V$ is the sheaf ${\bf Q}$ on $V$ extended by zero outside of $V$.
Now $H^i(S^2,{\bf Q}_V)=H^i_c(V,{\bf Q})$, so we see that
$$\chi(S^2,\F)=\chi(S^2,{\bf Q}^2)+\chi(S^2,{\bf Q}_V)=4+(2-n)=6-n.$$
Now in the $(2,2)$ case, $\chi(S^2,\F)=b_1+b_2=4$, so $n=2$, and the description
of the monodromy is then clear. In the $(1,2)$ or $(2,1)$ case,
$\chi(S^2,\F)=b_1+b_2=3$, so $n=3$.
Since $T_1T_2T_3=I$, we see $T_1$ and $T_2$ must generate the monodromy
group, and it then follows from Theorem 1.10 that they must be of the form given
in the statement of the theorem in some basis.

Next consider the $(1,1)$ case. Then ${\bf Q}$ is a subsheaf of $R^2g_*{\bf Q}$,
and let $\F$ be the quotient $R^2g_*{\bf Q}/{\bf Q}$. Now there are two
types of points in $S^2\cap\Delta$. In the first, the monodromy takes the
form $\pmatrix{1&*&*\cr 0&1&0\cr 0&0&1\cr}$, and here $\F|_V$ on an open
neighbourhood $V$ of such a point $b$ is of the form ${\bf Q}_V\oplus 
{\bf Q}_{V-\{b\}}$. At the second type of point, the monodromy takes
the form $\pmatrix{1&0&*\cr 0&1&*\cr 0&0&1\cr}$, in which case locally
$\F$ has a subsheaf ${\bf Q}$ and the quotient $\F/{\bf Q}
={\bf Q}_{V-\{b\}}$. Now let $V_1$ be the complement in $S^2$ of 
points of the first type, and $V_2$ be the complement in $S^2$ of points
of the second type. Then globally on $S^2$ there is an exact sequence
$$\exact{{\bf Q}_{V_1}}{\F}{{\bf Q}_{V_2}},$$
where the inclusion ${\bf Q}_{V_1}$ in $\F$ is given by $e_2$.
From this we now obtain
$$\eqalign{2=\chi(S^2,R^2g_*{\bf Q})=\chi(S^2,{\bf Q})+\chi(S^2,\F)
&=\chi(S^2,{\bf Q})+\chi(S^2,{\bf Q}_{V_1})+\chi(S^2,{\bf Q}_{V_2})\cr
&=6-\#(S^2-V_1)-\#(S^2-V_2)\cr}$$
Thus $n=4$, as desired.

Finally, consider $T_1,\ldots,T_4$ 
with $T_1T_2T_3T_4=I$. At least one of these must
be of the first type and one of the second, and making a change of
basis and reordering the points if necessary, we can assume
$T_1=\pmatrix{1&1&0\cr 0&1&0\cr 0&0&1\cr}$ and $T_2=\pmatrix
{1&0&0\cr 0&1&1\cr 0&0&1\cr}$. Interchanging $p_3$ and $p_4$ if necessary,
$T_3$ and $T_4$ must be of the given form. $\bullet$
\bigskip

{\hd \S 2. Semistable Compactifications.}

We will now describe how to construct well-behaved topological 
compactifications of suitable $T^3$ fibre bundles $f_0:X_0\rightarrow B_0$.
This serves several purposes. First, it allows us to construct
examples of $T^3$-fibrations with monodromy of the sort
classified in \S 1. Second, it allows us to build up a set of
techniques for constructing $T^3$-fibrations. Third, and most
importantly, it gives compactifications of the dual fibration of 
a $T^3$-fibration.

We first describe the initial data required. In the rest of this
section, $B_0$ will be an open subset of a three-manifold $B$, 
$\Delta=B-B_0$, and
$f_0:X_0\rightarrow B_0$ a $T^3$ fibre bundle, which we will always take to
have a section. One reason for insisting on the existence of a section
is that a torus bundle with a section is completely determined by
its monodromy; in particular, if $\E$ is the total space of the vector bundle
corresponding to the sheaf $(R^{n-1}f_{0*}{\bf R})\otimes C^{\infty}(B_0)$,
$f_0$ is isomorphic to $\E/R^{n-1}f_{0*}\boldz$. (We normally just abuse
notation and identify this space with $R^{n-1}f_{0*}({\bf R}/\boldz)$.)
Of course, the local system $R^{n-1}f_{0*}\boldz$ is completely determined
by its monodromy.

We will require in addition the following
conditions on $\Delta$ and the monodromy of the fibre bundle:
\item{(1)} $\Delta$ decomposes as a disjoint union  $\Delta_d\cup \Delta_g$,
with $\Delta_d$ discrete and $\Delta_g$ a union of 1-manifolds which are
open intervals or circles. Each endpoint of a component in $\Delta_g$
is contained in $\Delta_d$. 
Thus $\Delta$ is a graph, with the points
of $\Delta_d$ the vertices and components of $\Delta_g$ the edges.
Furthermore, because of Theorem 1.11, we insist that the
valency of each vertex is 3 or 4.
\item{(2)} The monodromy group of $f_0:X_0\rightarrow B_0$
at each point in $\Delta$ is unipotent.
In addition, $R^3f_{0*}\boldz=\boldz$. This latter condition is equivalent to
the monodromy representation $\pi_1(B_0,b)\rightarrow GL(H_1(X_b,\boldz))$
taking values in $SL(H_1(X_b,\boldz))$.
\item{(3)} In suitable bases, the monodromy at each point
of $\Delta_g$ is as described in Theorem 1.7, and the monodromy at each point
of $\Delta_d$ is as described in Theorem 1.11, (2)-(4).

The main theorem of this section is

\proclaim Theorem 2.1. If $f_0:X_0\rightarrow B_0$ satisfies the above
properties, then there exists a topological manifold $X$, $X_0\subseteq X$,
and a well-behaved $T^3$-fibration $f:X\rightarrow B$ such that
the diagram
$$\matrix{X_0&\hookrightarrow&X\cr
\mapdown{f_0}&&\mapdown{f}\cr
B_0&\hookrightarrow&B\cr}$$
commutes.

The proof of this will be given after some basic constructions
are discussed. But first the main application:

\proclaim Corollary 2.2. If $f:X\rightarrow B$ is a well-behaved
$T^3$-fibration with only semistable fibres, then $f$ has a well-behaved
dual $\check f:\check X\rightarrow B$.

Proof. Conditions (1)-(3) hold for $f_0:X_0\rightarrow B_0$, by Definitions
1.1 and 1.2 and Theorems 1.7 and 1.11. 
Construct $\check f_0:\check X_0\rightarrow B_0$ by identifying
this with the family of tori $R^1f_{0*}\bf R/\boldz$. The monodromy matrices
of the local system $R^2\check f_{0*}\boldz \cong R^1f_{0*}\boldz$ are the
transpose inverses
of the monodromy matrices of the local system $R^2f_{0*}\boldz$.
Thus all monodromy groups are still unipotent, and conditions (1)-(3) hold
for $\check f_0:\check X_0\rightarrow B_0$. Then Theorem 2.1 implies there is
a well-behaved compactification, which is the desired dual.
$\bullet$
\bigskip

The basic idea for producing semi-stable fibrations is as follows.
In a small neighbourhood $U$ of a point $b_0\in B$ such that $X_{b_0}$ is
a semi-stable fibre, and $b\in U-\Delta$, there will always be a non-zero
subspace $L\subseteq H_1(X_b,\boldz)^G$ of monodromy invariant
cycles. One then expects the torus $T(L)=L\otimes_{\boldz} {\bf R}/L$
to act fibrewise by translation
on the fibration $f^{-1}(U)\rightarrow U$, giving
a lower dimensional fibration $f^{-1}(U)/T(L)\rightarrow U$. 
If one has a guess as to what this latter fibration is, one can try
to reconstruct $f^{-1}(U)\rightarrow U$. Now $\pi:f^{-1}(U)\rightarrow
f^{-1}(U)/T(L)$ is not in general a $T(L)$-bundle because
there will be fixed points of the $T(L)$ action in $f^{-1}(U)$.
However, there will be some dense open subset $V\subseteq f^{-1}(U)/T(L)$
for which $\pi^{-1}(V)\rightarrow V$ is a $T(L)$-bundle, with Chern class
$c_1\in H^2(V,L)$. Thus, if we want to construct examples
of semistable $T^n$-fibrations, we want to understand the following question:
given a manifold $\bar Y$, an open subset $Y\subseteq \bar Y$,
and a principal $T(L)$-bundle
$\pi:X\rightarrow Y$ with Chern class $c_1\in
H^2(Y,L)$, when does there exist a manifold $\bar X$
and a commutative diagram
$$\matrix{X&\hookrightarrow&\bar X\cr
\mapdown{\pi}&&\mapdown{\bar\pi}\cr
Y&\hookrightarrow&\bar Y\cr}$$
such that the $T(L)$ action on $X$ extends to a $T(L)$ action on
$\bar X$?

For us, the basic example of this sort of compactification is

{\it Example 2.3.} In [9], III.3.A, it was shown that the map
$f:{\bf C}^n\rightarrow {\bf R}^n$ given by 
$$f(z_1,\ldots,z_n)=(\im \prod z_i,|z_1|^2-|z_2|^2,\ldots,|z_1|^2-|z_n|^2)$$
was a special Lagrangian fibration. Now the map
$\bar\pi:{\bf C}^n\rightarrow {\bf C}\times {\bf R}^{n-1}$ given by
$$\bar\pi(z_1,\ldots,z_n)=(\prod z_i,|z_1|^2-|z_2|^2,\ldots,|z_1|^2-|z_n|^2)$$
is invariant under the $T^{n-1}$ action
$$(z_1,\ldots,z_n)\mapsto 
(e^{i\theta_1}z_1,\ldots,e^{i\theta_n}z_n),$$
$\theta_1+\cdots+\theta_n=0$. It is convenient to identify $T^{n-1}$
with $T(L)$, where $L=\{(a_1,\ldots,a_n)\in\boldz^n|\sum a_i=0\}$.
Then $\bar\pi$ is the quotient map ${\bf C}^n\rightarrow {\bf C}^n/T(L)$.
Furthermore, the set of points fixed by some subgroup of $T(L)$
is $Crit(f)=\bigcup\{z_i=z_j=0\}$. Thus $\bar\pi$ can be viewed as
a (partial) compactification of the $T(L)$-bundle
$\pi:
{\bf C}^n-Crit(f)\rightarrow {\bf C}\times{\bf R}^{n-1}-\bar\pi(Crit(f))$.

We will in fact only make use of the cases $n=2$ and $3$, but
more generally, much of what is done below can be extended to higher
dimensions. However, in dimensions 2 and 3, we have a complete
description of semistable monodromy, so we can produce a list of examples
of all possible monodromies. We build up these examples using the following
basic building blocks:

\proclaim Proposition 2.4.
(1) Let $\pi:X\rightarrow Y={\bf R}^3-\{0\}$ be
a principal $S^1$-bundle with Chern class $c_1\in H^2(Y,\boldz)=\boldz$,
and let $\bar Y={\bf R}^3\supseteq Y$. Then there is a unique topological
space $\bar X=X\cup \{p\}$ extending the topology on $X$
such that there is a commutative diagram
of continuous maps
$$\matrix{X&\hookrightarrow&\bar X\cr
\mapdown{\pi}&&\mapdown{\bar\pi}\cr
Y&\hookrightarrow&\bar Y\cr}$$
with $\bar\pi$ proper and $\bar\pi(p)=0$. Furthermore, $\bar X$ is a
topological manifold if and only if $c_1=\pm 1$.

Proof. A bundle over $Y$ with Chern class $c_1$
can be constructed as follows. Choose a hermitian
metric on $\O_{\Pone}(c_1)$, and let $X$ be the total
space of $\O_{\Pone}(c_1)$ minus the zero section. Then define
$\pi:X\rightarrow\Pone\times {\bf R}_{>0}\cong Y$ by taking
an element $s$ in the fibre of $\O_{\Pone}(c_1)$ at a point $x\in\Pone$
to $(x,||s||)\in Y$. The map $\pi$ defines an $S^1$-bundle with first
Chern class $c_1$.

It is easy to see that the topology on $\bar X=X\cup \{p\}$ is 
uniquely determined by the condition that $\bar\pi$ be a proper continuous
map. This topology has as basis the open sets in $X$ plus the inverse 
images of open sets in $\bar Y$ under $\bar\pi$. The space $\bar X$ can then be 
described as being obtained from $\O_{\Pone}(c_1)$ by contracting
the zero section. 
This is well-known to have
a quotient singularity if $c_1<-1$; if $c_1=-1$, then $\bar X$
is a manifold. If $c_1>0$, then $\pi:X\rightarrow Y$ is topologically
equivalent to the $S^1$-bundle with Chern class $-c_1$, but with the
orientation of the $S^1$ action reversed. If $c_1=0$, the fibration
is trivial and $\bar X$ is not a manifold. Thus $\bar X$ is a manifold
if and only if $c_1=\pm 1$. $\bullet$

We note that in the above Proposition, it is now clear that if 
$c_1=\pm 1$, then $\bar\pi:\bar X\rightarrow\bar Y$
coincides with $\bar\pi$ of Example 2.3 with $n=2$.

\proclaim Proposition 2.5. Let $Y$ be a manifold, $\pi:X\rightarrow Y$ a 
principal
$S^1$-bundle with $Y\subseteq\bar Y$ a manifold and $\bar Y-Y=S$ an
oriented connected submanifold of codimension three. Let $c_1\in H^2(Y,\boldz)$
be the Chern class of $\pi$. Now $H_S^3(\bar Y,\boldz)
=H^0(S,\boldz)
=\boldz$, and suppose the image of $c_1$ under the natural map
$H^2(Y,\boldz)\rightarrow H^3_S(\bar Y,\boldz)$ is $\pm 1$. Then there
exists a unique topological space $\bar X=X\cup S$ extending the topology
on $X$ such that there
is a commutative diagram of continuous maps
$$\matrix{X&\hookrightarrow&\bar X\cr
\mapdown{\pi}&&\mapdown{\bar\pi}\cr
Y&\hookrightarrow&\bar Y\cr}$$
with $\bar\pi$ proper and $\bar\pi|_S:S\rightarrow S$ the identity.
Furthermore, $\bar X$ is a topological manifold.

Proof. Again, the topology on $\bar X$ with basis consisting of
open sets in $X$ and inverse images of open sets in $\bar Y$ under
$\bar\pi$ is easily seen to be the unique topology in which $\bar\pi$
is proper and $\bar\pi|_S$ is the identity. Now if $p\in S$ is a point,
there exists an open neighbourhood $U\subseteq\bar Y$ of $p$ such
that $U\cong {\bf R}^{n-3}\times {\bf R}^3$ with $U\cap S$
identified with ${\bf R}^{n-3}\times\{0\}$. Then the restriction of
$\pi$ to $\{q\}\times ({\bf R}^3-\{0\})$ for any point
$q\in U\cap S$ is an $S^1$-bundle with Chern class $\pm 1$.
Thus by Proposition 2.4, we have a commutative diagram
$$\matrix{\pi^{-1}(U-S)&\hookrightarrow&{\bf R}^{n-3}\times{\bf C}^2\cr
\mapdown{\pi}&&\mapdown{\bar\pi'}\cr
U-S&\hookrightarrow&{\bf R}^{n-3}\times {\bf R}^3\cr}$$ where
$\bar\pi'$ is the identity on the first factor and the map $\bar\pi$ of
Example 2.3 on the second factor. Thus, by uniqueness of the compactification,
$\bar\pi'$ coincides with $\bar\pi^{-1}(U)\rightarrow U$, and thus in 
particular, $\bar X$ is a manifold. $\bullet$

{\it Example 2.6.} (1)  Let $D^*$ be a punctured disk,
$f:X_0\rightarrow D^*$ be a $T^2$-fibre bundle, with a section and
monodromy $\pmatrix{1&1\cr 0&1}$ in
some suitable basis $e_1,e_2$ of $H_1$ of a fibre. Then translation
on fibres by elements of
${\bf R}e_1$ induces an $S^1$ action on $X_0$, yielding an $S^1$-bundle
$\pi_0:X_0\rightarrow Y_0=S^1\times D^*$. One easy way to see what is the
Chern class of the bundle $\pi_0$ is via the Leray spectral sequence
for $\pi_0$: the image of $1\in H^0(Y_0,R^1\pi_{0*}\boldz)$ via the
differential $d:H^0(Y_0,R^1\pi_{0*}\boldz)\rightarrow H^2(Y_0,\boldz)$
is the Chern class. The group 
$H^2(X_0,\boldz)$ is easily seen to be torsion-free,
and thus $\coker d$ must be also. Since $H^2(Y_0,\boldz)\cong H^1(S^1,\boldz)
\otimes H^1(D^*,\boldz)\cong\boldz$, the Chern class $c_1$ of the bundle
$\pi_0$ must be $\pm 1$. 

To extend the fibration $f:X_0\rightarrow D^*$ to $\bar f:\bar X\rightarrow
D$, we proceed as follows. We have
$Y_0\subseteq \bar Y=S^1\times D$.
Let $Y=\bar Y-\{(p,0)\}$, where $p$ is a point in $S^1$, $0\in D$
the puncture point. Then $H^2(Y,\boldz)\cong H^3_{\{(p,0)\}}(\bar Y,\boldz)
=\boldz$ and the restriction map $H^2(Y,\boldz)\rightarrow H^2(Y_0,\boldz)$
is an isomorphism. Thus the $S^1$-bundle $\pi_0$ extends to an $S^1$-bundle
$\pi:X\rightarrow Y$. Then Proposition 2.5 
can be applied to obtain a manifold $\bar X=X\cup
\{pt\}$ and a map
$\bar\pi:\bar X\rightarrow\bar Y$. When $\bar\pi$ is 
composed with the projection
$\bar Y\rightarrow D$, we get a map $\bar f:\bar X
\rightarrow D$. This map 
is proper, and $f^{-1}(0)$ is an $I_1$ fibre.

It is worth noting that the effect of changing the sign of $c_1$ in
this construction is simply to reverse the direction of the
$S^1$ action on $X$. This is equivalent to replacing $e_1$ with
$-e_1$, and hence the monodromy changes from $\pmatrix{1&1\cr 0&1\cr}$
to $\pmatrix{1&-1\cr 0&1\cr}$. Clearly though it gives the same
$T^2$-fibration.

We also note that it is easy to see the fibration $\bar f:\bar X\rightarrow
D$ is well-behaved and has a section.

(2) This example will not be used, but demonstrates that the topology
of the singular fibre may not be reflected in the monodromy about it.
Let $Y=S^1\times D-\{(p_1,0),(p_2,0)\}$, $\bar Y=S^1\times D$,
 $p_1$ and $p_2$ two distinct
points in $S^1$, $D$ a two-dimensional disk, $0\in D$. 
Then $H^2(Y,\boldz)=H^3_{\{(p_1,0)\}}(\bar Y,\boldz)
\oplus H^3_{\{(p_2,0)\}}(\bar Y,\boldz)$, and the restriction map
$H^2(Y,\boldz)\rightarrow H^2(S^1\times D^*,\boldz)$ is given by
summation. Construct the $S^1$-bundle $\pi:X\rightarrow Y$
with Chern class $(1,-1)\in H^2(Y,\boldz)$. Then $\pi$ restricted to
$S^1\times D^*$ is trivial. Also, $\pi$ can be compactified
as above to obtain a fibration $\bar f:\bar X\rightarrow D$,
with central fibre being a union of two $S^2$'s. However, the
monodromy about the central fibre is trivial. This of course could
not happen
for a holomorphic elliptic fibration.

(3) {\it Type $(2,2)$ fibres}.
Let $B=D\times (0,1)$, $B_0=D^*\times (0,1)$, $\bar Y=T^2\times B$.
Let $S=S^1\times \{0\}\times
(0,1)\subseteq T^2\times \{0\}\times (0,1)\subseteq \bar Y$
be a surface, where $S^1\subseteq T^2$ is a circle whose homology
class is a primitive
element $ae_2+be_3\in H_1(T^2,\boldz)$, and $e_2,e_3$ is a basis
for $H_1(T^2,\boldz)$. Let $Y=\bar Y-S$. There is an exact sequence
$$\exact{H^2(\bar Y,\boldz)}{H^2(Y,\boldz)}{H^3_{S}(\bar Y,\boldz)}.$$
$H^2(\bar Y,\boldz)$ is isomorphic to $H^2(T^2,\boldz)$,
and the exact sequence is split via the restriction map $H^2(Y,\boldz)
\rightarrow H^2(T^2\times\{p\},\boldz)$ for a point $p\in B_0$.
Thus $H^2(Y,\boldz)=H^2(\bar Y,\boldz)\oplus H^3_S(\bar Y,\boldz)$. Also
$H^3_S(\bar Y,\boldz)=H^0(S,\boldz)=\boldz$.

Now let $\pi:X\rightarrow Y$ be an $S^1$-bundle with first Chern class $c_1
\in H^2(Y,\boldz)$ given by either $1$ or $-1$ in $H^3_S(\bar Y,\boldz)
=\boldz$. (Note this equality depends on
choices of orientation for $\bar Y$ and $S$, so we won't
worry about the sign.) Then applying Proposition 2.5, 
we can construct a compactification
$\bar\pi:\bar X\rightarrow\bar Y$. By composing this with
the projection $\bar Y\rightarrow B$, we obtain a map
$\bar f:\bar X\rightarrow B$ which is a $T^3$-bundle over $B_0$, and has
fibres over $\Delta=B-B_0$ homeomorphic to $I_1\times S^1$.
Furthermore, it is clear from the construction that we can write
$\bar X=\bar X'\times S^1\times (0,1)$, where $\bar X'\rightarrow
D$ is a fibration of the type constructed in (3). In particular,
choosing a basis for $H_1(\bar X_b,\boldz)$, $b\in B_0$ of
the form $e_1,e_2,e_3$, where $e_1$ is a fibre
of $\pi$ and $e_2,e_3$ project to $e_2,e_3\in H_1(\bar Y_{b},\boldz)$
under the map $\pi$, we see that $e_1$ and $ae_2+be_3$ are
a basis for the monodromy invariant cycles. Furthermore, the
monodromy about a loop around $\Delta$ is $\pmatrix{1&\pm b&\mp a\cr
0&1&0\cr 0&0&1\cr}$, where the choice of sign can
be changed by replacing $e_1$ with $-e_1$. In particular, this provides
a compactification of a $T^3$-fibre bundle $f:X_0=\bar f^{-1}(B_0)
\rightarrow B_0$
with this monodromy.

(4) {\it Type $(2,1)$ fibres.}
Let $B_0=(S^2-\{p_1,p_2,p_3\})\times (0,1)$, and let $B=B^3$, 
the three-ball. Identify $B-\{0\}$ with $S^2\times (0,1)$ and embed
$B_0\subseteq B$ via $S^2-\{p_1,p_2,p_3\}\subseteq S^2$. Then
$\Delta=B-B_0$ is the cone over three points (best-visualized as a ``Y''). 
Write $\Delta=\Delta_1\cup \Delta_2\cup\Delta_3\cup\{0\}$, where
$\Delta_i=\{p_i\}\times (0,1)$ are the ``legs'' of the ``Y''.
Let $\bar Y=T^2\times B$, and let $e_2,e_3$ be a basis for
$H_1(T^2,\boldz)$. Let $S\subseteq
\bar Y$ be a ``pair of pants'' sitting over $\Delta$ constructed
as follows: $S\cap (T^2\times\Delta_i)$ will be a cylinder 
$S^1\times \Delta_i$, with the $S^1$ of class $-e_3$, $e_2$ and $-e_2+e_3$
for $i=1,2,3$ respectively. These three cylinders are then joined in
$T^2\times\{0\}$ to form the pair of pants:
$$\epsfbox{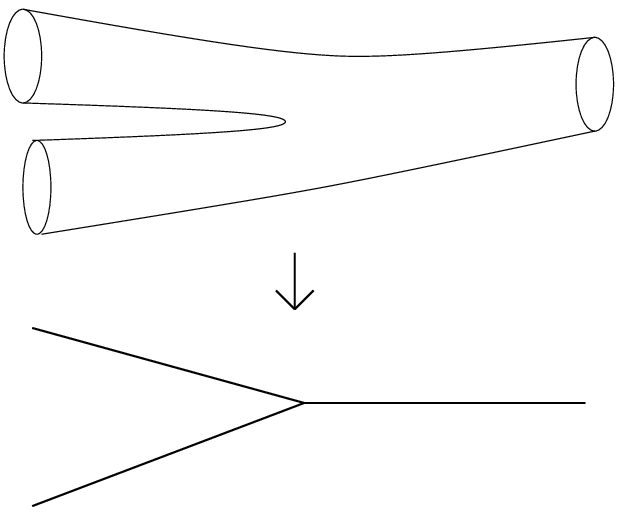}\leqno{\hbox{Figure $2.1$}}$$
Note that there is some ambiguity as to the fibre of $S\rightarrow \Delta$
over $0\in \Delta$. This will affect the precise shape of the singular
fibre, but not the topology of the total space. 
This could be a figure eight or else my preferred choice,
suggested by Ruan
(illustrated in [17], Figure 6), and shown in the following figure:
$$\epsfbox{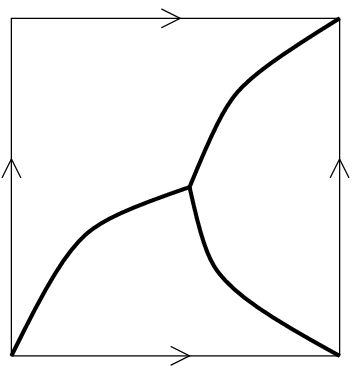}\leqno{\hbox{Figure $2.2$}}$$
Here the torus has been depicted as a rectangle with opposite sides identified,
and the darker lines represent the fibre of $S$ over $0\in\Delta$.
We will not limit our choices at this point, and instead allow
any choice, as we are currently just interested in the topology of
the fibration. This choice might become clearer if a local model for
a special Lagrangian fibration with this sort of singular fibre is found.

In any event, put $Y=\bar Y-S$, and then $H^2(Y,\boldz)=
H^2(\bar Y,\boldz)\oplus H^3_S(Y,\boldz)=\boldz\oplus\boldz$ as in
(3). Take
$c_1=1\in H^3_S(Y,\boldz)$, and we obtain an $S^1$-bundle $X\rightarrow Y$,
which by Proposition 2.5 compactifies to $\bar X\rightarrow \bar Y$. By
composing with the projection $\bar Y\rightarrow B$, we obtain
a proper map $\bar f:\bar X\rightarrow B$. Now if $\gamma_1,\gamma_2,\gamma_3$
are loops in $B_0$ around $\Delta_1,\Delta_2$ and $\Delta_3$, with
$\gamma_1\gamma_2\gamma_3=1$ in $\pi_1(B_0,b)$, then by (3)
it follows that the monodromy matrices $T_i$ about $\gamma_i$ take
the form (in the basis $e_1,e_2,e_3$ of $H_1(\bar X_b,\boldz)$ as in (3))
$$\hbox{$T_1=\pmatrix{1&\pm 1&0\cr 0&1&0\cr 0&0&1\cr},T_2=\pmatrix{1&0&\pm 
1\cr 0&1&0\cr 
0&0&1\cr},T_3=\pmatrix{1&1&1\cr 0&1&0\cr 0&0&1\cr}$ or 
$\pmatrix{1&-1&-1\cr 0&1&0\cr 0&0&1\cr}$}$$
Using the additional requirement that $T_1T_2T_3=I$, and possibly changing
the sign of $e_1$, we see we must have
$$T_1=\pmatrix{1&1&0\cr 0&1&0\cr 0&0&1\cr},T_2=\pmatrix{1&0&1\cr 0&1&0\cr 
0&0&1\cr},T_3=\pmatrix{1&-1&-1\cr 0&1&0\cr 0&0&1\cr}$$
Thus $\bar f:\bar X\rightarrow B$ provides a compactification of
a fibration $f_0:X_0\rightarrow B_0$ with the above monodromy.

The fibre over $0\in B$, $\bar X_0$, is described as follows. Let $S_0
\subseteq T^2$ be the fibre of $S\rightarrow\Delta$ over $0\in B$. Then
$\bar X_0=S^1\times T^2/\sim$, where $(x,y)\sim (x',y')$ if
$(x,y)=(x',y')$ or $y=y'\in S_0$. Thus $S^1\times S_0$ is contracted
to $S_0$.

(5) {\it Type $(1,1)$ fibres.}
Let $B_0=(S^2-\{p_1,\ldots,p_4\})\times (0,1)$, 
and let $B=B^3$, the three-ball, with
$\Delta=B-B_0$ the cone over four points. Let $\Delta_i=\{p_i\}\times (0,1)$,
so that $\Delta=\Delta_1\cup\cdots\cup\Delta_4\cup\{0\}$.
Let $l_1$ be the
cone over $\{p_1,p_3\}$
and $l_2$ the cone over $\{p_2,p_4\}$,
so that $\Delta=l_1\cup l_2$ and
$l_1\cap l_2=\{0\}\subseteq B$. 
Let $X'\rightarrow D$ be a $T^2$-fibration as constructed in (1)
with monodromy $\pmatrix{1&1\cr 0&1\cr}$, and identifying
$D\times (0,1)$ with $B$ so that $\{0\}\times (0,1)$ is
identified with $l_2$, we obtain a $T^2$-fibration $g:\bar Y=X'\times (0,1)
\rightarrow B$
with discriminant locus $l_2$ and monodromy $\pmatrix{1&1\cr 0&1\cr}$
about $l_2$ in a suitable basis $e_2,e_3$ of $H_1(\bar Y_b,\boldz)$,
for some $b\in B_0$.

Now construct a surface $S\subseteq \bar Y$ as follows.
$S$ will be a cylinder fibred in circles over $l_1$. To specify the cohomology
class of this $S^1$ over $\Delta_1$ and $\Delta_3$, we have to be careful
about how we identify these classes. In particular, as
$e_3$ is not invariant under monodromy,
it is not well-defined globally. We picture here the sphere $S^2-
\{p_1,p_2,p_3,p_4\}$, along with loops $\gamma_1,\ldots,\gamma_4$
and paths $c_1$ and $c_3$ joining $p_1$ and $p_3$ to $b$:
$$\epsfbox{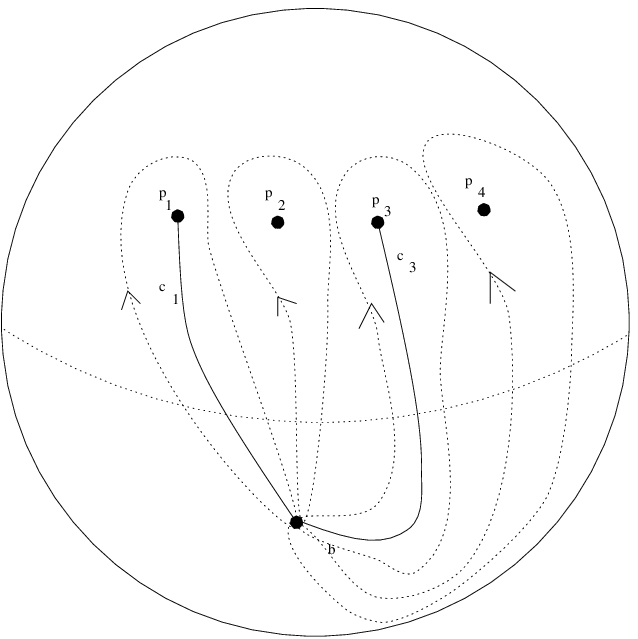}\leqno{\hbox{Figure $2.3$}}$$
We can then identify $H_1(\bar Y_{p_1},\boldz)$ and $H_1(\bar Y_{p_3},
\boldz)$ with $H_1(\bar Y_b,\boldz)$ via the paths $c_1$ and $c_3$.
Using these identifications,
$S\cap g^{-1}(\Delta_1)
=S^1\times\Delta_1$, with $S^1\subseteq T^2$ of class $e_3$, while 
$S\cap g^{-1}(\Delta_3)=S^1\times \Delta_3$ with $S^1\subseteq T^2$ of class
$ae_2+e_3$. These two circles specialise to circles of the same class
in $\bar Y_0$, as $e_2$ is the vanishing cycle.
Thus these two cylinders can be joined up.

Now let $Y=\bar Y-S$. We have as usual a splitting $H^2(Y,\boldz)
=H^2(\bar Y,\boldz)\oplus H^3_S(\bar Y,\boldz)$, and take an $S^1$-bundle
$\pi:X\rightarrow Y$ given by $c_1=1\in H^3_S(\bar Y,\boldz)$. 
As usual, we compactify to obtain $\bar\pi:\bar X\rightarrow\bar Y$,
with the composed map $\bar f=g\circ\bar\pi:\bar X\rightarrow B$. Then
$X_0=\bar f^{-1}(B_0)\rightarrow B_0$ has monodromy as given in Theorem 1.11,
(3), and we have constructed a compactification of such a $T^3$-bundle.

The fibre of type $(1,1)$ constructed here can then be described as
$S^1\times I_1/\sim$, where $(x,y)\sim(x',y')$ if $(x,y)=(x',y')$ or
$y=y'\in S\cap g^{-1}(0)$. The singular locus of this fibre is a circle.
\bigskip

We still need a construction for type $(1,2)$ fibres. An example
of such a fibration with a type $(1,2)$ fibre has been given
in [6], \S 1. However, I would like to give a different construction
of that example, more in keeping with the previous construction techniques
given here. To do so, we need to compactify $T^2$ bundles.

Let $L\cong\boldz^2$ and let
$T(L):=L\otimes_{\boldz}{\bf R}/L$. The Chern class of a principal
$T(L)$-bundle over $Y$ is an element of $H^2(Y,L)$.  Indeed, if ${\cal R}$
denotes the sheaf of continuous real-valued functions on $Y$, then
a principal $T(L)$-bundle is determined by an element of $H^1(Y,
L\otimes_{\boldz}{\cal R}/L)$, which is isomorphic to $H^2(Y,L)$.
The corresponding element in this latter group is the Chern class. It 
is also convenient to think of a $T^2$-bundle as a product of $S^1$-bundles:
if $e_1,e_2$ is a basis for $L$, giving a decomposition $H^2(Y,L)
\cong H^2(Y,\boldz)e_1\oplus H^2(Y,\boldz)e_2$, we can write
$c\in H^2(Y,L)$ as $c=(a_1,a_2)$, $a_i\in H^2(Y,\boldz)$. If $X_i\rightarrow Y$
is an $S^1$-bundle with Chern class $a_i$, then $X=X_1\times_Y X_2\rightarrow
Y$ is a $T(L)$-bundle with Chern class $c$. Note that $X_1=X/T(\boldz e_2)$,
$X_2=X/T(\boldz e_1)$.
Here, $T(\boldz e_i)$ is the subtorus of $T(L)$ generated by $e_i$.

\proclaim Proposition 2.7. Let $Y$ be a manifold, $\pi:X\rightarrow Y$
a principal $T(L)$-bundle with $Y\subseteq\bar Y$ a manifold and $\bar Y-Y=S$
an oriented connected submanifold of codimension three. Let $c_1
\in H^2(Y,L)$ be the Chern class of $\pi$, and suppose the image of
$c_1$ under the natural map $H^2(Y,L)\rightarrow H^3_S(\bar Y,L)=H^0(S,L)
=L$ is a primitive class $e_1\in L$. Then there exists a unique topological
space $\bar X=X\cup S'$ such that there is a commutative diagram of continuous
maps
$$\matrix{X&\hookrightarrow&\bar X\cr
\mapdown{\pi}&&\mapdown{\bar\pi}\cr
Y&\hookrightarrow&\bar Y\cr}$$
with $\bar\pi$ proper, the $T(L)$ action on $X$ extending to
a $T(L)$ action on $\bar X$ so that $\bar\pi|_{S'}:S'\rightarrow S$
is a $T(L/\langle e_1\rangle)$-bundle. Furthermore, $\bar X$ is
a manifold.

Proof. Choose $e_2$ so that $e_1,e_2$ forms a basis for $L$.
Then $c_1\in H^2(Y,L)=H^2(Y,\boldz)e_1\oplus H^2(Y,\boldz)e_2$
can be written as $c_1=(a_1,a_2)$, with the image of $a_1$
under the natural map $H^2(Y,\boldz)\rightarrow H^3_S(\bar Y,\boldz)$
being 1 and the image of $a_2$ under this same map being zero. Thus 
$a_2$ lifts to an element $\bar a_2\in H^2(\bar Y,\boldz)$.
Note this lifting is unique because $H^2_S(\bar Y,\boldz)=0$.
Let $\pi_i:X_i\rightarrow Y$ be the $S^1$-bundle with Chern class $a_i$.
Then by Proposition 2.5, $\pi_1$ can be extended to $\bar\pi_1:
\bar X_1\rightarrow\bar Y$, and we also have an $S^1$-bundle
$\bar\pi_2:\bar X_2\rightarrow\bar Y$ with Chern class $\bar a_2$,
extending $\pi_2$. Then we can take $\bar\pi:\bar X=
\bar X_1\times_{\bar Y}\times \bar X_2\rightarrow\bar Y$ as the
desired compactification of $\pi$.

To show uniqueness, suppose we are given $\bar\pi:\bar X\rightarrow\bar Y$
with the desired properties. Letting $\bar X_1=\bar X/T(\boldz e_2)$
and $\bar X_2=\bar X/T(\boldz e_1)$, $\bar\pi_i:\bar X_i\rightarrow
\bar Y$ the projection, it is clear that $\bar\pi_1:\bar X_1\rightarrow
\bar Y$ restricts to the $S^1$-bundle $\pi_1:X_1\rightarrow Y$,
and $\bar\pi_1^{-1}(S)\cong S$, so $\bar X_1$ is uniquely determined
by Proposition 2.5. On the other hand, $\bar\pi_2:\bar X_2\rightarrow\bar Y$
is an $S^1$-bundle, necessarily with Chern class $\bar a_2$. Finally,
$\bar X=\bar X_1\times_{\bar Y}\bar X_2$, so $\bar X$ is necessarily
as constructed as above. $\bullet$

{\it Example 2.8. Type $(2,2)$ fibres again.} One can apply Proposition 2.7 
to give an alternate construction of type $(2,2)$ fibres.
Take $Y=S^1\times D\times (0,1)-\{p\}\times\{0\}\times (0,1)
\subseteq \bar Y=S^1\times D\times (0,1)$. Then $H^2(Y,L)=H^3_{\{p\}\times
\{0\}\times (0,1)}(\bar Y,L)=L$, and choosing $c_1\in L$ 
primitive, one obtains a $T(L)$-bundle over $Y$ which can be 
compactified to give a map $\bar\pi:\bar X\rightarrow \bar Y$.
Composing with the projection $\bar Y\rightarrow B=D\times (0,1)$
gives a $T^3$-fibration $\bar f:\bar X\rightarrow B$. 
Let $e_1,e_2$ be a basis for $L$, and take 
for $H_1(\bar X_b,\boldz)$ a basis $e_1,e_2,e_3$, with $e_3$ mapping
under $\bar\pi$ to a generator of $H_1(S^1,\boldz)$. Then it is easy
to see that if $c_1=ae_1+be_2$,
then the monodromy about a loop around $\{0\}\times (0,1)$ is
$$\hbox{$\pmatrix{1&0&a\cr 0&1&b\cr 0&0&1\cr}$
or $\pmatrix{1&0&-a\cr 0&1&-b\cr 0&0&1\cr}$.}$$
Either can occur by changing the sign of the $T(L)$ action.

The final variant of this idea we need is

\proclaim Proposition 2.9. Let $Y$ be a four-manifold, $\pi:X\rightarrow Y$
a principal 
$T(L$)-bundle with $Y\subseteq\bar Y$ a manifold and $\bar Y-Y=\Delta$
a trivalent graph, i.e. a union of disjoint curves and vertices with each
vertex being in the closure of precisely three curves. We do not
assume $\Delta$ is compact. Write $\Delta=\Delta_g\cup\Delta_d$,
with $\Delta_g$ the union of (open) edges and $\Delta_d$ the union of
vertices. Let $c_1\in H^2(Y,L)$ be the Chern class of $\pi$, and denote
also by $c_1$ its image under the natural map $H^2(Y,L)\rightarrow
H^3_{\Delta}(\bar Y,L)$. Then $H^3_{\Delta}(\bar Y,L)
=H_1^{BM}(\Delta, L)$, the Borel-Moore homology of $\Delta$ with coefficients
in $L$. If $\bar\Delta
=\Delta\cup\{\infty\}$ is the one-point compactification of $\Delta$,
then $H_1^{BM}(\Delta,L)\cong H_1(\bar\Delta,\{\infty\},L)$,
the simplicial relative homology group. Then if $c_1\in H^3_{\Delta}(\bar Y,L)$
is represented by a relative simplicial 1-cycle on $\bar\Delta$, assigning
to each oriented edge $\Delta_i$ of $\bar\Delta$ the element $a_i\in L$;
$a_i$ is primitive for each $i$; and $a_i,a_j,a_k$ span $L$ whenever
$\Delta_i,\Delta_j,\Delta_k$ have a common vertex, then there exists
a unique topological space $\bar X=X\cup S$ such that
there exists a commutative diagram of continuous maps
$$\matrix{X&\hookrightarrow&\bar X\cr
\mapdown{\pi}&&\mapdown{\bar\pi}\cr
Y&\hookrightarrow&\bar Y\cr}$$
so that $\bar\pi$ is proper, $T(L)$ acts on $\bar X$, $\bar\pi^{-1}(
\Delta_i\cap\Delta_g)\rightarrow\Delta_i\cap\Delta_g$ is a
$T(L/\langle a_i\rangle)$-bundle for each edge $\Delta_i$, and
$\bar\pi^{-1}(\Delta_d)\rightarrow\Delta_d$ is an isomorphism. Furthermore,
$\bar X$ is a manifold.

Proof. That $H^3_{\Delta}(\bar Y,L)=H_1^{BM}(\Delta,L)$ follows from a suitable
form of Poincar\'e duality (see [3], V, Theorem 9.3). The description
of $H_1^{BM}(\Delta,L)$ follows from [3], V. Now using Proposition 2.7,
$\pi:Y\rightarrow X$ can be extended to $\pi':X'\rightarrow Y'=\bar Y-\Delta_d$
in a unique way to maintain the desired properties. Putting $\bar X
=X'\cup \Delta_d$, there is a unique topology on $\bar X$ inducing
the topology on $X'$ and making the extension of $\pi'$ to
$\bar\pi:\bar X\rightarrow\bar Y$ proper. We only need to check $\bar X$ is
a manifold in a neighbourhood of a point $p\in\Delta_d$.
To do so, we use the
local model in Example 2.3 for such a fibration. First note that we can take
$L\subseteq\boldz^3$ with standard basis $e_1,e_2,e_3$ and
$L=\{\sum b_ie_i|\sum b_i=0\}$. Let $a_i$ be the coefficient of the oriented
edge $\Delta_i$ in the 1-cycle representing $c_1$.
If edges, say, $\Delta_1,\Delta_2,\Delta_3$ meet
at $p\in\Delta_d$, then possibly after changing orientations of these
three edges, we must have $a_1+a_2+a_3=0$. This is the cycle
condition. Since we also assume the $a_i$'s are primitive and $a_1,a_2,a_3$
span $L$,
without loss of generality we can take $a_1=e_2-e_3,
a_2=e_3-e_1,$ and $a_3=e_1-e_2$. 
Now in the $T(L) $-action on ${\bf C}^3$ given in Example 2.3, $\{z_i=z_j=0\}$
is fixed by $T(a_k\boldz)$, if $\{i,j,k\}=\{1,2,3\}$. This makes
it clear that $\pi_{2.3}:{\bf C}^3-\{0\}\rightarrow {\bf C}\times {\bf R}^2
-\{0\}$ given in Example 2.3 coincides with
$\bar\pi^{-1}(U-\{p\})\rightarrow U-\{p\}$, where $U\cong {\bf C}\times {\bf R}^2$
is a small open neighbourhood of $p$. Thus $\bar\pi_{2.3}$ coincides with
$\bar\pi:\bar\pi^{-1}(U)\rightarrow U$, and
$\bar X$ is a manifold. $\bullet$

{\it Example 2.10. Type $(1,2)$ fibres.} 
Let $B$ and $B_0$ be as in Example 2.6, (4).
Take $Y=S^1\times B-\{p\}\times\Delta\subseteq \bar Y=S^1\times B$,
where $p\in S^1$ is a point. Then $H^2(Y,L)=H^3_{\{p\}\times\Delta}(\bar Y,
L)\cong H_1(\bar\Delta,{\infty},L)$, and we can take $((1,0),(0,1),(-1,-1))\in H^2(Y,L)$ as the
Chern class of a $T^2$-bundle $\pi:X\rightarrow Y$, which we
then compactify using Proposition 2.9  to $\bar\pi:\bar X\rightarrow\bar Y$.
It is then easy to see that the composed fibration $\bar f:\bar X\rightarrow
B$ has monodromy  dual to that in Example 2.6, (4), 
and has a section. We note this construction 
coincides topologically with the example given in [6], \S 1.

{\it Proof of Theorem 2.1.} We simply piece together the examples above.
First suppose we are given $T^3$-fibrations 
$f_U:X_U\rightarrow U$, $f_V:X_V\rightarrow V$, $U,V\subseteq B$,
such that $f_U$ and $f_V$ have sections, and $f_U:f_U^{-1}(U\cap V)\rightarrow
U\cap V$, $f_V:f_V^{-1}(U\cap V)\rightarrow U\cap V$ are $T^3$-fibre bundles
yielding the same monodromy representation $\pi_1(U\cap V,b)\rightarrow
SL_3(\boldz)$, $b\in U\cap V$ a basepoint. Then as observed earlier,
because $f_U$ and $f_V$ have sections, there is an isomorphism $\varphi$
in a commutative diagram 
$$\matrix{f_U^{-1}(U\cap V)&\mapright{\varphi}&f_V^{-1}(U\cap V)\cr
\mapdown{f_U}&&\mapdown{f_V}\cr
U\cap V&\mapright{=}& U\cap V\cr}$$
which takes a given section of $f_U$ to a given section of $f_V$. Thus
we can glue $X_U$ and $X_V$ via $\varphi$ to obtain a $T^3$-fibration
$f:X\rightarrow U\cup V$, with a section.

This gluing construction extends easily to the case that
$f_U:f_U^{-1}(U\cap V)\rightarrow U\cap V$ and
$f_V:f_V^{-1}(U\cap V)\rightarrow U\cap V$ are well-behaved 
$T^3$-fibrations whose discriminant loci coincide, and
only have singular fibres of type $(2,2)$, assuming that locally
the fibrations are of the form given in Example 2.6, (3). In this case, as long
as the monodromy representations for $f_U$ and $f_V$ coincide, one
can still find an isomorphism $\phi:f_U^{-1}(U\cap V)
\rightarrow f_V^{-1}(U\cap V)$ as above.

Now given $f_0:X_0\rightarrow B_0$ satisfying the hypotheses of the
theorem, choose an open covering $\{U_i\}$ of $\Delta$, $U_i\subseteq B$,
such that $X_{0,i}=f_0^{-1}(U_i\cap B_0)\rightarrow U_i\cap B_0$ is
of the type given in Example 2.6, (3)-(5) or Example 2.10. In each case this is
compactified to $X_i\rightarrow U_i$ with a choice of section, 
and these compactifications can then be glued
into $f_0:X_0\rightarrow B_0$, matching up the chosen
sections, to give
the desired fibration. It is easy to see this fibration is 
well-behaved. $\bullet$

{\it Remark 2.11.} Theorem 2.1 has produced a topological manifold $X$,
but it will be important to know that this manifold has a smooth
structure if $B$ and $\Delta_g\subset B$ are smooth. 
It will not matter what the smooth structure is, because we will 
be interested in situations where [20], Theorem 6 applies, which implies
that given suitable hypotheses, all smooth structures on $X$ are diffeomorphic.

The construction of Proposition 2.5
produces a smooth manifold $\bar X$ as long as
$\bar Y$ and $S\subseteq\bar Y$ are smooth. Thus
the only places in the constructions of Example 2.6 where we leave the
smooth category are (4) and (5). In (4), the $S$ we constructed
in $\bar Y$ is not smooth. However, one can find a homeomorphism
$\varphi$ of $\bar Y$ which is the identity outside of
some neighbourhood of $T^2\times\{0\}$, such that $\varphi^{-1}(S)$ is smooth.
Then we can replace the projection $p:\bar Y\rightarrow B$ with
$p'=p\circ\varphi$ and replace $S\subseteq\bar Y$ with $S'=
\varphi^{-1}(S)\subseteq
\bar Y$. Using $S'$, we construct $\bar\pi:\bar X\rightarrow \bar Y$
as before. Then $\bar X$ does have a smooth
structure. However, the map $\bar f=p'\circ\bar\pi:\bar X\rightarrow B$ may
not be differentiable, as $\varphi$ will not be. It is not clear if we expect
$\bar f$ to be differentiable at such a singular fibre. However,
we can guarantee at least that $\bar f$ is differentiable on $X^{\#}$.

A similar argument produces a smooth structure in (5) on $X$.
\bigskip

We conclude this section with a few basic results about the topology
of the $T^3$-fibrations constructed in this section.

\proclaim Theorem 2.12. Let $f_0:X_0\rightarrow B_0$ be as in Theorem 2.1,
and $f:X\rightarrow B$ the compactification produced in the proof of
Theorem 2.1. If $B$ is simply connected and $H^0(B,R^1f_*\boldz/n\boldz)
=0$ for all $n$, then $X$ is simply connected.

Proof. Let $\mu:\tilde X\rightarrow X$ be the universal cover of $X$, 
and define a space $\tilde B=\tilde X/\sim$, where $x\sim y$ if
$f(\mu(x))=f(\mu(y))=b$ and $x$ and $y$ are in the same connected component
of $(f\circ\mu)^{-1}(b)$. Then $f\circ\mu$ factors as $\tilde X
\mapright{\tilde f} \tilde B\mapright{\gamma} B$. 

{\it Claim.} $\gamma:
\tilde B\rightarrow B$ is a covering, and thus since $B$ is simply
connected, $\gamma$ is an isomorphism. 

Proof.
First for any point $b\in B$, there is a sufficiently small
open neighbourhood $U$ of $b$ such that $V=f^{-1}(U)\rightarrow U$ looks like
one of the cases of Example 2.6, (3)-(5) or Example 2.10,
or else $V=T^3\times U$. Let $\tilde V
\rightarrow V$ be the universal covering. I claim $\tilde V\rightarrow U$
has connected fibres. Indeed, if $X_b$ is a non-singular fibre, then
$V=T^3\times U$ and $\tilde V={\bf R}^3\times U$. If $X_b$ is
a fibre of type $(2,2)$ or $(2,1)$, then by construction there is
an $S^1$ action on $V$ such that $V/S^1=T^2\times U$, and $\tilde
V=V\times_{T^2\times U} {\bf R}^2\times U$. If $X_b$ is of type
$(1,2)$, $\pi_1(V)=\boldz$, there is a local $T^2$-action on
$V$, $V/T^2=S^1\times U$, and $\tilde V=V\times_{S^1\times U}
{\bf R}\times U$. Finally, in the $(1,1)$ case, $\pi_1(V)=\boldz$
again, and if, as in Example 2.6, (5), 
$\pi:V\rightarrow \bar Y$ is the quotient by
the $S^1$-action and $\tilde Y\rightarrow\bar Y$ is the universal
cover of $\bar Y$, then $\tilde V=V\times_{\bar Y}\tilde Y$.
In all these cases, one sees
$\tilde V\rightarrow B$ has connected fibres. Now we can conclude that
$\gamma:\tilde B\rightarrow B$ is a covering: 
$(\mu\circ f)^{-1}(U)$
must consist of a number of connected components $V_1,\ldots,V_n$
which are quotients of the universal covering $\tilde V$, and hence the
maps $V_i\rightarrow U$ have connected fibres. We conclude that $\gamma^{-1}(U)$
consists of $n$ disjoint copies of $U$, so $\gamma$ is a covering space,
hence an isomorphism, as $B$ is simply connected. $\bullet$

It is now clear that $\pi_1(X)$ is simply the Galois group of the
covering $\tilde X_b\rightarrow X_b$ for $X_b$ a non-singular fibre. This
is an abelian group, hence if non-zero, one would have 
$H^1(X,\boldz/n\boldz)\not=0$ for some $n$. But by the Leray spectral
sequence, this is not possible if $H^0(B,R^1f_*\boldz/n\boldz)=0$ for all
$n$. $\bullet$

For the remainder of the section, we assume we have put
a differentiable structure on $X$ as in Remark 2.5.

\proclaim Proposition 2.13. If $f:X\rightarrow B$ is as in Theorem 2.1 and
$B\cong S^3$, then $w_2(X)$, the second Steifel-Whitney class of $X$,
vanishes.

Proof. Since $Crit(f)$ is a union of two-dimensional manifolds,
$$H^2(X,\boldz/2\boldz)\cong H^2(X^{\#},\boldz/2\boldz),$$ and thus
$w_2(X)\in H^2(X,\boldz/2\boldz)$ coincides with $w_2(X^{\#})
\in H^2(X^{\#},\boldz/2\boldz)$. However the tangent bundle of 
$X^{\#}$ fits into an exact sequence
$$\exact{\T_{X^{\#}/B}}{\T_{X^{\#}}}{f^{\# *}\T_{S^3}},$$
with $\T_{X^{\#}/B}|_{X_b^{\#}}$ the tangent bundle of
$X_b^{\#}$. Now both $\T_{S^3}$ and $\T_{X^{\#}/B}$ are trivial
vector bundles, and so $\T_{X^{\#}}$ is also. Thus $w_2(X^{\#})
=0=w_2(X)$. $\bullet$

The last calculation is that of the first Pontrjagin class of $X$.
Before carrying this out, we would like to put
an orientation on $Crit(f)$ given an orientation on $X$. This will
allow us to define a class $Crit(f)\in H^4(X,\boldz)$ as the
Poincar\'e dual of the corresponding oriented 2-cycle. We do this
using the following definition. 

\proclaim Definition 2.14. 
\item{(1)}
If $g:X'\rightarrow D$ is a $T^2$-fibration
over a disk as constructed in Example 2.6, (1), there is
an immersion $i:S^2\rightarrow X'$ onto the singular fibre $X_0'$.
This immersion fails to be an embedding precisely at the singular
point of $X_0'$, where two sheets of the immersed $S^2$ cross.
This local structure can be seen in the map
$f:{\bf C}^2\rightarrow {\bf R}^2$ of Example 2.3, or equivalently,
in the map ${\bf C}^2\rightarrow {\bf C}$ given by $(z_1,z_2)\mapsto
z_1z_2$. Given an orientation on $S^2$,
we call the orientation on $X'$ for which these two sheets intersect
positively the {\it positive orientation} on $X'$.
Note this is independent of the choice of orientation
on $S^2$. 
\item{(2)}
Given $f:X\rightarrow B$ a $T^3$-fibration produced by Theorem 2.1,
and an orientation on $X$, then $Crit(f)$ is a union of connected two-manifolds
$\bigcup S_i$ meeting at most at points. We can orient each $S_i$ as
follows.
For a point $b\in\Delta_g\cap f(S_i)$, there is a neighbourhood
$U=D\times (0,1)$ of $b$ such that $f^{-1}(U)=X'\times (0,1)\times S^1$,
so that $f$ is induced by a map $g:X'\rightarrow D$, as in Definition 1.2.
Take the positive orientation on $X'$ over $U$. Then $S_i\cap f^{-1}(U)$ is
the surface $S=Crit(g)\times
(0,1)\times S^1$ meeting $X'\times \{1/2\}\times \{p\}$
transversally. Orient $S$ so that it meets this latter surface positively.
If each $S_i$ is orientable, this gives an orientation on $S_i$ for
each $i$, and hence
makes $Crit(f)$ into an oriented two-cycle.
We call this orientation
on $Crit(f)$ the {\it canonical orientation}. We shall see in Theorem 2.17
that this orientation
does not depend on the choice of $b\in \Delta_g\cap f(S_i)$,
and that $S_i$ is orientable.

{\it Example 2.15.} 
Let $g:X'\rightarrow D$ be as constructed in Example 2.6, 
(1), a well-behaved $T^2$-fibration
with a singular $I_1$ fibre. Put an orientation $\eta_D$ on
$D$ and choose a basis $e_1,e_2$
of $H_1(X'_b,\boldz)$ so that a positive
(counter-clockwise) loop $\gamma$ in $D^*$ 
has monodromy $\pmatrix{1&1\cr 0&1\cr}$
in this basis. Take $e_1\wedge e_2$
to give an orientation on a fibre. Then the induced
orientation $\eta_D\wedge e_1\wedge e_2$ on 
$X'$ is positive. To see this, one can for example look
at the standard family of degenerating elliptic curves,
identifying $D$ with the unit disk in ${\bf C}$ with the standard
orientation, and taking
a family with period $\tau={1\over 2\pi i}\log z$. Then the above
described orientation on $X'$ is the same as that given
by the complex structure, and hence $X'$ is positively oriented.

Let $f:X=X'\times (0,1)\times S^1\rightarrow B=D\times (0,1)$
be the induced fibration, with $X'$ positively oriented by the orientation
$\eta_{X'}$, and
$(0,1)\times S^1$ given some orientation $\eta_{(0,1)\times S^1}$.
Suppose $X$ is assigned the induced product orientation $\eta_{X'}\wedge
\eta_{(0,1)\times S^1}$.
Then the canonical orientation on $Crit(f)=Crit(g)\times
(0,1)\times S^1$ is then clearly the induced orientation $\eta_{(0,1)\times
S^1}$.

{\it Example 2.16.} This is needed in the proof of Theorem 4.4. Consider
a $T^3$-fibration $f:\bar X\rightarrow B={\bf R}^3$ with a line of
type $(2,2)$ fibres constructed via the method of Example 2.8.
In other words, one takes $\bar Y=S^1\times B$, $p\in S^1$ a point,
$l\subseteq {\bf R}^3$ a line, and obtain $\bar X$ by compactifying
a $T^2=T(L)$-bundle over $Y=S^1\times B-\{p\}\times l$ with Chern class $c_1\in L$ primitive. This yields a map $\bar\pi:\bar X\rightarrow\bar Y$,
and $f$ is obtained as
the composition of $\bar\pi$ and the projection $\bar Y\rightarrow B$.

Let $H\subseteq {\bf R}^3$ be a plane containing $l$, and decompose
$H$ into two closed half-planes $H_1$ and $H_2$, with $H=H_1\cup H_2$
and $l=H_1\cap H_2$. Let $D_i=\bar\pi^{-1}(\{p\}\times H_i)$. Then $D_i$
is a smooth four-manifold, and $D_1$ and $D_2$ intersect transversally
along $Crit(f)$. Orient $D_1$ and $D_2$ as follows. Choose a basis
$v_1,v_2$ of the tangent space $\T_H$ of $H$,
and a basis $w_1,w_2$ of $L\otimes_{\boldz}{\bf R}$.
For $p_i\in \{p\}\times H_i$, the tangent space $\T_{D_i,y_i}$ to
a point $y_i\in \bar\pi^{-1}(p_i)$ can be decomposed as $\T_H\oplus
(L\otimes_{\boldz} {\bf R})$, with $\bar\pi_*$ mapping $\T_H$ isomorphically
to the tangent space of $\{p\}\times H$. Orient $D_i$ using
$\eta_{D_i}=v_1\wedge w_1\wedge v_2\wedge w_2$. The important thing here
is that this gives orientations on $D_1$ and $D_2$ which are consistent
in some sense.
I claim that the oriented intersection $D_1\cap D_2$ is the canonical
orientation on $Crit(f)$, irrespective of the orientation on
$\bar X$. (Of course, the orientation of $D_1\cap D_2$ depends on
the orientation on $\bar X$, as does the canonical orientation on
$Crit(f)$.)

To see this, write as usual $\bar X=X'\times {\bf R}\times S^1$,
with $g:X'\rightarrow {\bf R}^2$ as in Example 2.6, (1). By construction of
$X'$, an open neighbourhood $U$ of $Crit(g)$ can be identified
with ${\bf C}^2$, so that $g:U\rightarrow {\bf R}^2$ coincides with the
Harvey-Lawson map $f_{2.3}:{\bf C}^2\rightarrow {\bf R}^2$ of
Example 2.3. Furthermore, this can be done so that 
$D_i\cap (U\times \{q\})$ (for a point $q\in {\bf R}\times S^1$)
coincides with $\{z_i=0\}\subseteq {\bf C}^2$. Now the positive
orientation $\eta_{X'}$ on $X'$ coincides with the orientation induced
by the complex structure on ${\bf C}^2\cong U$. Thus in particular,
with orientations $\eta_{D_i\cap (U\times \{q\})}$ being
induced by the complex structure on $\{z_i=0\}$, $D_1\cap (U\times \{q\})$
and $D_2\cap (U\times \{q\})$ intersect positively at $Crit(g)$. Now
choose an orientation $\eta_{{\bf R}\times S^1}$ on ${\bf R}\times S^1$
so that $\eta_{X'}\wedge\eta_{{\bf R}\times S^1}$ is the given
orientation on $X$. Then with $D_i$ oriented by
$\eta_{D_i\cap (U\times\{q\})}\wedge \eta_{{\bf R}\times S^1}$, 
it is clear that the orientation on the intersection $D_1\cap D_2
=Crit(g)\times{\bf R}\times S^1=Crit(f)$ is $\eta_{{\bf R}\times S^1}$,
which by Example 2.15 is the canonical orientation on 
$Crit(f)$. On the other hand, the $\eta_{D_i}$ either both coincide with
$\eta_{D_i\cap (U\times\{q\})}\wedge \eta_{{\bf R}\times S^1}$
or both have the opposite sign. Hence, with the orientations $\eta_{D_1}$
and $\eta_{D_2}$, the induced intersection
orientation on $D_1\cap D_2$ coincides with the
canonical orientation on $Crit(f)$.

\proclaim Theorem 2.17. If $f:X\rightarrow B$ is as in Theorem 2.1,
$B\cong S^3$, and $Crit(f)=\bigcup S_i$ as in Definition 2.14,
then
\item{(1)} each $S_i$ is orientable and the canonical orientation
is well-defined.
\item{(2)} $p_1(X)=-2Crit(f)\in H^4(X,{\bf Q})$. 
\item{(3)} If $X$ has an almost complex structure in which 
$c_1(X)=0$, then $c_2(X)=Crit(f)$.

Proof. (3) follows from (2) and the
fact that $p_1(X)=p_1(\T_X)=-c_2(\T_X\otimes_{\bf R} {\bf C})
=-2c_2(\T_X)+c_1(\T_X)^2$ if $\T_X$ is already a complex vector
bundle.

For (2), the argument of the proof of Proposition 2.13.
already shows that $p_1(X)$ is a cohomology class supported on
$Crit(f)$. Thus we will be able to reduce the calculation
of $p_1(X)$ to a local calculation near $Crit(f)$. To compute
$p_1(\T_X)=-c_2(\T_X\otimes_{\bf R}{\bf C})$, note that
$ch(\T_X\otimes_{\bf R}{\bf C})=6-c_2(\T_X\otimes_{\bf R}{\bf C})
\in H^*(X,{\bf Q})$, as $c_1(\T_X\otimes_{\bf R}{\bf C})=
c_3(\T_X\otimes_{\bf R}{\bf C})=0$. Thus it is sufficient to calculate the
degree 4 part of $ch(\T_X\otimes_{\bf R}{\bf C})$.

Let $N\subseteq X$ be a small open
tubular neighbourhood of $Crit(f)$. One then has an exact sequence
in K-theory (see [13], I, \S 9 for notation used)
$$K(X,X-N)\rightarrow K(X)\rightarrow K(X-N).$$
From the proof of Proposition 2.13, there is an isomorphism
$\sigma:(\T_X\otimes {\bf C})|_{X-N}\rightarrow {\bf C}_{X-N}^6$
of vector bundles, where ${\bf C}_{X-N}$ denotes the trivial
vector bundle on $X-N$. Thus $[\T_X\otimes_{\bf R}{\bf C}]-[{\bf C}_X^6]
\in K(X)$ maps to zero in $K(X-N)$, and hence $[\T_X\otimes_{\bf R}{\bf C}]
-[{\bf C}_X^6]$ comes from an element of $K(X,X-N)$. Such an element is
represented by $[\T_X\otimes{\bf C},{\bf C}^6_X,\sigma]$.
Now there is a commutative square
$$\matrix{K(X,X-N)&\mapright{}&K(X)\cr
\mapdown{ch}&&\mapdown{ch}\cr
H^*(X,X-N,{\bf Q})&\mapright{}&H^*(X,{\bf Q})\cr}$$
A word of explanation is required for the first vertical map.
To define $ch:K(X,X-N)\rightarrow H^*(X,X-N,{\bf Q})$, use
the fact that $K(X,X-N)=\tilde K(X/(X-N))$, with base point $p$
the image of $X-N$ in $X/(X-N)$. Then there is the usual
$ch:\tilde K(X/(X-N))\rightarrow \tilde H^*(X/(X-N),{\bf Q})
=H^*(X/(X-N),\{p\},{\bf Q})$, and this can be composed with the pull-back
$H^*(X/(X-N),\{p\},{\bf Q})\rightarrow H^*(X,X-N,{\bf Q})$.

Thus, to compute the degree 4 part of $ch(\T_X\otimes_{\bf R}{\bf C})
=ch([\T_X\otimes_{\bf R}{\bf C}]-[{\bf C}_X^6])$, we compute
the degree 4 part of $ch([\T_X\otimes_{\bf R}{\bf C},{\bf C}_X^6,\sigma])$
in $H^4(X,X-N,{\bf Q})$, and take its image under the map to $H^4(X,{\bf Q})$.
Now $H^4(X,X-N,{\bf Q})\cong H_2(N,{\bf Q})$ by Lefschetz duality.
First suppose that
if $Crit(f)=\bigcup S_i$ is the decomposition of $Crit(f)$ into
connected 2-manifolds, all the $S_i$'s are orientable.
Then $H_2(N,{\bf Q})$ is generated by the
homology classes of the $S_i$ (with some chosen orientation).
In particular, given a class $\alpha\in H^4(X,X-N,{\bf Q})$,
it can be determined by evaluating it on suitable 4-cycles
in $H_4(X,X-N,{\bf Q})$. Explicitly, fixing $i$,
let $b_i\in f(S_i)\cap \Delta_g$, and let $U_i$ be an open neighbourhood
of $b$ as in Definition 1.2, so that $f^{-1}(U_i)=X'\times (0,1)\times S^1$.
Fixing $p\in S^1$, let $X_i'=(X'\times \{1/2\}\times\{p\})\cap \overline{N}$.
Furthermore, taking the positive orientation on $X'$, we obtain
an orientation on $X_i'$, defining a class $X_i'\in H_4(X,X-N,{\bf Q})$.
We take the canonical orientation on $S_i$, so that $X_i'\cdot S_i=+1$. Thus
given a cohomology class $\alpha\in H^4(X,X-N,{\bf Q})$,
we can write 
$$\alpha=\sum (X_i\cdot\alpha) S_i.$$

We apply this analysis to  $ch([\T_X\otimes_{\bf R}{\bf C},{\bf C}^6,\sigma])$,
noting that this lives entirely in $H^4(X,X-N,{\bf Q})$.
Restricting $\T_X\otimes_{\bf R} {\bf C}$ and $\sigma$ to $X_i'$,
one obtains an element of $K(X_i',\partial X_i')$, and one has an exact
sequence
$$\exact{\T_{X_i'}}{\T_X|_{X_i'}}{\N}.$$
The normal bundle $\N$ is trivial. In addition, via the fibration
$g_i:X_i'\rightarrow D$ one obtains a trivialization $\sigma'$
of $\T_{X_i'}$ away from $Crit(g_i)$, as in Proposition 2.13,
compatible with the trivialization of $\T_X|_{X_i'}$. Thus
$$[\T_X\otimes_{\bf R}{\bf C}|_{X_i'},{\bf C}^6_{X_i'},\sigma|_{\partial X_i'}]
=[\T_{X_i'}\otimes_{\bf R}{\bf C},{\bf C}^4_{X_i'},\sigma']$$
in $K(X_i',\partial X_i')$. Note that $H^p(X_i',\partial X_i',{\bf Q})=0$
if $p\not= 4$ and ${\bf Q}$ if $p=4$.
We observe that
$ch([\T_{X_i'}\otimes_{\bf R}{\bf C},{\bf C}^4_{X_i'},\sigma'])=
-2\in H^4(X_i',\partial X_i',{\bf Q})$
(with the positive orientation). To see this, one  either performs a 
direct calculation or observes
that on an elliptically fibred K3 surface $Y$ with 24 $I_1$ fibres,
$\int_Y ch(\T_Y\otimes {\bf C})=-48$, and thus
each singular fibre must contribute $-2$ to $ch$. Thus
$X_i'\cdot ch([\T_X\otimes {\bf C},{\bf C}^6,\sigma])=-2$,
so $ch([\T_X\otimes {\bf C}, {\bf C}^6,\sigma])
=-\sum 2S_i=-2Crit(f)$. (Note this makes it clear that the
orientation on $S_i$ is independent of the choice of $b_i\in\Delta_g
\cap f(S_i)$; otherwise we would get two different values
for $ch([\T_X\otimes_{\bf R}{\bf C},{\bf C}_X^6,\sigma])$.)
Thus the degree 4 piece of
$ch(\T_X\otimes {\bf C})$ is $-2Crit(f)
\in H^4(X,{\bf Q})$. Thus $p_1(X)=-2Crit(f)$. 

Finally suppose that some of the $S_i$'s were not orientable.
Then $H_2(N,{\bf Q})$ is generated only by those $S_i$'s which are
orientable. In particular, if $S_i$ is not orientable, then
we can choose an $X_i'$ as before intersecting $S_i$ transversally.
However, such an $X_i'$ is now homologous to zero, and thus
$ch([\T_X\otimes_{\bf R}{\bf C}|_{X_i'},{\bf C}^6_{X_i'},\sigma|_{\partial X_i'}
])=0$, which is a contradiction with the earlier calculation of this
Chern character. Thus all $S_i$'s must be orientable. This proves
(1) and (2). $\bullet$
\bigskip

{\hd \S 3. Warm-up: Local Mirror Symmetry?}

Let $N\cong \boldz^n$, and let $M=\hom(N,\boldz)$ be the dual
lattice. Put $N_{\bf R}:=N\otimes_{\boldz}{\bf R}$. Let $\sigma\subseteq
N_{\bf R}$ be a strongly convex rational polyhedral cone. We will
assume that $\sigma$ is a Gorenstein canonical cone. This means that
if $n_1,\ldots,n_s\in N$ are the set of generators of 1-dimensional
faces of $\sigma$, then there exists an $m_0\in M$ such that
$\langle m_0,n_i\rangle=1$ for all $i$ and $\langle m_0,n\rangle\ge 1$
for all $n\in\sigma\cap (N-\{0\})$. Denote by $Y_{\sigma}$ the
corresponding affine toric variety. $Y_{\sigma}$ has Gorenstein canonical
singularities. Let $P$ be the convex hull of $n_1,\ldots,n_s$, and suppose
there is a triangulation of $P$ such that the fan $\Sigma$ obtained
as the cone over this triangulation yields a non-singular toric
variety $Y_{\Sigma}$. Then $Y_{\Sigma}\rightarrow Y_{\sigma}$ is a crepant
resolution, and $K_{Y_{\Sigma}}=0$.

We will construct a topological $T^{n-1}\times{\bf R}$ fibration on 
$Y_{\Sigma}$ which appears to be a natural generalization of the 
Harvey-Lawson example given in Example 2.3. 
This will prove to be a useful guide
in constructing a $T^3$-fibration on the quintic in \S 4,
and also offer a possible insight into the phenomenon of local mirror
symmetry discussed in [4].

Let $T_{\bf C}(N):=N\otimes_{\boldz} {\bf C}^*=\hom(M,{\bf C}^*)$,
$T(N):=N\otimes_{\boldz} U(1)$. There is a map $-\log |\cdot|:T_{\bf C}(N)
\rightarrow N_{\bf R}$ with kernel $T(N)$.
There is the standard action of $T_{\bf C}(N)$ on $Y_{\Sigma}$, which
restricts to a $T(N)$ action. It is also standard (see [16], \S 1.3)
that $Mc(\Sigma):=Y_{\Sigma}/T(N)$ is a real $n$-dimensional manifold
with corners. Let $q:Y_{\Sigma}\rightarrow Mc(\Sigma)$ be the quotient
map. By [16], Proposition 1.8, $N_{\bf R}$ acts continuously
on  $Mc(\Sigma)$ in such a way that $q$ is equivariant with respect
to $-\log|\cdot|:T_{\bf C}(N)\rightarrow N_{\bf R}$. The orbits of the 
$N_{\bf R}$ action on $Mc(\Sigma)$ are isomorphic to $N_{\bf R}/{\bf R}\tau$
for $\tau$ ranging through all cones in the fan $\Sigma$.
Furthermore, if $\O_{\tau}$ is the orbit corresponding to $\tau$,
then $\O_{\tau}\subseteq \overline{\O_{\tau'}}$ if and only if
$\tau'\subseteq\tau$.

{\it Example 3.1.} (1) Let $\sigma=\Sigma$ 
be generated by $e_1,\ldots,e_n\in N$ a basis for $N$
so that $Y_{\Sigma}=Y_{\sigma}={\bf C}^n$. Then
$Mc(\Sigma)$ is naturally identified with 
${\bf R}^n_{\ge 0}$, via $\mu:Y_{\Sigma}\rightarrow
{\bf R}^n_{\ge 0}$ given by $\mu(z_1,\ldots,z_n)=(|z_1|,\ldots,|z_n|)$.
The orbit ${\bf R}^n_{>0}$ is $\O_{\{0\}}$, $\{0\}\times {\bf R}^{n-1}_{>0}
=\O_{{\bf R}_{\ge 0}e_1}$, $\{(0,0)\}\times {\bf R}_{>0}^{n-2}
=\O_{{\bf R}_{\ge 0}e_1+{\bf R}_{\ge 0}e_2}$, etc.

(2) Let $Y_{\sigma}={\bf C}^3/\boldz_3$, where the group action is
generated by $(z_1,z_2,z_3)\mapsto (\xi z_1,\xi z_2,\xi z_3)$ with
$\xi$ a primitive third root of unity.
The cone $\sigma$ is the cone generated by $$(1,0,0),(0,1,0),(0,0,1)$$
in the lattice $\boldz^3+{1\over 3}(1,1,1)\boldz$. The resolution is
given by the fan $\Sigma$ determined by the triangulation with interior
vertex ${1\over 3}(1,1,1)$ depicted here:
$$\epsfbox{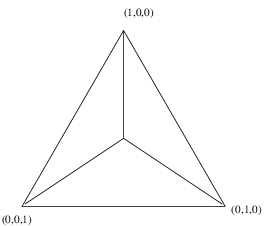}\leqno{\hbox{Figure $3.1$}}$$
$Mc(\Sigma)$ is homeomorphic (but not as a manifold with corners)
to ${\bf R}_{\ge 0}\times {\bf R}^2$, where ${\bf R}_{>0}\times
{\bf R}^2$ is the dense orbit corresponding to the cone $\{0\}$,
and $\{0\}\times {\bf R}^2$ decomposes into orbits as depicted in
Figure 3.2.
$$\epsfbox{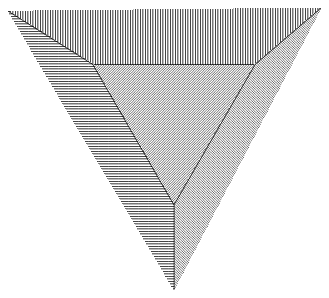}\leqno{\hbox{Figure $3.2$}}$$

In general, $Mc(\Sigma)$ is homeomorphic to ${\bf R}_{\ge 0}\times
{\bf R}^{n-1}$, with the orbit $\O_{\{0\}}={\bf R}_{>0}\times {\bf R}^{n-1}$.
The boundary of $Mc(\Sigma)$ is the union of codimension $\ge 1$
orbits. It will be useful to identify $\partial Mc(\Sigma)$
with the $n-1$-plane containing the polyhedron $P$ in such a way that
the union of codimension $\ge 2$ orbits
of $N_{\bf R}$ is the $n-2$-skeleton of the dual cell complex of
the triangulation of $P$ used to produce $\Sigma$. For example, we picture
the 1-skeleton of Figure 3.2 as the dark lines in
$$\epsfbox{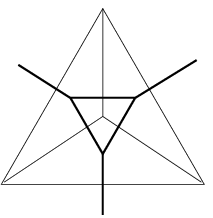}\leqno{\hbox{Figure $3.3$}}$$

Now let $N_{m_0}\subseteq N$ be the sublattice given by
$$N_{m_0}=\{n\in N|\langle m_0,n\rangle=0\}.$$
Then $q$ factors via
$$Y_{\Sigma}\mapright{q_1} Y_{\Sigma}/T(N_{m_0})\mapright{q_2}
Mc(\Sigma).$$
Consider first the map $q_2$. The fibres of $q$ over $\O_{\{0\}}$
are isomorphic to $T(N)$; thus the fibres of $q_2$ over $\O_{\{0\}}$
are circles. On the other hand, for a point $y\in Mc(\Sigma)$
in a codimension one orbit $\O_{\tau}$, $\tau$ a 1-dimensional cone
in $\Sigma$, $y'\in q_2^{-1}(y)$, $y''\in q_1^{-1}(y')$, 
the stabilizer of $y''$ in the group $T(N)$ is ${\bf R}\tau/({\bf R}\tau\cap
N)$. Since $({\bf R}\tau)\cap N_{m_0}=\{0\}$, the fibre
$q_1^{-1}(y')$ is an $n-1$ dimensional torus, as is the fibre
$q^{-1}(y)$. Thus $q_2^{-1}(y)=\{y'\}$. Thus we see that
$q_2^{-1}(\partial Mc(\Sigma))\rightarrow\partial Mc(\Sigma)$
is a homeomorphism. In particular, $Y_{\Sigma}/T(N_{m_0})$
can be identified with ${\bf C}\times {\bf R}^{n-1}$, with $q_2$
given by $(z,x)\mapsto (|z|,x)\in {\bf R}_{\ge 0}\times {\bf R}^{n-1}$.
We stress of course that this identification is purely topological,
and this may not be the most natural thing to do.

Now let $B={\bf R}^n$ and let $f:Y_{\Sigma}\rightarrow B$ be
the composition of $q_1:Y_{\Sigma}\rightarrow Y_{\Sigma}/T(N_{m_0})
\cong {\bf C}\times {\bf R}^{n-1}$ with $r:(z,x)\mapsto (\im z,x)
\in {\bf R}\times {\bf R}^{n-1}=B$. Since $q_2$ identifies
$\{0\}\times{\bf R}^{n-1}\subseteq {\bf C}\times {\bf R}^{n-1}$
with $\{0\}\times {\bf R}^{n-1}\subseteq {\bf R}_{\ge 0}\times {\bf R}^{n-1}$,
we can use $q_2$ to
identify the union of codimension $\ge 2$ orbits in $Mc(\Sigma)$
with a subset $\Delta'$ of $Y_{\Sigma}/T(N_{m_0})$.
Let $\Delta=r(\Delta')\subseteq B$. It is then clear that $f:Y_{\Sigma}
\rightarrow B$ is a fibration with general fibre $T^{n-1}\times{\bf R}$,
but for a point $b\in\Delta$ in a dimension $k$ orbit, 
$f^{-1}(b)\cong T^{n-1}\times{\bf R}/T^{n-1-k}\times \{0\}$,
by which we mean we divide the torus $T^{n-1}\times \{0\}$
out by the action of a $T^{n-1-k}$. 

{\it Example 3.2.} (1) Continuing Example 3.1, (1),
let $\gamma:{\bf R}^n_{\ge 0}\rightarrow {\bf R}_{\ge 0}\times {\bf R}^{n-1}$
be given by $\gamma(x_1,\ldots,x_n)=(\prod x_i,x_1^2-x_2^2,\ldots,
x_1^2-x_n^2)$.
It is easy to see that $\gamma$ is a homeomorphism. We can identify
$Y_{\Sigma}/T(N_{m_0})$ with ${\bf C}\times {\bf R}^{n-1}$ so that the
map $q_1$ is given by $q_1(z_1,\ldots,z_n)=(\prod z_i,|z_1|^2-|z_2|^2,\ldots,
|z_1|^2-|z_n|^2)\in {\bf C}\times {\bf R}^{n-1}$. Thus $f:{\bf C}^n\rightarrow
B$ is given by
$$f(z_1,\ldots,z_n)=(\im \prod z_i,|z_1|^2-|z_2|^2,\ldots,|z_1|^2-|z_n|^2).$$
This is precisely the Harvey-Lawson example given in Example 2.3.

(2) Continuing Example 3.1, (2), we note we obtain a $T^2\times {\bf R}$
fibration with discriminant locus $\Delta$ as depicted by the dark lines in
Figure 3.3.
\bigskip

In the three-dimensional case, the map $Y_{\Sigma}-q_1^{-1}(\Delta')
\rightarrow Y_{\Sigma}/T(N_{m_0})-\Delta'$ is a principal $T^2$-fibre bundle,
with fibres naturally isomorphic to $T(N_{m_0})$, and $Y_{\Sigma}
\rightarrow Y_{\Sigma}/T_{N_{m_0}}$ a compactification of this.
Recall from Proposition 2.9 that $H^2(Y_{\Sigma}/T(N_{m_0})-\Delta',
N_{m_0})
=H^3_{\Delta'}(Y_{\Sigma}/T(N_{m_0}),N_{m_0})\cong H_1^{BM}(\Delta',N_{m_0})$,
whose elements can be identified with simplicial 1-chains on $\Delta'\cup
\{\infty\}$ whose boundary is supported on $\infty$.

\proclaim Proposition 3.3. Up to sign, the 
Chern class of the $T(N_{m_0})$-fibre bundle
$q_1:Y_{\Sigma}-q_1^{-1}(\Delta')\rightarrow Y_{\Sigma}/T(N_{m_0})
-\Delta'$ is the 1-chain $c_1$ on $\Delta'\cup\{\infty\}$
described as follows. Choose an orientation on $P$.
View $\Delta'$ as the 1-skeleton
of the dual cell complex of the triangulation of $P$ determining
$\Sigma$.
Each oriented edge $E$ of $\Delta'$ intersects
a unique edge $\langle\tau_i,\tau_j\rangle$ in the triangulation of $P$,
oriented so that $E$ and $\langle\tau_i,\tau_j\rangle$ intersect positively
with respect to the chosen orientation on $P$. The 1-chain $c_1$ assigns
to the oriented edge $E$ the element $\tau_j-\tau_i$ of $N_{m_0}$.
(See Figure 3.4)
This $1$-chain has the property that all coefficients of edges
are primitive in $N_{m_0}$, and the three coefficients associated
to edges with a common vertex span $N_{m_0}$. In particular, the 
compactification $q_1:Y_{\Sigma}\rightarrow Y_{\Sigma}/T(N_{m_0})$
of $Y_{\Sigma}-q_1^{-1}(\Delta')\rightarrow Y_{\Sigma}/T(N_{m_0})-\Delta'$
coincides with that given in Proposition 2.9.

$$\epsfbox{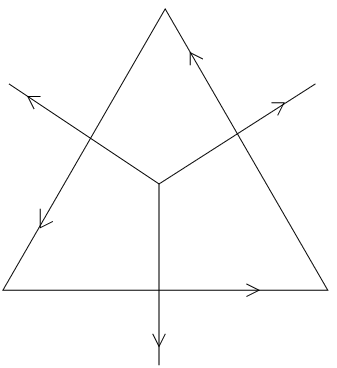}\leqno{\hbox{Figure $3.4$}}$$

Proof. An edge $E$ of $\Delta'$ corresponds to
a one-dimensional cone $\tau$ generated by some
$\tau_i$ and $\tau_j$, and $\langle\tau_i,\tau_j\rangle$
is the unique edge of the triangulation of $P$ intersecting $E$. 
Now if $y\in\O_{\tau}\subseteq Mc(\Sigma)$,
$y'\in q^{-1}(y)$, then the stabilizer of $y'$ in $T(N_{m_0})$ is
$${(N_{m_0}\otimes{\bf R})\cap {\bf R}\tau\over N_{m_0}\cap {\bf R}\tau}
=T((\tau_i-\tau_j)\boldz).$$
Thus the coefficient in $c_1$ of the
oriented edge $E$ is $\pm(\tau_i-\tau_j)\in N_{m_0}$. Once
the sign is chosen for one edge, there is a unique uniform choice
of sign for every oriented edge to ensure that $c_1$ is an element of
$H_1^{BM}(\Delta',N_{m_0})$. The method of assigning signs
given in the statement of the Proposition is one of the two possible choices,
with the opposite sign arising from the opposite choice of orientation
on $P$.

For the last statement, observe that
given a three-dimensional cone $\tau\subseteq\Sigma$, $\tau$ is
a simplicial cone, as $Y_{\Sigma}$ is non-singular. In particular,
$\tau$ is spanned by $e_1,e_2,e_3$, a basis of $N$, and $m_0=e_1^*
+e_2^*+e_3^*$. Let $\tau_i$ be the cone spanned by $e_j$ and $e_k$,
$\{i,j,k\}=\{1,2,3\}$. The coefficients associated to the three edges
corresponding to $\tau_i$ are $e_2-e_3$, $e_3-e_1$ and $e_1-e_2$.
These span $N_{m_0}$, as desired. $\bullet$

{\it Remark 3.4.} This is a rather technical remark which will only
be used in the proof of Theorem 4.4.
Note that the one-dimensional cones $\tau$
of $\Sigma$ correspond to divisors in $Y_{\Sigma}$. One way to identify
these divisors is as follows. Identifying $\O_{\tau}\subseteq
Mc(\Sigma)$ with $q_2^{-1}(\O_{\tau})\subseteq Y_{\Sigma}/T(N_{m_0})$,
we see that the corresponding divisor $D_{\tau}$ is 
$q_1^{-1}(\overline{\O_{\tau}})$.

Now $D_{\tau}$ has an orientation induced by the complex structure
on $D_{\tau}$. It is easy to see that in the $n=3$ case
this orientation can be described as follows: Let $\langle \tau,
\tau_1,\tau_2\rangle$ be a 2-simplex in the triangulation of $P$.
Recall that we identified $\partial Mc(\Sigma)$ with the plane
containing $P$ in such a way that the union of codimension $\ge 2$
orbits is the 1-skeleton of the dual complex of $P$. Thus it is natural
to think of $\tau$ as a point in $\O_{\tau}$ and
$\tau_1-\tau$ and $\tau_2-\tau$ as tangent vectors
$v_1,v_2$ respectively in $\T_{\O_{\tau},\tau}$, the tangent space
of the orbit $\O_{\tau}$ at the point $\tau\in\O_{\tau}$.
On the other hand, the tangent space of $D_{\tau}$ at a point
$y\in q_1^{-1}(\tau)$ can be written as $H\oplus (N_{m_0}\otimes_{\boldz}
{\bf R})$, where $H$ is a subspace mapped isomorphically to $\T_{\O_{\tau},
\tau}$ by $q_{1*}$. Then the orientation on $D_{\tau}$ is
given by $v_1\wedge (\tau_1-\tau)\wedge v_2\wedge (\tau_2-\tau)$,
thinking of $v_1,v_2\in H$, $\tau_1-\tau,\tau_2-\tau\in N_{m_0}\otimes_{\boldz}
{\bf R}$. 

To see this, note this is a local issue, allowing one to reduce to the case
that $\sigma$ is as in Examples 3.1 and 3.2, (1) with $n=3$. Then taking
$\tau=e_1$, $D_{\tau}$ is the divisor $z_1=0$, and then $q_1|_{D_{\tau}}$
is $(z_2,z_3)\mapsto (-|z_2|^2,-|z_3|^2)$. It is then easy to see
that the orientation on $D_{\tau}$ is as described.

Note that changing the sign of the Chern class $c_1$ of the $T(N_{m_0})$-bundle 
$q_1:Y_{\Sigma}-q_1^{-1}(\Delta')\rightarrow Y_{\Sigma}/T(N_{m_0})
-\Delta'$ has the affect of changing the sign of the $T(N_{m_0})$
action, but doesn't change the above orientation (because $n-1=2$ is even).
This is useful because we can get the right orientation on $D_{\tau}$ even
if we haven't gotten the sign of $c_1$ right.

{\it Remark 3.5.}
As a final bit of speculation in the $n=3$ case, 
suppose that there is an inclusion
$Y_{\Sigma}\subseteq \bar Y_{\Sigma}$ with an extension
$\bar f:\bar Y_{\Sigma}\rightarrow B$ of $f$ with $\bar f$
a well-behaved $T^3$-fibration.
Then by Corollary 2.2, a dual fibration $\check f:\check Y_{\Sigma}
\rightarrow B$ can be constructed. Since $\bar f$ only has
singular fibres of type $(2,2)$ and $(1,2)$, all singular fibres of
$\check f$ will be of type $(2,2)$ and $(2,1)$.  By the construction of Example 2.6, (4) of fibres of type $(2,1)$, the
critical locus of $\check f$ is a topological two-manifold $S$ fibred
over $\Delta$. If $\Sigma$ is a cone over a triangulation of the
polygon $P$ contained
in $\{n\in N|\langle m_0,n\rangle=1\}$, then the number of punctures
of $S$ is the number of edges of $P$, and the genus of $S$ is the
number of interior vertices of $P$. Here is a genus 2 example, the
triangulated polygon $P$ shown in lighter lines and $\Delta$ shown
in darker lines:
$$\epsfbox{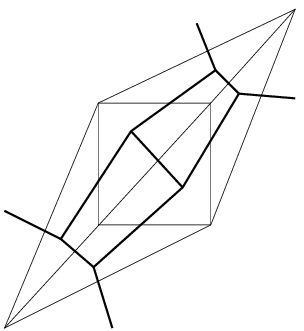}\leqno{\hbox{Figure $3.5$}}$$
(see [4], Figure 3). Now [4] suggests through other arguments
that the ``mirror'' of $Y_{\Sigma}$ is a curve whose genus is precisely
the number of interior vertices of $P$. This suggests that in this
local setting, this ``mirror curve'' is the critical locus of $\check f$.
Because it is quite likely one can write down explicit special Lagrangian
versions of the fibrations given above in some simple cases, this
may prove to be a useful laboratory in which to study the
SYZ conjecture. This will be taken up elsewhere.
\bigskip

{\hd \S 4. A Well-behaved Fibration on the Quintic.}

We will now construct a torus fibration $f:X\rightarrow B$ on the quintic
threefold $X\subseteq \Pfour$. We adopt a different approach to that given
in [17] and [21], essentially constructing a $T^3$-fibration by hand and then
showing the total space of this fibration is diffeomorphic to $X$.
To do so,
we make use of Ruan's description of the monodromy in his fibration,
and construct the fibration directly from this data.

Recall the moment map $\mu:\P^n\rightarrow \Xi$ given by
$$(z_0,\ldots,z_n)\rightarrow {\sum_{i=0}^n |z_i|^2 P_i\over \sum_{i=0}^n |z_i|^2},$$ 
where $P_0,\ldots,P_n$ are the vertices of an $n$-simplex $\Xi\subseteq
{\bf R}^n$. Then it is obvious that the restriction 
of $\mu$ to each linear subspace $\{z_{i_1}=\cdots=
z_{i_{n-p}}=0\}$ is just the moment map of $\P^{p}$,
and the fibre of $\mu$ over a $p$-dimensional face of $\Xi$ is
$T^p$. In particular, with $n=4$, the moment map of $\Pfour$ restricted to 
$F=\{z_0z_1z_2z_3z_4=0\}$ 
has image $\partial\Xi$, and $\mu:F\rightarrow\partial\Xi$
is a three-torus fibration, with fibres dropping rank on lower dimensional
simplices of $\partial\Xi$. We will construct $X$ so as to be a topological
resolution of $F$, and so that $X$ fibres over $\partial\Xi$, with $n=4$.

Begin by considering the graph
$\Gamma\subseteq\partial\Xi$ described by Ruan in [17], \S 4. 
Namely, if $P_{ij}$ is the barycenter
of the 1-simplex $\langle P_i,P_j\rangle$, and $P_{ijk}$ is the
barycenter of the 2-simplex $\langle P_i,P_j,P_k\rangle$, the
vertices of $\Gamma$ are $\{P_{ij},P_{ijk}\}$, and the edges of
$\Gamma$ connect $P_{ij}$ with $P_{i'j'k'}$ whenever $\{i,j\}
\subseteq \{i',j',k'\}$. This graph is trivalent.
Now take a small tubular neighbourhood $N_{\Gamma}$ of $\Gamma$,
and let $U=\partial\Xi-\overline{N_{\Gamma}}$. Let $\sigma_i$
be the maximal face of $\partial\Xi$ given by $\langle P_0,\ldots,
\hat P_i,\ldots,P_4\rangle$, and let 
$$U_i=Int(\sigma_i)-\overline{N_{\Gamma}}.$$
(Here, $Int(\sigma)=\sigma-\partial\sigma$.)
Let 
$$U^i=\{P=\sum a_iP_i\in \partial \Xi|
\hbox{$a_i>a_j$ for all $j\not=i$}\}-\overline{N_{\Gamma}}.$$
Then $U_i\cap U_j=\phi$, $U^i\cap U^j=\phi$, for all $i$ and $j$,
and $U^i\cap U_j=\phi$
if and only if $i=j$. Also $U_i$ and $U^j$ are contractible.
Furthermore, 
$$U=\bigcup_{i=0}^4 (U^i\cup U_i).$$
 
We will begin
by constructing a $T^3$-fibre bundle over $U$, $f_1:X_1\rightarrow U$.
$f_1$ will be trivial over each $U_i$ and $U^i$, as these sets
are contractible, so
$f_1^{-1}(U_i)=T^3\times U_i$ and $f_1^{-1}(U^i)=T^3\times U^i$. However,
we must be precise in identifying these tori. Fix a lattice
$N=\boldz^5$ with basis $e_0,\ldots,e_4$, and let $N_i=N
/\langle e_i,\sum_{j=0}^4 e_j\rangle$. Then $\dual{N_i}$ is naturally
identified with $\{\sum a_je_j^*\in \dual{N}| a_i=\sum a_j=0\}$.
For a lattice $N$, denote by $T(N)$ as before
the torus $(N\otimes_{\boldz}{\bf R})/N$.
Identify $f_1^{-1}(U_i)$ with $T(N_i)\times U_i$
and $f_1^{-1}(U^i)$ with $T(\dual{N_i})\times U^i$. 
Define isomorphisms $T_{ij}:N_i\rightarrow\dual{N_j}$ for $i\not=j$
by taking $e_0,\ldots,\hat e_i,\ldots,\hat e_j,\ldots,e_4$
as a basis for $N_i$ and setting
$$T_{ij}(e_k)=e_k^*-e_i^*.$$
$T_{ij}$ then induces homeomorphisms $T(N_i)\rightarrow T(\dual{N_j})$.
Glue $T(N_i)\times U_i$ with $T(\dual{N_j})\times U^j$ over $U_i\cap U^j$
via this homeomorphism. This gives a $T^3$-bundle $f_1:X_1\rightarrow U$ with
a section.

We compute the monodromy of $f_1:X_1\rightarrow U$
about a loop $\gamma_{ij,k}$
passing around an edge $\langle P_{ij},P_{ijk}\rangle$ of the
graph $\Gamma$. This calculation was also carried out in [17], \S 4.
We first note that this edge lies in a two-simplex which is the common face
of $\sigma_l$ and $\sigma_m$, where $\{i,j,k,l,m\}=\{0,\ldots,4\}$.
In addition, $P_{ij}$ is contained in the 1-simplex spanned by
$P_i$ and $P_j$. Thus we can take $\gamma_{ij,k}$ to be a simple loop based
at $b\in U_l$, passing first into $U^i$, then $U_m$, then $U^j$ and then
back into $U_l$. In Figure 4.1 we have depicted the
loops $\gamma_{ij,k},\gamma_{jk,i}$ and $\gamma_{ik,j}$.
Note that $\gamma_{ij,k}\cdot\gamma_{jk,i}=\gamma_{ik,j}$
in $\pi_1(U,b)$. 
$$\epsfbox{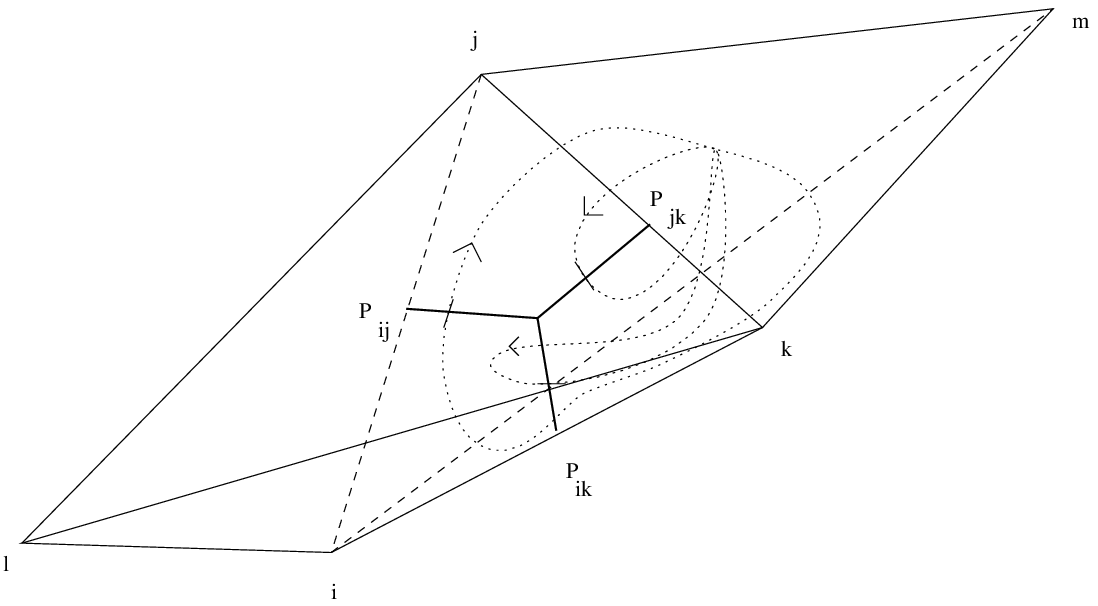}\leqno{\hbox{Figure $4.1$}}$$
Thus the monodromy transformation associated to the loop
$\gamma_{ij,k}$ is $T_{ij,k}=T_{lj}^{-1}T_{mj}T_{mi}^{-1}T_{li}:N_l\rightarrow
N_l$.
In the basis $e_j,e_k,e_m$ of $N_l$, it is easy to calculate that
$$\eqalign{
T_{ij,k}(e_j)&=e_j+5e_m,\cr
T_{ij,k}(e_k)&=e_k,\cr
T_{ij,k}(e_m)&=e_m.\cr}\leqno{(4.1)}$$
Note that from this general formula, we see that
$$\eqalign{
T_{ik,j}(e_j)&=e_j,\cr
T_{ik,j}(e_k)&=e_k+5e_m,\cr
T_{ik,j}(e_m)&=e_m\cr}$$
and
$$\eqalign{
T_{jk,i}(e_j)&=e_j-5e_m,\cr
T_{jk,i}(e_k)&=e_k+5e_m,\cr
T_{jk,i}(e_m)&=e_m.\cr}$$
Now let $N'_{\Gamma}$ be a slightly larger tubular neighbourhood
of $\Gamma$ then $N_{\Gamma}$, and let $U_{ijk}$ be a ball with
center $P_{ijk}$ small enough so it doesn't intersect any dimension 1 faces
of $\partial\Xi$.
Let $V_{ijk}:=N'_{\Gamma}\cap U_{ijk}$. 
To visualise this, we depict the intersection of $V_{ijk}$ with the
face $\langle P_i,P_j,P_k\rangle$ of $\partial\Xi$:
$$\epsfbox{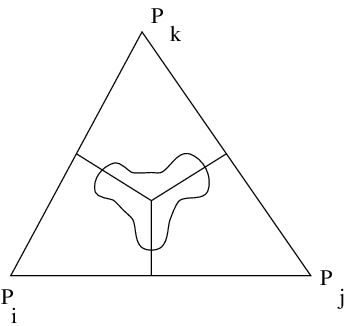}\leqno{\hbox{Figure 4.2}}$$

Now we build a $T^3$-fibration over $V_{ijk}$ as follows: construct
a surface $S_{ijk}\subseteq T(N_l/\langle e_m\rangle)
\times V_{ijk}
=T^2\times V_{ijk}$. This will be a surface
sitting over the following graph $\Delta_{ijk}$ 
in $V_{ijk}\cap \langle P_i,P_j,P_k\rangle$,
represented by the dark lines in the following picture:
$$\epsfbox{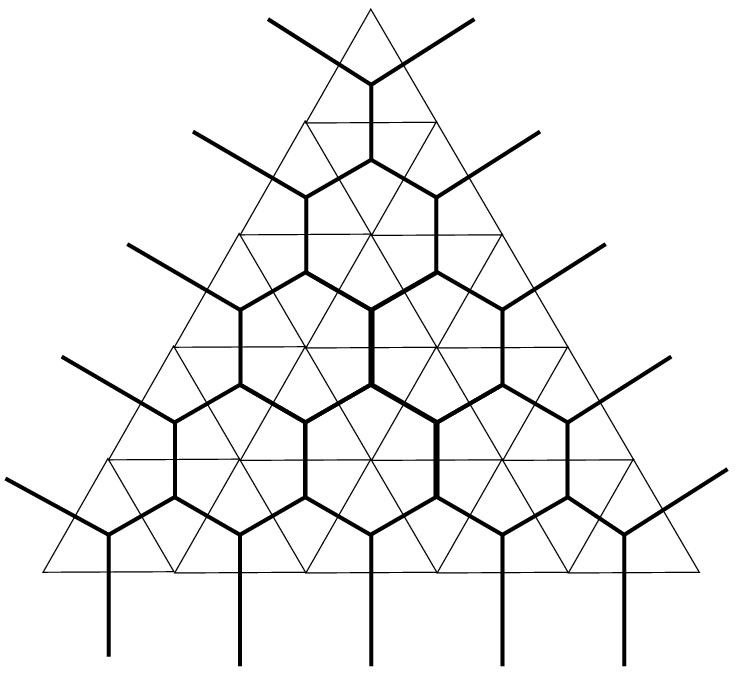}\leqno{\hbox{Figure 4.3}}$$
Over each edge $E$ of the graph, $S_{ijk}$ will be $E\times S^1$,
with $S^1\subseteq T(N_l/\langle e_m\rangle)$ 
of class $e_i$, $e_j$ or $e_k$ depending
on whether $E$ is running northeast-southwest (parallel to the edge
$\langle P_{jk},P_{ijk}\rangle$ of $\Gamma$), northwest-southeast (parallel
to $\langle P_{ik},P_{ijk}\rangle$)
or south-north (parallel to $\langle P_{ij},P_{ijk}\rangle$)
respectively. At each vertex we glue up these cylinders
as in Example 2.6 (4), which we can do by the relation $e_i+e_j+e_k=0$ in
$N_l/\langle e_m\rangle$. The surface $S_{ijk}$ has genus 6 with
15 punctures. Now build an $S^1$-bundle over $T^2\times V_{ijk}
-S$ whose first Chern class is
$1\in H^3_S(T^2\times V_{ijk},\boldz)\subseteq
H^2(T^2\times V_{ijk}-S,\boldz)$. Using Proposition 2.5, this can be
compactified to obtain a $T^3$-fibration $f_{ijk}:X_{ijk}\rightarrow V_{ijk}$ 
with discriminant locus $\Delta_{ijk}$.

Let $b\in V_{ijk}\cap U_l$.
Take as a basis of $H_1(f_{ijk}^{-1}(b),
\boldz)$ elements $e_j,e_k,e_m$, the first two mapping to $e_j$ and
$e_k$ in $T(N_l/\langle e_m\rangle)$, and $e_m$ the class of the $S^1$ fibre
of $X_{ijk}\rightarrow T(N_l/\langle e_m\rangle)\times V_{ijk}$. Then
by construction, possibly after changing the sign of $e_m$, the monodromy
$T$ about a suitably oriented loop around one of the edges of 
$\Delta_{ijk}$ running parallel to $\langle P_{ij},P_{ijk}\rangle$ is given
by $T(e_j)=e_j+e_m, T(e_k)=e_k, T(e_m)=e_m$. The monodromy
about an edge parallel to $\langle P_{jk},P_{ijk}\rangle$ is given
by $T(e_j)=e_j-e_m, T(e_k)=e_k+e_m, T(e_m)=e_m$. Finally, the
monodromy about an edge parallel to $\langle P_{ik}, P_{ijk}\rangle$ is
given by $T(e_j)=e_j, T(e_k)=e_k+e_m,$ and $T(e_m)=e_m$. 

Now $V_{ijk}\cap U$ is homotopic to a thrice-punctured sphere, and $\pi_1(
V_{ijk}\cap U,b)$ is generated by $\gamma_{ij,k},\gamma_{jk,i}$ and
$\gamma_{ik,j}$, with the relation $\gamma_{ij,k}\cdot\gamma_{jk,i}=
\gamma_{ik,j}$. 
Then by the above discussion of the monodromy about edges of $\Delta_{ijk}$,
the monodromy of $f_{ijk}$ around the three loops $\gamma_{ij,k}$,
$\gamma_{jk,i}$ and $\gamma_{ik,j}$ coincides precisely with that of $f_1:X_1
\rightarrow U$. Thus we can glue in $f_{ijk}:X_{ijk}\rightarrow V_{ijk}$
to $f_1:X_1\rightarrow U$.
Doing this for each simplex $\langle P_i,P_j,P_k\rangle$, we obtain a new
fibration $f_2:X_2\rightarrow U'$, with $U'=U\cup \bigcup V_{ijk}$.

Next let $U_{ij}$ be a ball centred at $P_{ij}$, so that 
$\bigcup_{i,j} U_{ij}\cup \bigcup_{i,j,k} U_{ijk}$ covers 
$\overline{N_{\Gamma}}$. Let $V_{ij}:=U_{ij}\cap N'_{\Gamma}$. The situation
near
$V_{ij}$ is shown in Figure 4.4.
$$\epsfbox{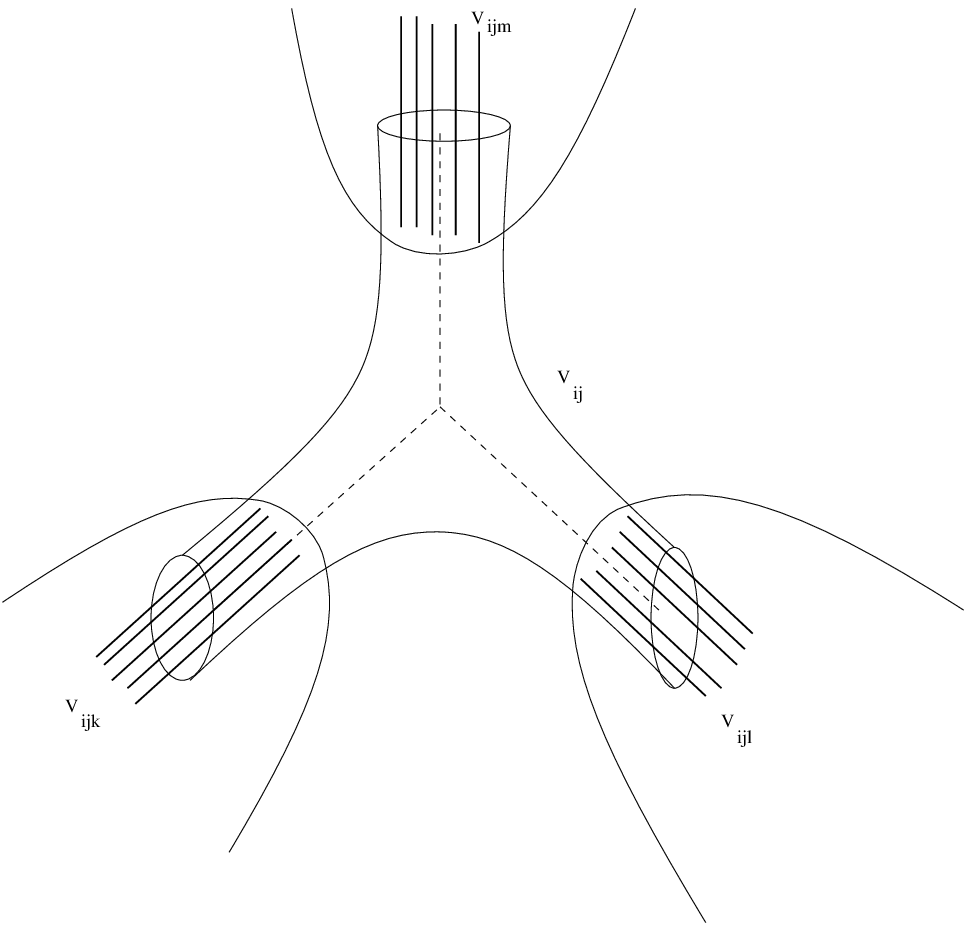}\leqno{\hbox{Figure 4.4}}$$
Note that the relevant monodromy transformations of loops based at $b\in U_l$
 about the three edges
of $\Gamma$ containing $P_{ij}$ are $T_{ij,k}$, $T_{ij,m}$, and
$T_{ij,m}^{-1}\circ T_{ij,k}^{-1}$. Note that
$$\eqalign{
T_{ij,m}(e_j)&=e_j+5e_k,\cr
T_{ij,m}(e_k)&=e_k,\cr
T_{ij,m}(e_m)&=e_m.\cr}$$
We then construct $X_{ij}\rightarrow V_{ij}$ by producing a fibration
with discriminant locus  $\Delta_{ij}$ as in the dark lines of Figure 4.5:
$$\epsfbox{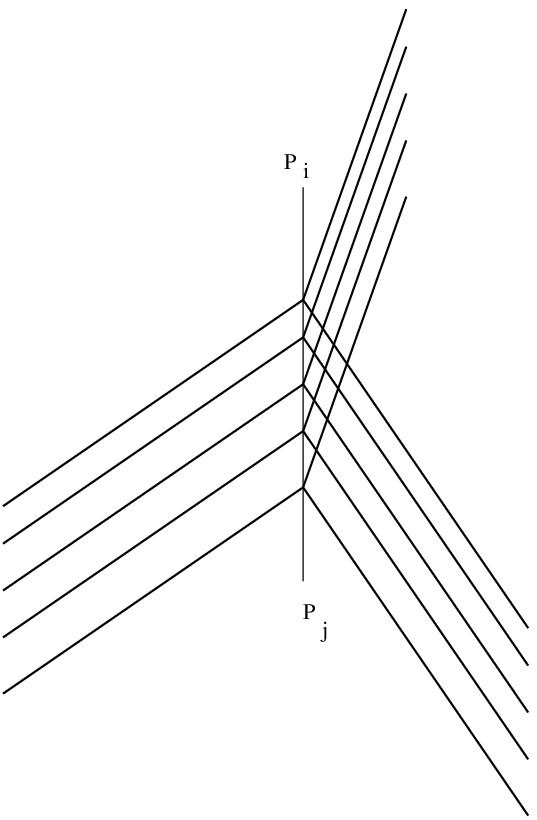}\leqno{\hbox{Figure 4.5}}$$
This is done so that $\Delta_{ijk}\cap V_{ij}=\Delta_{ij}\cap V_{ijk}$.
We build this by taking 5 identical copies
of the construction of Example 2.10 and gluing them in
a chain so that $X_{ij}\rightarrow V_{ij}$ has five
fibres of type $(1,2)$ over $V_{ij}\cap \langle P_i,P_j\rangle$. 
As $e_k$ and $e_m$ are invariant under $T_{ij,k}$ and $T_{ij,m}$,
this can also be viewed as the compactification of a $T^2$-bundle
over $(T(N_l/\langle e_k,e_m\rangle)\times V_{ij})-(\{p\}\times\Delta_{ij})$
for some point $p\in T(N_l/\langle e_k,e_m\rangle)$.
This can be done
with the desired monodromy so that
$X_{ij}\rightarrow V_{ij}$ can be glued in to $f_2:X_2\rightarrow U'$
to obtain finally a fibration $f:X\rightarrow B:=\partial\Xi$ with a section.

We note that as in Remark 2.5, putting a differentiable structure
on $\partial\Xi\cong S^3$, one can put a differentiable structure on $X$.

\proclaim Theorem 4.1. $f:X\rightarrow B$ is
a well-behaved $T^3$-fibration, and there is a continuous
map $\xi:X\rightarrow F=\{z_0z_1z_2z_3z_4=0\}$ 
such that $f=\mu\circ\xi$. Furthermore,
$X$ is diffeomorphic to a non-singular quintic threefold. 

Proof. The only part of admissibility which doesn't follow immediately
from the construction is that $R^3f_*\boldz=\boldz$. However, this is
true locally, and since $B$ is simply connected, it is true globally.
Then well-behavedness is clear by construction.

Before continuing, we introduce some notation. Let $F_{i_0\ldots i_p}:=
F_{i_0}\cap\cdots\cap F_{i_p}$ 
where $F_i
=\{z_i=0\}\subseteq\Pfour$. 
Let $F^{[p]}=
\coprod_{i_0<\cdots<i_p} F_{i_0\cdots i_p}$.
In particular, $F^{[0]}$ is the normalization
of $F$. Let $F^{(p)}=\bigcup_{i_0<\cdots<i_p} F_{i_0\cdots i_p}$.
Let $F_0^{(p)}=F^{(p)}-F^{(p+1)}$. Let
$i_p:F_0^{(p)}\rightarrow F^{(p)}$ be the inclusion map.

We now construct a map $\xi:X\rightarrow F$. Choose a maximal face
$\sigma_l$ of $\partial\Xi$. Then by construction,
$f^{-1}(Int(\sigma_l))=T(N_l)\times Int(\sigma_l)$.
Note that the five-torus
$T(N)$ acts naturally on $\Pfour$ via the standard
diagonal action, and $T(N_l)$ has
a natural induced action on $F_l$. The quotient map $F_l-Sing(F)
\rightarrow (F_l-Sing(F))/T(N_l)$ coincides with the moment map
$\mu:F_l-Sing(F)\rightarrow Int(\sigma_l)$, and allows $F_l-Sing(F)$
to be identified naturally with $T(N_l)\times Int(\sigma_l)
=f^{-1}(Int(\sigma_l))$.

Now also by construction, if $\sigma=\langle P_i,P_j,P_k\rangle$
is a face of $\sigma_i$, and $\{i,j,k,l,m\}=\{0,1,2,3,4\}$,
then there is a natural $S^1$-action on $f^{-1}(Int(\sigma))$
induced by $e_m\in N_l$. Again $f^{-1}(Int(\sigma))/S^1=
T(N_l/\langle e_m\rangle)\times Int(\sigma)$ is naturally identified
with $F_{lm}-Sing(F)$, so that the projection $f^{-1}(Int(\sigma)
)/S^1\rightarrow Int(\sigma)$ coincides with $\mu:F_{lm}-Sing(F)\rightarrow Int(\sigma)$.
If $\tau=\langle P_i,P_j\rangle$ is a dimension
one facet of $\sigma_l$, there is a natural $T^2$ action on $f^{-1}(
Int(\tau))$ induced by $e_k$ and $e_m$ in $N_l$ and $f^{-1}(Int(\tau))
/T^2=T(N_l/\langle e_k,e_m\rangle)\times Int(\tau)$ is
naturally identified with $F_{klm}-Sing(F)$ so that the projection 
$f^{-1}(Int(\tau))/T^2\rightarrow Int(\tau)$ coincides with
$\mu:F_{klm}-Sing(F)\rightarrow Int(\tau)$. 

Now take the topological
quotient of $X$ with respect to the following equivalence relation:
$$\vbox{\narrower\narrower $x\sim y$ if and only if $f(x)=f(y)$ and either (1)
$f(x)$ is in the interior of a three-dimensional face of $\partial\Xi$
and $x=y$, or (2) $f(x)$ is in the interior of a dimension two face of 
$\partial\Xi$ and $x$ and $y$ are identified under the natural
$S^1$-action described above, or (3) $f(x)$ is in the interior of a dimension
one face of $\partial\Xi$ and $x$ and $y$ are identified under
the natural $T^2$ action described above, or (4) $f(x)$ is a vertex
of $\Xi$.

}$$
It is then clear that the quotient of $X$ by this equivalence relation
is homeomorphic to $F$.
Take $\xi:X\rightarrow F$ to be the quotient map. It is also
clear by construction that $f=\mu\circ\xi$.

We note at this point that while $X$ is orientable since $H^6(X,\boldz)
=H^3(B,R^3f_*\boldz)=\boldz$, it is
worth remarking that this orientation may be chosen to be compatible
via the map $\xi$
with the orientation on $F-Sing(F)$ induced by the complex
structure. Let us be completely explicit about
this orientation, for future use. First, consider 
$\mu:F_l-Sing(F)\rightarrow Int(\sigma_l)$, $\sigma_l
=\langle P_i,P_j,P_k,P_m\rangle$.
 Then the fibre
of $\mu$ is $T(N_l)$, and the natural orientation coming from
the complex structure on $F_l$ can be represented as
$$(P_j-P_i)\wedge e_j\wedge (P_k-P_i)\wedge e_k\wedge (P_m-P_i)\wedge e_m.
\leqno{(4.2)}$$
Here, as always, we give orientations by giving an element in the
top exterior power of the tangent bundle.
We abuse notation by identifying $N_l\otimes {\bf R}$ with
the tangent space of a fibre, and identify the tangent space
of $\partial\Xi$ at a point in the interior of a maximal face
of $\partial\Xi$ with the obvious three-dimensional subspace of
${\bf R}^4$. We note this orientation is independent of
the order of $i,j,k$ and $m$.
It is easy to check these orientations are compatible using
the identifications $T_{kj}^{-1}\circ T_{ij}:N_i\rightarrow N_k$.

To show $X$ is diffeomorphic to the quintic, we use Wall's classification
of simply-connected 6-manifolds. To apply this result, we need
to know: (1) $X$ is simply connected. (2) $H^*(X,\boldz)$
is torsion-free. (3) $H^2(X,\boldz)=\boldz$, generated by
a class $H$ with $H^3=5$. (4) $\rank H^3(X,\boldz)=204$.
(5) $w_2(X)=0$. (6) $p_1(X)$ coincides with the first Pontrjagin class
of the quintic. If these six conditions hold, then it follows from
[20], Theorem 5 that $X$ is diffeomorphic to the quintic. We treat each of these
in turn.

(1) This follows from Theorem 2.12: the precise description of the monodromy
shows there are no invariant 1-cycles modulo any $n$. 

(2)-(4) 
The map $\xi:X\rightarrow F$ is proper, with fibres over points of
$F^{(p)}_0$ being $p$-tori except for points of $\xi(Crit(f))$. Furthermore,
$\xi|_{Crit(f)}$ is a homeomorphism onto its image.  Note $\xi(Crit(f))$ is
a union of two-manifolds $\bigcup C_{ij}$, $C_{ij}\subseteq F_{ij}$
of genus 6, $C_{ijk}:=C_{ij}\cap F_{ijk}=C_{ik}\cap F_{ijk}
=C_{jk}\cap F_{ijk}$ a set of five points, and $C_{ij}\cap F^{(3)}
=\phi$.

We wish to use the Leray spectral sequence for $\xi$ to compute the cohomology 
of $X$. We will need to understand the sheaves $R^n\xi_*\boldz$.
First, let $p\in F_0^{(i)}$ and let $q$ be a nearby point in $F_0^{(i-l)}$.
Then $(R^n\xi_*\boldz)_p\cong H^n(\xi^{-1}(U),\boldz)$
for some sufficiently small open neighbourhood $U$ of $p$. If $q\in U$,
then we have a restriction map $H^n(f^{-1}(p),\boldz)\cong
H^n(f^{-1}(U),\boldz)\rightarrow H^n(f^{-1}(q),\boldz)$, which we call
the specialisation map. To understand what this map is,
connect $p$ and $q$ with a path
$\gamma:[0,1]\rightarrow F$, $\gamma(0)=q$ and $\gamma(1)=p$. Then
$\xi^{-1}(\gamma([0,1]))\cong [0,1]\times T^{i-l}\cup \{1\}\times T^i$,
with $\{1\}\times T^{i-l}$ sitting inside $\{1\}\times T^i$. Then the
specialisation map is just the restriction map
$H^n(T^i,\boldz)\rightarrow H^n(T^{i-l},\boldz)$. This should be viewed as a
map $(R^n\xi_*\boldz)_p\rightarrow (R^n\xi_*\boldz)_q$.

To understand this map in practice, we just need to see how
$T^{i-l}$ sits inside $T^i$. Explicitly, first let $p=\mu^{-1}(P_i)
\in F^{(3)}$, and let $q$ be a nearby point in 
$\mu^{-1}(\langle P_i,P_j\rangle)$. Then $\xi^{-1}(p)=T(\dual{N_i})$,
while if $\langle P_i,P_j\rangle$ is an edge of $\sigma_l$, then
$\xi^{-1}(q)$ can be identified with $T(\langle e_k,e_m\rangle)
\subseteq T(N_l)$ (where as usual $\{i,j,k,l,m\}=
\{0,\ldots,4\}$). Using $T_{li}$ to identify $N_l$ with $\dual{N_i}$,
we see the specialisation map is induced by the inclusion
$T(\langle e_k^*-e_l^*,e_m^*-e_l^*\rangle)\subseteq T(\dual{N_i})$.
If $r\in\mu^{-1}(\langle P_i,P_j,P_k\rangle)$ near the point
$p$, then the specialisation map $(R^n\xi_*\boldz)_q
\rightarrow (R^n\xi_*\boldz)_r$ is induced by the inclusion
$$T(\langle e_m^*-e_l^*\rangle)\subseteq T(\dual{N_i}).$$

We now note the following facts:

(a) $R^3\xi_*\boldz=\boldz_{F^{(3)}}$. In particular, $H^0(F,R^3\xi_*\boldz)
=\boldz^5$.

(b) There is an exact sequence
$$\exact{R^2\xi_*\boldz}{\bigoplus_{i<j<k} \boldz_{F_{ijk}-C_{ijk}}}
{\boldz_{F^{(3)}}},$$
where $\boldz_{F_{ijk}-C_{ijk}}$ is the sheaf supported on $F_{ijk}$
obtained by extension by zero of the constant sheaf $\boldz$ on $F_{ijk}-
C_{ijk}$. In particular, $H^0(F,R^2\xi_*\boldz)=0$ and $H^1(F,R^2\xi_*\boldz)$
is torsion-free.

Proof. First note that $(R^2\xi_*\boldz)|_{F_0^{(2)}}=
\bigoplus_{i<j<k} \boldz_{F_{ijk}-C_{ijk}}|_{F_0^{(2)}}$. Thus in particular,
$i_{2*}i_2^* R^2\xi_*\boldz=\bigoplus_{i<j<k} \boldz_{F_{ijk}-C_{ijk}}$.
The exact sequence given in (b) will be obtained from the natural
map $R^2\xi_*\boldz\rightarrow i_{2*}i_2^* R^2\xi_*\boldz$.
To find the kernel and cokernel of this map, one studies the map
on stalks. Let $p\in F^{(3)}$; without loss of generality we can take
$p=(1:0:0:0:0)=\mu^{-1}(P_0)$. Then
$(i_{2*}i_2^*R^2\xi_*\boldz)_p=\bigoplus_{i=1}^4 (R^2\xi_*\boldz)_{p_i}$,
where $p_i\in\mu^{-1}(\langle P_0,P_i\rangle)$ is a point near $p$, and the
map $\psi:(R^2\xi_*\boldz)_p
\rightarrow\bigoplus_{i=1}^4 (R^2\xi_*\boldz)_{p_i}$
is given by specialisation. By the discussion of specialisation
above, this coincides with the restriction map
$$H^2(T(\dual{N_0}),\boldz)\rightarrow \bigoplus_{i=1}^4 H^2(T_i,\boldz),$$
with 
$$\eqalign{T_1&=T(\langle e_2^*-e_4^*,e_3^*-e_4^*\rangle)\cr
T_2&=T(\langle e_1^*-e_4^*,e_3^*-e_4^*\rangle)\cr
T_3&=T(\langle e_1^*-e_4^*,e_2^*-e_4^*\rangle)\cr
T_4&=T(\langle e_2^*-e_1^*,e_3^*-e_1^*\rangle)\cr}$$
From this, one
sees the kernel of $\psi$ is zero and the cokernel of $\psi$ is $\boldz$.
This yields the desired
exact sequence. Now since 
$$H^i(F,\boldz_{F_{ijk}-C_{ijk}})=H^i_c(F_{ijk}-C_{ijk},\boldz),$$
we see from the exact sequence of (b) that $H^0(F,R^2\xi_*\boldz)=0$
and $H^1(F,R^2\xi_*\boldz)$ is torsion-free. $\bullet$

(c) There is an exact sequence
$$0\rightarrow R^1\xi_*\boldz\rightarrow\bigoplus_{i<j}\boldz_{F_{ij}-C_{ij}}
\rightarrow
\bigoplus_{i<j<k} \boldz_{F_{ijk}-C_{ijk}}\rightarrow \boldz_{F^{(3)}}
\rightarrow 0.$$
In particular, $H^0(F,R^1\xi_*\boldz)=H^1(F,R^1\xi_*\boldz)=0$ and
$H^2(F,R^1\xi_*\boldz)$ is torsion-free.

Proof. Again $(R^1\xi_*\boldz)|_{F_0^{(1)}}
=\bigoplus_{i<j} \boldz_{F_{ij}-C_{ij}}|_{F_0^{(1)}}$, so
$i_{1*}i_1^*R^1\xi_*\boldz=\bigoplus_{i<j}\boldz_{F_{ij}-C_{ij}}$.
We study the kernel and cokernel of the map $\varphi:R^1\xi_*\boldz\rightarrow
i_{1*}i_1^*R^1\xi_*\boldz$. Let $\F=\coker \varphi$.
Let $p$ be a point in, without loss of generality, $(F_{234}-C_{234})\cap
F_0^{(2)}$ near $(1:0:0:0:0)$. 
If $p_i\in F_{2\cdots \hat \imath\cdots
4}-F^{(2)}$ are points near $p$, $i=2,3,4$, then the map on stalks
$\varphi_p:
(R^1\xi_*\boldz)_p\rightarrow (i_{1*}i_1^*R^1\xi_*\boldz)_p$ coincides
with the specialisation map $H^1(\xi^{-1}(p),\boldz)\rightarrow
\bigoplus_{i=2}^4 H^1(\xi^{-1}(p_i),\boldz)$. 
Now $\xi^{-1}(p)$ is identified with $T(\langle e_2^*-e_4^*,
e_3^*-e_4^*\rangle)$. Under this
identification, $\xi^{-1}(p_2)=T(\langle e_3^*-e_4^*\rangle)$,
$\xi^{-1}(p_3)=T(\langle e_2^*-e_4^*\rangle)$, and $\xi^{-1}(p_4)=
T(\langle e_3^*-e_2^*\rangle)$. From this one sees that $\varphi_p$ 
is injective with cokernel $\boldz$. A similar calculation
for $p=(1:0:0:0:0)$ yields $\ker\varphi_p=0$.
In particular, one obtains an exact sequence
$$\exact{R^1\xi_*\boldz}{\bigoplus_{i<j}\boldz_{F_{ij}-C_{ij}}}
{\F}.\eqno{(4.3)}$$
Also, $\F|_{F_0^{(2)}}=\bigoplus_{i<j<k}
\boldz_{F_{ijk}-C_{ijk}}|_{F_0^{(2)}}$. Next one studies the map
$$\F\rightarrow i_{3*}i_3^*\F=\bigoplus_{i<j<k} \boldz_{F_{ijk}-C_{ijk}}.$$
We have to study the map on stalks at $p\in F^{(3)}$. Without loss of
generality take $p=(1:0:0:0:0)$. Choose $p_{0i}
\in F_{1\ldots\hat \imath\ldots 4}-F^{(3)}$ near $p$ and $p_{0ij}\in F_{1\ldots\hat 
\imath
\ldots\hat\jmath\ldots 4}-F^{(2)}$ near $p$. We then have a diagram
$$\matrix{0&\mapright{}&H^1(T(\dual{N_0}),\boldz)&\mapright{}&
\bigoplus_{1\le i<j\le 4} H^1(\xi^{-1}(p_{0ij}),\boldz)&\mapright{}&\F_p&
\mapright{}&0\cr
&&\mapdown{\varphi_1}&&\mapdown{\varphi_2}&&\mapdown{\varphi_3}&&\cr
0&\mapright{}&\bigoplus_{i=1}^4 H^1(\xi^{-1}(p_{0i}),\boldz)
&\mapright{}& \bigoplus_{i,j=1\atop i\not=j}^4
H^1(\xi^{-1}(p_{0ij}),\boldz)&\mapright{}&
\bigoplus_{i=1}^4 \F_{p_{0i}}&\mapright{}&0\cr}$$
Here the first row is obtained by taking the stalks of (4.3) at $p$,
and the second is obtained by taking the direct sum of the stalks of
the exact sequence (4.3) at $p_{0i}$. The map 
$\varphi_1$ is again specialisation,
and $\varphi_2$ is given by mapping $\alpha\in H^1(\xi^{-1}(p_{0ij}),\boldz)
\mapsto \alpha\oplus\alpha\in 
H^1(\xi^{-1}(p_{0ij}),\boldz)\oplus H^1(\xi^{-1}(p_{0ij}),\boldz)$.
Explicit calculation then shows that the map $\varphi_3:\F_p\rightarrow
(i_{3*}i_3^*\F)_p$ is injective with cokernel $\boldz$. We omit the
details. This shows that we have an exact sequence
$$\exact{\F}{\bigoplus_{i<j<k}\boldz_{F_{ijk}-C_{ijk}}}{\boldz_{F^{(3)}}}.
\leqno{(4.4)}$$
Pasting this in with (4.3) gives the desired long
exact sequence.

Finally, to show the desired facts about $H^k(F,R^1\xi_*\boldz)$, just keep in
mind that $$H^k(F,\boldz_{F_{ij}-C_{ij}})=H^k_c(F_{ij}-C_{ij},\boldz).$$
This is zero if $k=0$ or $1$ and is torsion-free is $k=2$, from which
it follows from the exact sequence of (c) that $H^k(F,R^1\xi_*\boldz)$
is zero if $k=0$ or $1$ and is torsion-free for $k=2$.
$\bullet$

(d) $\xi_*\boldz=\boldz$, $H^0(F,\boldz)=H^2(F,\boldz)
=H^3(F,\boldz)=\boldz$, and $H^1(F,\boldz)=0$.

Proof. The first statement follows from the fact that $\xi$ has connected
fibres. To compute the cohomology of $F$, use the spectral sequence
$$E_1^{pq}=H^q(F^{[p]},\boldz)\Rightarrow H^n(F,\boldz)$$
(see e.g. [14]). Displaying the entries of this sequence for
$0\le p,q\le 3$, we obtain
$$\matrix{0&0&0&0\cr
\boldz^5&\boldz^{10}&\boldz^{10}&0\cr
0&0&0&0\cr
\boldz^5&\boldz^{10}&\boldz^{10}&\boldz^5\cr}$$
Only zeroes occur to the right. It is standard that the bottom row
coincides with the \v Cech complex for the intersection graph of
$F$, which is the simplicial complex $\partial\Xi$ (see e.g. [14],
pg. 105). The
second row is a truncation of the same complex. Thus the $E_2$ and
$E_3$ term
of this spectral sequence, for $0\le p,q\le 3$, is
$$\matrix{0&0&0&0\cr
\boldz&0&\boldz^4&0\cr
0&0&0&0\cr
\boldz&0&0&\boldz\cr}$$
Since $\Pic(F)$ certainly contains a copy of $\boldz$ generated by
the hyperplane section, so does $H^2(F,\boldz)$, and thus
the differential $d:E^3_{02}\rightarrow E^3_{30}$ must be zero.
From this one sees the desired identity on cohomology. $\bullet$

Now to show (2)-(4). From (a)-(d), the Leray spectral sequence for $\xi$,
for $0\le p+q\le 3$, looks like 
$$\matrix{\boldz^5&&&\cr
0&M_1&&\cr
0&0&M_2&\cr
\boldz&0&\boldz&\boldz\cr}$$
where $M_1$ and $M_2$ are torsion-free. This makes it clear
that $H^2(X,\boldz)=\boldz$ and $H^3(X,\boldz)$ is torsion-free.
Furthermore, $H^2(F,\boldz)\cong H^2(X,\boldz)$. If $H$ is a hyperplane
section of $F$, then $H$ pulls back to a generator $\xi^*H$ of
$H^2(X,\boldz)$. Furthemore, $(\xi^*H)^3=5$ since $H^3=5$.
Finally, $\rank H^3(X,\boldz)=204$ since $\chi(X)=-200$.
This is, for example, easily computed by counting singular fibres:
$$\eqalign{\chi(X)&=\hbox{\# of type $(1,2)$ fibres}-\hbox{\# of type $(2,1)$
fibres}\cr
&=10\cdot 5-10\cdot 25\cr
&=-200.\cr}$$

(5) follows from Proposition 2.13.

(6) will follow from Theorem 2.17. Indeed, $H\cdot c_2(quintic)=50$,
so we just need to show that $\xi^*H\cdot p_1(X)=-100$. 
Since $p_1(X)=-2Crit(f)$ by Theorem 2.17, and 
$\xi(Crit(f))=\bigcup C_{lm}\subseteq F\subseteq \Pfour$,
we just need to show that $H\cdot C_{lm}=5$ for each $l$ and $m$,
with the orientation on $C_{lm}$ induced by the canonical orientation
on $Crit(f)$. To
do so, consider $C_{lm}\subseteq F_{lm}=\Ptwo$. Then
$F_{lkm}$ is a $\Pone$ contained inside $\Ptwo$,
and $F_{klm}\cap C_{lm}$ consists of five points. If we show
this intersection is transversal and positive, then $H\cdot C_{lm}=5$.
If we take coordinates $(z_i,z_j,z_k)$ on $F_{lm}$,
then $F_{klm}$ is given by $z_k=0$, and near a point $p\in F_{klm}\cap C_{lm}$
the set
$\mu(C_{lm})$ is a line segment in $\langle P_i,P_j,P_k\rangle$
meeting $\langle P_i,P_j\rangle$. Adjusting this line segment
a bit if necessary and perturbing $C_{lm}$, 
we can assume that near $p$, $C_{lm}$ is given
by $z_i=1$, $z_j=constant$, making it clear the intersection with $F_{klm}$ is 
transversal. If $C_{lm}$ has the same orientation as that
induced by the complex structure, then the intersection will be positive.

To see this, let $q\in C_{lm}$ be a point near $p$. Take a small neighbourhood
of $\mu(q)$, so that $U=D\times (0,1)$ with $\mu(C_{lm})\cap D
=\{0\}\times (0,1)$. Let $\gamma$ be a small loop in $D^*$ in the same
direction as 
$\gamma_{ij,k}$ as shown in Figure 4.1. If $D$ is given an orientation
$\eta_D$ so
that $\gamma$ is counterclockwise, and $(0,1)$ is given the orientation
$\eta_{(0,1)}$ in the direction
from $\mu(p)$ to $\mu(q)$, then the product orientation on
$D\times (0,1)$ is the orientation 
$\eta_D\wedge\eta_{(0,1)}=
(P_j-P_i)\wedge (P_m-P_i)\wedge (P_k-P_i)$ on $\sigma_l$. The
monodromy about $\gamma$ in the basis $e_m,e_j,e_k$ of
$N_l$ is $\pmatrix{1&1&0\cr 0&1&0\cr 0&0&1\cr}$.
Then if $g:X'\rightarrow D$ is the fibration obtained by
$X'=f^{-1}(D\times\{1/2\})/T(\langle e_k\rangle)\rightarrow D$,
then by
Example 2.15, the orientation $\eta_D\wedge e_m\wedge e_j$ is the
positive orientation on $X'$.
In addition, taking the orientation $\eta_{(0,1)}\wedge e_k$ on
$(0,1)\times S^1$, the product orientation on
$f^{-1}(U)=
X'\times (0,1)\times S^1$ is
$$\eqalign{\eta_D\wedge e_m\wedge e_j\wedge \eta_{(0,1)}\wedge e_k
&=\eta_D\wedge\eta_{(0,1)}\wedge e_j\wedge e_k\wedge e_m\cr
&=-(P_j-P_i)\wedge (P_k-P_i)\wedge (P_m-P_i)\wedge e_j\wedge e_k\wedge e_m,
\cr}$$ 
which from (4.2) is the chosen orientation
on $X$. Thus by Example 2.15, $Crit(f)\cap f^{-1}(U)$ is oriented
by the product orientation $\eta_{(0,1)}\wedge e_k$, as is
$C_{lm}\cap \mu^{-1}(U)$. But this is precisely the orientation
induced by the complex structure.
$\bullet$
\bigskip

Having constructed a well-behaved fibration $f:X\rightarrow
B$ on the quintic, we can now construct a dual $\check f:\check X\rightarrow
B$ from Corollary 2.2. 
It follows immediately from [7], Theorem 3.10, which applies
as $f$ and $\check f$ are $\boldz$ and $\boldz/n\boldz$-simple,
that $H^2(\check X,\boldz)=\boldz^{101}$ and $H^3(\check X,\boldz)
=\boldz^4$. This coincides with the cohomology of the mirror quintic.
To show that $\check X$ is diffeomorphic to the mirror quintic,
however, we need to show that the cubic form on 
$H^2(\check X,\boldz)$ coincides with that of the mirror
quintic, and this will require some work. We begin by 
reviewing the geometry of the mirror quintic.

Let $\Xi$ be the simplex in ${\bf R}^4$ with vertices $P_0=(-1,-1,-1,-1)$,
$P_1=(4,-1,-1,-1)$,
$P_2=(-1,4,-1,-1)$,
$P_3=(-1,-1,4,-1)$,
$P_4=(-1,-1,-1,4)$. Then as is well-known (see [2], 5.1),
$\Xi$ is a reflexive polyhedron and
$\P_{\Xi}$, the toric four-fold corresponding to $\Xi$, is $\Pfour$.
Also the dual
$\Xi^*$ is the convex hull of $(-1,-1,-1,-1),(1,0,0,0),\ldots,(0,0,0,1)$, and
$\P_{\Xi^*}$ is naturally identified with the hypersurface in $\Pfive$ with
equation $x_0\cdots x_4=x_5^5$. A general hyperplane section of $\P_{\Xi^*}$ is
then a singular model for the mirror quintic. Fix such a hyperplane
section $W'$. The singularities of $W'$ are induced
by the singularities of $\P_{\Xi^*}$ and are as follows: there are ten 
curves of singularities, $C_{ij}$, isomorphic to $\Pone$, corresponding to
the edges $\langle P_i,P_j\rangle$ of $\Xi$. These are generically
$\boldz_5$-quotient singularities. There are in addition 10 singular points
$P_{ijk}$, corresponding to the faces $\langle P_i,P_j,P_k\rangle$
of $\Xi$, which are $\boldz_5^2$-quotient singularities. Furthermore,
$C_{ij}\cap C_{ik}\cap C_{jk}=\{P_{ijk}\}$.

To resolve these singularities, note that $\P_{\Xi^*}$ can also be
thought of as the toric variety $Y_{\Sigma}$ associated to the fan
$\Sigma$ which is the union of cones over faces of $\Xi$ with vertex
the origin. A resolution of $\P_{\Xi^*}$ is then
obtained by choosing a suitable subdivision of $\Sigma$, or equivalently
by choosing a suitable triangulation of $\partial\Xi$. Such a resolution
will also yield a crepant resolution $\eta:W\rightarrow W'$.
There is much choice in doing so, but we will make these choices so
that each two-dimensional face of $\Xi$ 
is triangulated as follows:
$$\epsfbox{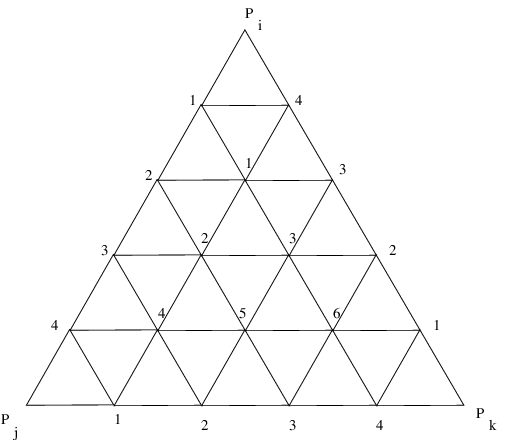}\leqno{\hbox{Figure 4.6}}$$
(See [15] or [2] for a more detailed description of the
method of resolution of $W'$.)
Each vertex of this triangulation gives rise to a toric divisor on
$W$. We let $L_0,\ldots,L_4$ be the divisors corresponding
to the points $P_0,\ldots,P_4$. These are just the proper transforms of
$W'\cap\{x_i=x_5=0\}\subseteq \{x_0\cdots x_4=x_5^5\}$. It is clear
from the choice of toric resolution that $\eta:L_i\rightarrow\Ptwo$
is an isomorphism. Let $E^l_{ij}$, $l=1,\ldots,4$ be the
divisors corresponding to the interior vertices of the edge
$\langle P_i,P_j\rangle$ as numbered in Figure 4.6. For each $i,j$, these
four divisors look like
$$\epsfbox{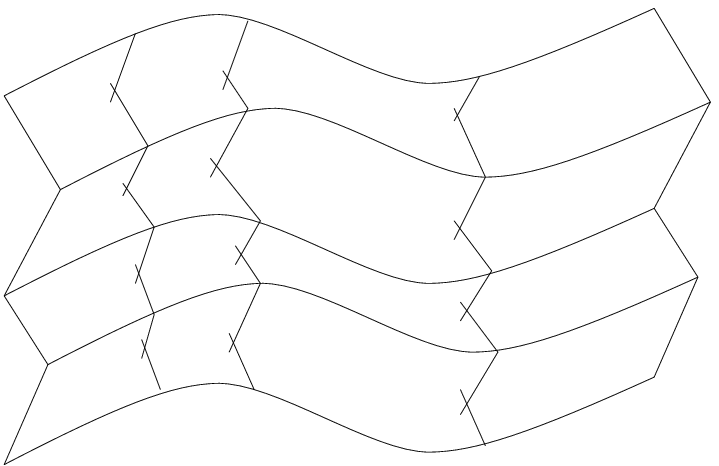}\leqno{\hbox{Figure 4.7}}$$
We will find it convenient to use the convention that $E_{ij}^l
=E_{ji}^{5-l}$.
Let $E^l_{ijk}$ be the divisors corresponding to the interior
vertices of $\langle P_i,P_j,P_k\rangle$, as numbered in Figure 4.6.
Then $\eta(E_{ijk}^l)=P_{ijk}$, $\eta(E^l_{ij})=C_{ij}$, and
the $E^l_{ijk}$, $E^l_{ij}$ are the exceptional divisors of
the resolution $\eta:W\rightarrow W'$. Furthermore, the pull-back
of a hyperplane section of $\P_{\Xi^*}$, giving the embedding of
$\P_{\Xi^*}$ in $\Pfive$, is just the anti-canonical
divisor of the toric resolution of $\P_{\Xi^*}$ restricted to
$W$, and as such is the sum of all the toric divisors.
We call this divisor
$$H:=\sum_i L_i+\sum_{i<j\atop l} E_{ij}^l+\sum_{i<j<k\atop l} E^l_{ijk}.$$

We summarize the intersection data:

\proclaim Proposition 4.2. On $W$, the following holds:
\item{(1)} $H^3=5$, $H\cdot E^l_{ij}
\cdot E^{l+1}_{ij}=1$, $H\cdot (E^l_{ij})^2=-2$, $H^2\cdot L_i=1$,
$H\cdot L_i\cdot E^1_{ij}=1$,
$H\cdot L_i^2=-3$.
\item{(2)}
$$\eqalign{
L_i\cdot (E^1_{ij})^2=1,\quad &L_i^2\cdot E^1_{ij}=-3,\cr
E^1_{ij}\cdot (E^{2}_{ij})^2=0,\quad& (E^1_{ij})^2\cdot E^{2}_{ij}=-2,\cr
E^2_{ij}\cdot (E^{3}_{ij})^2=-1,\quad& (E^2_{ij})^2\cdot E^{3}_{ij}=-1,\cr
E^3_{ij}\cdot (E^{4}_{ij})^2=-2,\quad& (E^3_{ij})^2\cdot E^{4}_{ij}=0,\cr}$$
and $L_i^3=9$, 
$(E^l_{ij})^3=5$. 
\item{(3)} Fixing $i$, $j$ and $k$, let 
$$\{D_1,\ldots,D_{21}\}
=\{E^1_{ijk},\ldots,E^6_{ijk},E^1_{ij},\ldots,E^4_{jk},L_i,L_j,L_k\}$$
be the 21 divisors corresponding to the 21 vertices in the triangulation
of Figure 4.6. Then $D_i^2\cdot D_j=-1$
if there is an edge $\langle i,j\rangle$ joining vertex $i$
and vertex $j$ and $\langle i,j\rangle$ is not an exterior
edge.
Also $D_i\cdot D_j\cdot
D_k=1$ if vertices $i,j,k$ are the vertices of a 2-simplex in Figure 4.6,
and $D_i^3=6$ if $D_i$ corresponds to an interior vertex.
\item{(4)}
All other intersection numbers involving $H$, the $L_i$'s,
the $E^l_{ij}$'s and the $E^l_{ijk}$'s are zero. 
\item{(5)}
If $c_2\in H^4(W,\boldz)$
is the second Chern class, then
$$\eqalign{H.c_2&=50,\cr
E^l_{ij}\cdot c_2&=2,\cr
E^l_{ijk}\cdot c_2&=0,\cr
L_i\cdot c_2&=-6.\cr}$$

Proof. This sort of computation is standard, and can be extracted,
for example, from [5], \S 5. We sketch it here: all values for cubes
of divisors follow from the formula $D^3=K_D^2$ for a smooth
divisor $D$ in a Calabi-Yau; all formulae involving $H$ are immediate
from the fact that $H$ induces the map $\eta:W\rightarrow W'\subseteq\Pfour$.
The formulae in (3) are clear from the toric
geometry, and the rest follows from,
for $A$ and $B$ two divisors meeting transversally in a Calabi-Yau threefold,
the formulae
$A\cdot B^2=(A\cap B)^2$, where the square is computed on $A$,
and $A\cdot B^2=(A\cap B)\cdot K_B$, calculated on $B$. Finally
the formulae for $c_2$ follow from Riemann-Roch.
$\bullet$

\proclaim Lemma 4.3. Let $M$ be an oriented, compact six-manifold with
$H^2(M,\boldz)=\boldz^{101}$, and let $S$ be a subset of 105 elements
of $H^2(M,\boldz)$, named so that
$$S=\{E^l_{ij}|0\le i<j\le 4,1\le l\le 4\}
\cup \{E^l_{ijk}|0\le i<j<k\le 4,1\le l\le 6\}
\cup \{L_0,\ldots,L_4\}$$
Suppose the elements of $S$ have
intersection numbers identical to those given in Proposition 4.2 (2)-(4).
Then the set $S$ generates $H^2(M,\boldz)$, and there is an isomorphism
$H^2(M,\boldz)\cong H^2(W,\boldz)$ preserving the cubic form, taking each
divisor in $S$ to the identically named divisor in $H^2(W,\boldz)$.

Proof.
Let $H=\sum_{D\in S} D$, so that $H$ has the same intersection numbers
as in Proposition 4.2. Let $L\subseteq H^2(M,\boldz)$ be the
integral span of the 101 divisors $\{H,E^l_{ij},E^l_{ijk}\}$.
Let $N\subseteq H^4(M,\boldz)$ be the span of
$\{D_1\cdot D_2|D_1,D_2\in S\}$. Let $L_{sat}\subseteq L\otimes_{\boldz}{\bf Q}$
be defined by 
$$L_{sat}=\{D\in L\otimes_{\boldz}{\bf Q}|D\cdot N\subseteq\boldz\}.$$

{\it Claim.} $L_{sat}/L\subseteq\boldz_5^4$, and $L_{sat}/L$ is generated
by $L_0,\ldots,L_4$.

Proof. Let $D\in L\otimes_{\boldz}{\bf Q}/L$, with
$$D=bH+\sum b_{ij}^l E_{ij}^l+\sum b_{ijk}^l E_{ijk}^l,$$
for $b,b_{ij}^l,b_{ijk}^l\in {\bf Q}/\boldz$.
Now $D\in L_{sat}/L$ if and only if $D\cdot D_1\cdot D_2\equiv 0\mod \boldz$
for all $D_1,D_2\in S$. We study what this condition tells us.

First focus on the divisors $E^l_{ij}$, $1\le l\le 4$. 
If
$$\eqalign{f_1&=E^1_{ij}\cdot(E^4_{ki}+E^1_{ijk})\cr
f_2&=E^2_{ij}\cdot(E^1_{ijk}+E^2_{ijk})\cr
f_3&=E^3_{ij}\cdot(E^2_{ijk}+E^4_{ijk})\cr
f_4&=E^4_{ij}\cdot(E^4_{ijk}+E^1_{jk})\cr}$$
then these are elements of $N$. If $a_{lm}=E^l_{ij}\cdot f_m$,
then from Figure 4.6 and the intersection numbers given in Proposition 4.2 (3),
the matrix $(a_{lm})$ is just the Cartan
matrix of $A_4$. (In $W$, $f_l$ would be the class of a fibre
of $\eta:E^l_{ij}\rightarrow C_{ij}$).
Furthermore, $f_1,\ldots,f_4$ intersect
all other divisors $H,E^l_{jk},E^l_{ijk}$ trivially.
Thus by inverting the Cartan matrix $A_4$, one sees that
$D\cdot f_m\equiv 0\mod \boldz$ for $m=1\ldots 4$ 
if and only if there exists a $b_{ij}\in {1\over 5}\boldz/\boldz$
such that $b^l_{ij}=lb_{ij}$.

Next focus on the simplex $\langle P_i,P_j,P_k\rangle$ 
as depicted in Figure 4.6. For each trapezoid
$$\epsfbox{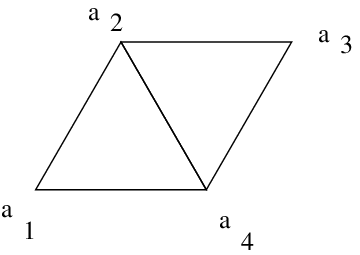}\leqno{\hbox{Figure 4.8}}$$
appearing in Figure 4.6, with corresponding divisors
$D_1,\ldots,D_4\in S$ and coefficients $a_1,\ldots,a_4$
in the expression for $D$, the requirement that 
$D\cdot D_2\cdot D_4\equiv 0\mod\boldz$ is equivalent
to $a_1+a_3\equiv a_2+a_4\mod\boldz$. Using this
relation, and given that $b_{ij}^l=lb_{ij}$, 
$b_{jk}^l=lb_{jk}$ and $b_{ki}^l=lb_{ki}$ for
$b_{ij},b_{jk},b_{ki}\in {1\over 5}\boldz/\boldz$,
one finds that the coefficients of all
classes corresponding to interior vertices are determined uniquely,
subject to the additional constraint that $b_{ij}+b_{jk}+b_{ki}
\equiv 0\mod\boldz$. (Here we use the convention
that $b_{ij}=-b_{ji}$.)

How many choices of $b_{ij}$'s are there which satisfy
this congruence? The set of all such $b_{ij}$'s is the kernel
of the map $\partial_1:\boldz_5^{10}\rightarrow\boldz_5^{10}$
defined by
$\partial_1((b_{ij})_{i<j})=
(b_{ij}+b_{jk}+b_{ki})_{i<j<k}$.
But this map $\partial_1$ is precisely the coboundary map
from 1-cochains to 2-cochains (with coefficients in
$\boldz_5$) of the two-skeleton $\Xi_2$ of $\Xi$.
The full cochain complex of the simplicial complex
$\Xi_2$ is
$$\matrix{\boldz_5^5&\mapright{\partial_0}&
\boldz_5^{10}&\mapright{\partial_1}&\boldz_5^{10},\cr}$$
and as $H^1(\Xi_2,\boldz_5)=0$, $\ker\partial_1=\im\partial_0$.
But $\ker\partial_0=H^0(\Xi_2,\boldz_5)=\boldz_5$, so 
$\ker\partial_1\cong\boldz_5^4$.

Next, we note that the coefficient of $H$ in $D$ is
completely determined from the relation $D\cdot
E^1_{ij}\cdot L_i\equiv 0\mod\boldz$ and the values of
the $b_{ij}$'s, as $H\cdot E^1_{ij}\cdot L_i=1$ for all
$i$ and $j$.
This makes it clear that the map $L_{sat}/L\rightarrow
\boldz_5^{10}$ given by $D\mapsto (b_{ij})_{0\le i,j\le 4}$
is injective, with image contained in $\ker\partial_1$.
Thus $L_{sat}/L\subseteq\boldz_5^4$.

Observe that these calculations imply that the 101 
generators of $L$ are in fact linearly independent.
If not, then in the calculations above, we would 
have found some non-zero 
choice of coefficients $b,b^l_{ij},b^l_{ijk}$
such that $\lambda D\cdot N
\equiv 0\mod\boldz$ for all $\lambda\in{\bf Q}$. We did not find such.
Thus $L\otimes_{\boldz}{\bf Q}=H^2(M,{\bf Q})$. Thus
in particular $L_p\in L_{sat}$, $0\le p\le 4$, and we can write
$$L_p=b^pH+
\sum lb_{ij}^pE^l_{ij}+\sum b^{l,p}_{ijk}
E^l_{ijk}$$
for suitable rational coefficients, $b^p,
b_{ij}^p,b^{l,p}_{ijk}$. 
It is then easy to see that $(b^p_{ij})_{i<j}\in
({1\over 5}\boldz)^{10}$ reduced modulo $\boldz$ is just $\partial_0(e_p)$,
where $e_0,\ldots,e_4$ are the generators of the
0-cochain group of $\Xi_2$ corresponding to the
vertices $P_0,\ldots,P_4$. Thus $L_0,\ldots,L_4$ generate
$L_{sat}/L$, proving the claim. $\bullet$

To prove the theorem, we just observe now that since
$L\otimes_{\boldz}{\bf Q}=H^2(M,{\bf Q})$,
$H^2(M,\boldz)\subseteq L_{sat}$. On the other hand
the claim shows $L_{sat}$ is generated by $S$, and hence
$L_{sat}\subseteq H^2(M,\boldz)$.
Thus $L_{sat}=H^2(M,\boldz)$. As the former is generated
by $S$, so is the latter. It is then clear that the
linear map $H^2(M,{\bf Q})\rightarrow H^2(W,{\bf Q})$
taking $H,E^l_{ij},E^l_{ijk}$ to the corresponding divisors
in $H^2(W,{\bf Q})$ is defined over $\boldz$, is an isomorphism over
$\boldz$, and under this
map the cubic forms coincide.$\bullet$

\proclaim Theorem 4.4. $\check X$ is diffeomorphic to $W$.

Proof. We again use Wall's classification theorem [20]. Because $\check X$
and $W$ are both simply connected and we already know that
$H^*(\check X,\boldz)\cong H^*(W,\boldz)$ as $\boldz$-modules,
we just need to show there is a choice of isomorphism
$H^2(\check X,\boldz)\cong H^2(W,\boldz)$ preserving the cubic form
and intersection with $p_1$. To do so, by Lemma 4.3,
we need to identify the relevant set $S\subseteq H^2(\check X,\boldz)$ 
of cohomology classes.

Before doing so, we need to study the geometry of $\check f$ in greater
detail. Let $U_i, U^i\subseteq\partial\Xi$ be the open sets defined
in the construction of $f$. Then as $f^{-1}(U_l)=T(N_l)\times U$, 
we can naturally identify $\check f^{-1}(U_l)$ with $T(\dual{N_l})\times U$.
Next, let $U_{lm}$ be an open set homeomorphic to ${\bf R}^3$ contained
in $\sigma_l\cup\sigma_m$, with $U_{lm}\cap\sigma_l\cap\sigma_m=
Int(\sigma_l\cap\sigma_m)$. Here
$\sigma_l\cap\sigma_m=\sigma:=\langle P_i,P_j,P_k\rangle$,
with $\{0,1,2,3,4\}=\{i,j,k,l,m\}$. Then the fibration
$\check f^{-1}(U_{lm})\rightarrow U_{lm}$ has a two-dimensional
space of cycles invariant under monodromy, and from the transpose inverse
of the description of the monodromy given in (4.1), this space of 
cycles is given by $\dual{N_l}\cap \dual{N_m}\subseteq \dual{N}$.
Thus there is a $T(\dual{N_l}\cap\dual{N_m})$ action on $\check f^{-1}(U_{lm})$,
with quotient map $q_{lm}:\check f^{-1}(U_{lm})\rightarrow
\check f^{-1}(U_{lm})/T(\dual{N_l}\cap\dual{N_m})
=S^1\times U_{lm}$. From the construction
method in Proposition 2.9, it is clear that the image of $Crit(\check f)$
under $q_{lm}$ can be identified with $\Delta':=\{p\}\times
(\Delta\cap U_{lm})$ for some $p\in S^1$. The $T(\dual{N_l}\cap
\dual{N_m})$-bundle
$$q_{lm}^{-1}(S^1\times U_{lm}-\Delta')\rightarrow S^1\times U_{lm}-\Delta'$$
has Chern class $$c_1\in H^3_{\Delta'}(S^1\times U_{lm},\dual{N_l}\cap\dual{N_m})
=H_1^{BM}(\Delta',\dual{N_l}\cap \dual{N_m}).$$ Using the known monodromy
about each edge of $\Delta\cap U_{lm}$,
we can describe $c_1$ up to
sign as a 1-chain on $\Delta'\cup\{\infty\}$.
Each edge of $\Delta'$ parallel (and in the same direction)
to the oriented edge $\langle P_{ijk},P_{ij}\rangle$ of $\Gamma$
has coefficient $e_j^*-e_i^*$; each edge of $\Delta'$ parallel
to $\langle P_{ijk},P_{ik}\rangle$ of $\Gamma$ has coefficient
$e_i^*-e_k^*$; and each edge of $\Delta'$ parallel to
$\langle P_{ijk},P_{jk}\rangle$ has coefficient $e_k^*-e_j^*$.

On the other hand, let $N'\cong\boldz^3\subseteq\dual{N}$ be a lattice with 
basis $e_i^*,e_j^*,e_k^*$, and let $\sigma$ be the cone in 
$N'\otimes_{\boldz}{\bf R}$ generated by $e_i^*$, $5e_j^*-4e_i^*$
and $5e_k^*-4e_i^*$. Then $\sigma$ can be subdivided to obtain
a fan $\Sigma$ which is a cone over the diagram in Figure 4.6. In
particular, $Y_{\Sigma}\rightarrow Y_{\sigma}$ is a crepant resolution
of $Y_{\sigma}$, as in \S 3. Also, in the notation of \S 3,
$N'_{m_0}=\dual{N_m}\cap\dual{N_l}$. Let $q_1:Y_{\Sigma}\rightarrow
{\bf C}\times{\bf R}^2$ be the $T(N_{m_0}')$-fibration
constructed in \S 3, $\Delta'_{\Sigma}\subseteq {\bf C}\times {\bf R}^2$
the locus over which $q_1$ fails to be a $T^2$-bundle. By construction,
$\Delta'_{\Sigma}$ is homeomorphic to $\Delta'$, and by Proposition
3.3, the first Chern class $c_1^{\Sigma}\in H_1^{BM}(\Delta'_{\Sigma},
\dual{N_m}\cap\dual{N_l})$ of the $T^2$-bundle $q_1:Y_{\Sigma}-
q_1^{-1}(\Delta'_{\Sigma})\rightarrow {\bf C}\times{\bf R}^2-\Delta'_{\Sigma}$
coincides with $c_1\in H_1^{BM}(\Delta',\dual{N_m}\cap\dual{N_l})$.
Now choose, for some point $p'\in S^1-\{p\}$, a homeomorphism
$\varphi_{lm}:(S^1-\{p'\})\times U_{lm}\rightarrow {\bf C}\times{\bf R}^2$
so that $\varphi^{-1}_{lm}(\{0\}\times {\bf R}^2)=\{p\}\times Int(\sigma)$ and
$\varphi_{lm}^{-1}(\Delta'_{\Sigma})=\Delta'$. Then this induces a homeomorphism
$$\psi_{lm}:q_{lm}^{-1}((S^1-\{p'\})\times U_{lm})\rightarrow Y_{\Sigma}.$$

Now $Y_{\Sigma}$ contains a number of toric divisors which are
submanifolds of $Y_{\Sigma}$, and these yield (possibly open) submanifolds
of $\check X$ via $\psi_{lm}$. These can be identified
as follows. There is one connected component of $Int(\sigma)-\Delta$
for each vertex $\tau$ of Figure 4.6 (see Figure 4.3). The component
containing $\tau$ we will call $\O_{\tau}\subseteq U_{lm}$. 
Then $q_{lm}^{-1}(\{p\}\times\overline{\O_{\tau}})$ is a smooth submanifold
of $\check f^{-1}(U_{lm})$, and $\psi_{lm}(q_{lm}^{-1}(\{p\}\times\overline
{\O_{\tau}}))$ 
is the toric divisor corresponding
to $\tau$ (modulo orientation, to be discussed shortly).

Of these 21 submanifolds, 6 correspond to interior vertices, and these
submanifolds are compact. We label them $E^p_{ijk}$, $1\le p\le 6$,
as in Figure 4.6. An additional 12 correspond to the interior vertices of
the three edges in Figure 4.6. These are non-compact submanifolds
of $\check X$, and we label these $E^p_{ij,ijk}, E^p_{ik,ijk}$ and
$E^p_{jk,ijk}$, $1\le p\le 4$, as numbered in Figure 4.6. Finally,
the three exterior vertices, $P_i,P_j$ and $P_k$, correspond to non-compact
submanifolds $L_{i,ijk}, L_{j,ijk}$ and $L_{k,ijk}$ of $\check X$.

Next we want to orient these submanifolds.
If $D_{\tau}$ is one of the submanifolds introduced above, which is a 
$T(\dual{N_l}\cap \dual{N_m})$-fibration over $\overline{\O_{\tau}}$, 
then $e_j^*-e_i^*, e_k^*-e_i^*$ can be taken as a basis for
$\dual{N_l}\cap \dual{N_m}$, and we can orient $D$ via
$$(P_j-P_i)\wedge (e_j^*-e_i^*)\wedge (P_k-P_i)\wedge (e_k^*-e_i^*).$$
This orientation is independent of the ordering of $i,j$
and $k$. By Remark 3.4, these orientations coincide with
the orientations induced by the complex structure on $\psi_{lm}(D_{\tau})
\subseteq Y_{\Sigma}$.

We now wish to glue together the three submanifolds $E^p_{ij,ijk}$,
$E^p_{ij,ijl}$ and $E^p_{ij,ijm}$. These are sitting over the 
shaded region in $B$ depicted in Figure 4.9:
$$\epsfbox{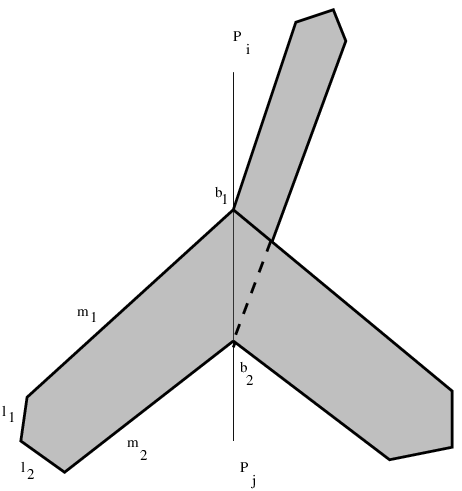}\leqno{\hbox{Figure $4.9$}}$$
For notational simplicity, we refer to these three submanifolds as
$E_k$, $E_l$ and $E_m$. In this diagram $b_1$ and $b_2$ are points
on $\langle P_i,P_j\rangle$ where three branches of $\Delta$ intersect.

First we look at a local model in a neighbourhood $V\subseteq B$
of the line segment $[b_1,b_2]\subseteq\langle P_i,P_j\rangle$.
This neighbourhood can be written as $V\cong D\times (0,1)$, $D$ a disk,
so that $V\cap \langle P_i,P_j\rangle=\{0\}\times (0,1)$,
$b_1,b_2\in (0,1)$, and there is a ``Y'' shape $\Delta'\subseteq
D$ with vertex at 0 such that $\Delta\cap V=\Delta'\times \{b_1,b_2\}$,
and $V\cap\check f(\bar E_k\cup \bar E_m\cup \bar E_n)=\Delta'\times
[b_1,b_2]$. Then by the construction of Example 2.6, (4), 
there is an $S^1$-fibration
$\pi:\check f^{-1}(V)\rightarrow T^2\times V$, with
$\pi(Crit(\check f)\cap \check f^{-1}(V))=S_1\cup S_2$, with
$S_1$ and $S_2$ pairs of pants sitting over $\Delta'\times \{b_1\}$
and $\Delta'\times \{b_2\}$ respectively. Furthermore, we can
assume that there is an $S\subseteq T^2\times\Delta'\subseteq T^2\times
D$ such that
$S_1=S\times \{b_1\}$ and $S_2=S\times \{b_2\}$. Let $T=
S\times [b_1,b_2]\subseteq T^2\times V$. Then $\pi^{-1}(T)$ is a submanifold of
$\check f^{-1}(V)$, and $E_k\cap \check f^{-1}(V)$ can be glued up
to $\pi^{-1}(T)$, possibly by perturbing $E_k$ a bit. The same goes
for $E_l$ and $E_m$. Thus using $\pi^{-1}(T)$, we can glue together
$E_k,E_l$ and $E_m$ to obtain a closed submanifold without
boundary $E^p_{ij}$ of $\check X$.

We need to make sure that this gluing can be done compatibly
with the choice of orientations made on these submanifolds.
To do this, we need to check that with this choice of orientation,
$\partial\bar E_k+\partial\bar E_l+\partial\bar E_m$ is homologous
to zero. First consider $\partial\bar E_k$. This is a fibration
over the line segment $[b_1,b_2]\subseteq\langle P_i,P_j\rangle$, with
fibres over points of $(b_1,b_2)$
homologous to the torus $T(\langle e_j^*-e_i^*,e_k^*-e_i^*\rangle)
\subseteq T(\dual{N_l})$, using a path from a chosen base-point
$b_l\in U_l$ to $b\in (b_1,b_2)$ to identify the fibre
$X_b$ with $T(\dual{N_l})$. The orientation $\eta_{\partial\bar E_k}$
on $\partial\bar E_k$ is induced by the orientation
$\eta_{\bar E_k}$ on $\bar E_k$ by the rule that
$\nu\wedge\eta_{\partial\bar E_k}=\eta_{\bar E_k}$, where $\nu$ is
an outward pointing normal vector on $\partial\bar E_k$.
Thus the orientation on $\partial\bar E_k$ is
$$(P_i-P_j)\wedge (e_j^*-e_i^*)\wedge (e_k^*-e_i^*),$$
as $P_i-P_k$ is an outward pointing normal vector.
Similarly, $\partial\bar E_l\cap X_b=T(\langle e_j^*-e_i^*,
e_l^*-e_i^*\rangle)\subseteq T(\dual N_k)$, (using
a different basepoint $b_k\in U_k$) and $\partial\bar E_l$ has orientation
$$(P_i-P_j)\wedge (e_j^*-e_i^*)\wedge (e_l^*-e_i^*).$$
Finally $\partial\bar E_m\cap X_b=T(\langle e_j^*-e_i^*,e_m^*-e_i^*\rangle)
\subseteq
T(\dual{N_l})$, using the same basepoint $b_l\in U_l$ and
path joining $b_l$ and $b$, and $\partial\bar E_m$ has orientation
$$(P_i-P_j)\wedge (e_j^*-e_i^*)\wedge (e_m^*-e_i^*).$$

To compare these three boundaries, we choose a path joining $b_k$
and $b_l$, and identify $T(\dual{N_k})$ with $T(\dual{N_l})$
using this path. Choose this path to pass from $U_k$ into
$U^i$ into $U_l$, thus identifying $\dual{N_k}$ and $\dual{N_l}$
via $T_{il}\circ T_{ik}^{-1}$. Now $T_{il}\circ T_{ik}^{-1}(e^*_j
-e_i^*)=e_j^*-e_i^*$, and $T_{il}\circ T^{-1}_{ik}(e_l^*-e_i^*)
=-(e_j^*-e_i^*)-(e_k^*-e_i^*)-(e_m^*-e_i^*)$. Thus the orientation on 
$\partial\bar E_l$ is 
$$(P_i-P_j)\wedge (e_j^*-e_i^*)\wedge (-(e_k^*-e_i^*)-(e_m^*-e_i^*)).$$
Since
$(e_j^*-e_i^*)\wedge (e_k^*-e_i^*)$, $(e_j^*-e_i^*)\wedge (e_m^*-e_i^*)$
and $(e_j^*-e_i^*)\wedge (-(e_k^*-e_i^*)-(e_m^*-e_i^*))$ 
add up to zero in $H_2(X_b,\boldz)$, it is then clear that the orientations
on $E_k$, $E_l$ and $E_m$ are compatible with the gluing process.
Thus $E^p_{ij}$ is oriented by these compatible orientations.

Next for each $i$ we glue together the 6 submanifolds $\{L_{i,ijk}|
\{j,k\}\subseteq \{0,\ldots,4\}-\{i\}\}$. The triple $L_{i,ijk}, L_{i,ijl}$
and $L_{i,ijm}$ glue together as above over $\langle P_i,P_j\rangle$, 
with compatible orientations, so all six of these submanifolds
can be glued together to obtain $L_i$ which is a manifold
except possibly over $\check f^{-1}(P_i)$. We need to clarify this
structure. I claim in fact that this gluing can be done in such a way
that there is an orientation preserving homeomorphism of $L_i$ with $\Ptwo$. 

To do this, we divide $\Ptwo$ into six pieces which can be identified
with the $L_{i,ijk}$'s. Without loss of generality take $i=4$,
and identify $\Ptwo$ with $\{z_0+\cdots+z_3=0\}\subseteq\Pthree$.
Let $L_{jk}\subseteq\Ptwo$ be defined by
$$L_{jk}=\{(z_0,\ldots,z_3)\in\Ptwo| |z_j|,|z_k|<|z_l|,|z_m|, \{j,k,l,m\}
=\{0,1,2,3\}\}.$$
Then it is not difficult to see that $L_{jk}$ can be identified with
$L_{i,ijk}$ in such a way that
$(\{z_j=0\}\cap L_{jk})\cup (\{z_k=0\}\cap L_{jk})$ is identified
with $Crit(\check f)\cap L_{i,ijk}$.
We leave as an exercise to the reader to verify that the $L_{i,ijk}$'s
can then be glued together inside $\check X$ in the same way
that the $L_{jk}$'s glue together to form $\Ptwo$. This gives a
homeomorphism between $L_i$ and $\Ptwo$ which identifies
$Crit(\check f)$ with the four lines $\bigcup_{j=0}^3 \{z_j=0\}\subseteq
\Ptwo$.

We have now identified the relevant set $S\subseteq H^2(\check X,\boldz)$
of elements, namely the Poincar\'e dual
cohomology classes of the oriented manifolds
$E^p_{ij}$, $E^p_{ijk}$ and $L_i$ defined above. We will also
think of $S$ as the set of these oriented manifolds. Note that while
the cohomology classes in $S$ still depend on a choice of orientation
on $\check X$, the intersection numbers $D_1\cdot D_2\cdot D_3
\in\boldz$, for $D_1,D_2,D_3\in S$, do not depend on the choice
of orientation: a change of orientation changes the signs of $D_1,D_2$
and $D_3$, and the sign of the identification $H^6(\check X,\boldz)
\cong\boldz$. Thus in particular, if we wish to compute one of the intersection
numbers corresponding to those described in Proposition 4.2 (3), this 
intersection number can in fact be calculated by performing this
intersection in $Y_{\Sigma}$. It is then clear that all intersection
numbers described in Proposition 4.2 (3) are valid in $\check X$.
The same is also true for those intersection numbers given in
Proposition 4.2, (4): all of the intersection numbers
$D_1\cdot D_2\cdot D_3$ not mentioned explicitly in Proposition 4.2 (2) 
and (3) are zero because the intersection
of the corresponding submanifolds is empty.

Thus there remains the intersection numbers in Proposition 4.2, (2);
these require more work.
It is easy to see that all intersection numbers given in Proposition
4.2, (2) follow from the following five claims:

{\it Claim 1.} 
$$(E^{p-1}_{ij})^2\cdot E^p_{ij}+E^p_{ij}\cdot (E^{p+1}_{ij})^2=-3,
\quad 1\le p\le 4,$$
with the convention that $E^0_{ij}=L_i$ and $E^5_{ij}=L_j$.

Proof. As toric divisors, $E^p_{ij,ijk}$,$E^p_{ij,ijl}$ and $E^p_{ij,ijm}$
are ruled surfaces over the affine line ${\bf C}$, each with one
reducible fibre consisting of a union of two $-1$ curves. 
One way to see this structure topologically is as follows.
Recall $E^p_{ij,ijk}=q_{lm}^{-1}(\{p\}\times\overline{\O_{\tau}})$ for
suitable $\tau$. The set
$\overline{\O_{\tau}}$ can be seen as one of the three sheets in Figure 4.9.
Write the boundary of $\overline{\O_{\tau}}$ as the union of five line segments:
$l_1$ and $l_2$ being the image of the two $-1$ spheres under the map
$q_{lm}:E^p_{ij,ijk}\rightarrow\{p\}\times\overline{\O_{\tau}}$, 
$m_1$ and $m_2$
being the upper and lower edges of $\overline{\O_{\tau}}$, as depicted in
Figure 4.9, and $[b_1,b_2]$ being the fifth line segment.
Now fibre $\{p\}\times\overline{\O_{\tau}}$ in line segments by a map
$\gamma:\{p\}\times\overline{\O_{\tau}}
\rightarrow [0,1]$ so that $\gamma^{-1}(0)=
l_1\cup l_2$, $\gamma^{-1}(x)$ for $x>0$ is a line segment 
connecting $m_1$ and $m_2$, and $\gamma^{-1}(1)=[b_1,b_2]$. 
Clearly $q_{lm}^{-1}(\gamma^{-1}(0))$
is the union of the two $-1$ spheres, while $q_{lm}^{-1}(\gamma^{-1}(x))$
for $x\in (0,1)$ is easily seen to be $S^1\times S^2$, in such a way that
$q_{lm}:\{q\}\times S^2\rightarrow \gamma^{-1}(x)$ is the
standard projection of a two-sphere to a line segment for any $q\in S^1$.
Thus the map $\gamma\circ q_{lm}:E^p_{ij,ijk}\rightarrow [0,1)$
factors through an $S^2$-fibration $\rho:E^p_{ij,ijk}\rightarrow D$,
where $D$ is a disk. It is clear that these three fibrations on
$E^p_{ij,ijk}$, $E^p_{ij,ijl}$ and $E^p_{ij,ijm}$ glue together to 
yield a map $\rho:E^p_{ij}
\rightarrow S^2$ whose general fibre is an $S^2$, 
and three fibres being the union
of two $-1$ spheres. In addition, $C_1=
E^{p-1}_{ij}\cap E^p_{ij}$ and $C_2=E^p_{ij}\cap E^{p+1}_{ij}$
are sections of $\rho$. Let $\beta:E^p_{ij}\rightarrow E$ be the
blow-down of three $-1$ spheres, one on each singular fibre. Then
$\rho$ induces a map $\rho':E\rightarrow S^2$ with $\rho=\rho'\circ\beta$,
and $E$ is an $S^2$ bundle over $S^2$, hence homeomorphic to 
$\P(\O_{\Pone}\oplus\O_{\Pone}(n))$ for some $n$. In particular,
$\beta(C_1)$ and $\beta(C_2)$ are still disjoint sections of $\rho'$,
so $\beta(C_1)^2+\beta(C_2)^2=0$. Thus, on the blow-up $E^p_{ij}$,
$C_1^2+C_2^2=-3$. But $C_1^2=E^p_{ij}\cdot (E^{p-1}_{ij})^2$ and
$C_2^2=E^p_{ij}\cdot (E^{p+1}_{ij})^2$, hence Claim 1 follows. $\bullet$

{\it Claim 2.} $(E^p_{ij})^3=5$, $1\le p\le 4$.

Proof. Before starting, let us note that by Example 2.16, the
canonical orientation on $Crit(\check f)$ can be described as follows.
Every component of $Crit(\check f)$ is the intersection $D_1\cap D_2$
of some two oriented submanifolds $D_1,D_2\in S$. The canonical
orientation of this component of $Crit(\check f)$ is
that of the oriented intersection $D_1\cap D_2$. As remarked in
Example 2.16, this is independent of the choice of orientation on
$\check X$.

Now perturb $E^p_{ij}$ to $\tilde E^p_{ij}$ by perturbing
$\check f(E^p_{ij})$. This can be done as pictured in Figure 4.10:
$$\epsfbox{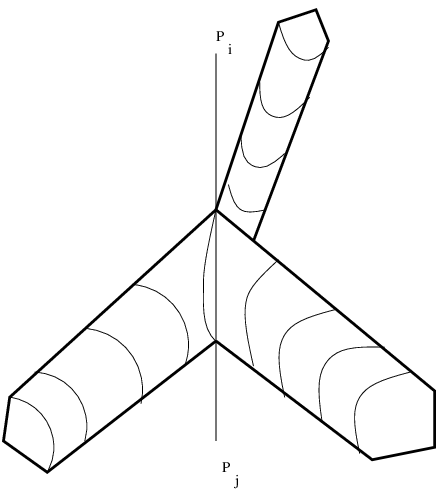}\leqno{\hbox{Figure $4.10$}}$$
We keep $\partial\check f(E^p_{ij})$ fixed but move the interior.
This can be done in such a way that $\check f(E^p_{ij})\cap \check f(\tilde
E^p_{ij})$ consists of the union of
$\check f(E^p_{ij})\cap \Delta$ and a line segment
$s$ with endpoints $b_1$ and $b_2$.
It is clear that the contribution to $E^p_{ij}\cap \tilde
E^p_{ij}$ sitting over $\check f(E^p_{ij})\cap \Delta$
is $Crit(\check f)\cap E_{ij}^p=C_1\cup C_2\cup
e_1\cup\cdots\cup e_6$, where $C_1$ and $C_2$ are the two sections
of $\rho$ as in the proof of Claim 1, and $e_1,\ldots,e_6$ are the
6 toric $-1$ $\Pone$'s. The orientation on these components of 
the oriented intersection $E^p_{ij}\cap \tilde E^p_{ij}$ is the
opposite orientation to the canonical orientation on $Crit(\check f)$.
One way to see this is to do a similar procedure with one of the divisors
$E^p_{ijk}$ fibering over $\overline{\O_{\tau}}$, 
$\tau$ a vertex in the interior
of the triangulation of Figure 4.6. The self-intersection
$(E^p_{ijk})^2$ can be computed by perturbing
$\overline{\O_{\tau}}$, to obtain $\tilde E^p_{ijk}$, so that
$E^p_{ijk}\cap\tilde E^p_{ijk}=Crit(\check f)\cap E^p_{ijk}$. 
However, from the toric geometry of $Y_{\Sigma}$ it is known
that $(E^p_{ijk})^2$ can be represented by the
two-cycle supported on $Crit(\check f)\cap E^p_{ijk}$ with
the opposite orientation to the canonical orientation. Since
the $E^p_{ij}$'s were oriented in the same manner
as the $E^p_{ijk}$'s, the same must be true for the former
self-intersection.

The contribution to $E^p_{ij}\cap\tilde E^p_{ij}$ sitting over the line
segment $s$ is seen to
be an $S^2$, homologous to a fibre $f$
of $\rho$, oriented so that $f\cdot C_1=f\cdot C_2=1$ with
$C_1,C_2$ having the canonical orientation. Thus as a cohomology
class in $H^2(E^p_{ij},\boldz)$, with $e_1,\ldots,e_6$ also
oriented by the canonical orientation, $E^p_{ij}\cap \tilde E^p_{ij}
=-C_1-C_2-\sum e_i+f$. Now $(E^p_{ij})^3$ coincides with
$(-C_1-C_2-\sum e_i+f)^2$ computed on $E^p_{ij}$. This
is
$$\eqalign{&C_1^2+C_2^2+\sum e_i^2+f^2+2\sum_{i<j} e_i\cdot e_j
+2\sum e_i\cdot (C_1+C_2)-2f\cdot C_1-2f\cdot C_2\cr
=& C_1^2+C_2^2-6+0+6+12-2-2\cr
=&C_1^2+C_2^2+8.\cr}$$
But in the proof of Claim 1 it was shown that $C_1^2+C_2^2=-3$,
hence $(E^p_{ij})^3=5$. $\bullet$

{\it Claim 3.} $E^{p-1}_{ij}\cdot (E^{p}_{ij})^2+
(E^p_{ij})^2\cdot E^{p+1}_{ij}=-1$, $1\le p\le 4$.

Proof. Recall the index $I(M)$ of an oriented 4-dimensional manifold $M$
is $q_+-q_-$, where the intersection form on $H^2(M,{\bf Q})$ has
signature $(q_+,q_-)$. By [10], 9.2, page 88, $I(E^p_{ij})
={1\over 3}(-(E^p_{ij})^3+p_1(\check X)\cdot E^p_{ij})$. Now from the
proof of Claim 1 it is clear that $I(E^p_{ij})=1-4=-3,$ so by Claim 2,
$p_1(\check X)\cdot E^p_{ij}=-4$. On the other hand, 
$p_1(\check X)=-2Crit(\check f)$. It is easy to see from the
description of the canonical orientation on $Crit(\check f)$
at the beginning of the proof of Claim 2 that
$$p_1(\check X)\cdot E^p_{ij}=
-2(E^{p-1}_{ij}\cdot E^p_{ij}+E^p_{ij}\cdot E^{p+1}_{ij})\cdot E^p_{ij}
-2(-6+9)$$
where the second term is the contribution from the six toric
$-1$ $\Pone$'s on $E^p_{ij}$ and nine additional curves of
$Crit(\check f)$ not contained in $E^p_{ij}$, but intersecting it. From this one
derives the formula of Claim 3. $\bullet$

{\it Claim 4.} $L_i\cdot (E^1_{ij})^2=1$.

Proof. $L_i$ is homeomorphic to $\Ptwo$, and by construction of this
homeomorphism,
$L_i\cap E^1_{ij}$ is a generator of $H^2(L_i,\boldz)$. Thus
in $L_i$, $(L_i\cap E^1_{ij})^2=1$. But $(L_i\cap E^1_{ij})^2
=L_i\cdot (E^1_{ij})^2$. $\bullet$

{\it Claim 5.} $L_i^3=9$.

Proof. As in the proof of Claim 3, 
$$1=I(L_i)={1\over 3}(-L_i^3+p_1(\check X)\cdot L_i)$$
so 
$$\eqalign{L_i^3&=p_1(\check X)\cdot L_i-3\cr
&=-2Crit(\check f)\cdot L_i-3\cr
&=-2\left( \sum_{j\not=i} E^1_{ij}\cdot L_i^2+6\right)-3.\cr}$$
Here the contribution of 6 is for the six components of $Crit(\check f)$
intersecting $L_i$ but not contained in it. It follows from Claims 1, 3 
and 4 that $E^1_{ij}\cdot L_i^2=-3$, and hence
$$L_i^3=-2(4\cdot(-3)+6)-3=9.\bullet$$ 

Finally, it is clear that all intersections of $p_1(\check X)$ agree
with those given in Proposition 4.2, (5) (with $p_1(W)=
-2c_2(W)$). Indeed, these were already calculated in the proof of
Claims 3 and 5 for $E^l_{ij}$ and $L_i$, and one sees easily that
$E^l_{ijk}\cdot Crit(\check f)=0$.
Thus we can conclude, by Wall's theorem, that $\check X$ is diffeomorphic
to the mirror quintic $W$.
$\bullet$
\bigskip
{\it Remark 4.5.} We end this paper with the question:
what if we had chosen a different toric resolution
of the singular mirror quintic
$W'$? It is clear that our particular choice of resolution
was already present in the construction of our $T^3$-fibration
on the quintic, in Figure 4.3, as dictated by the construction of
\S 3. Now two different toric resolutions
differ by a sequence of toric flops. A toric flop is given by 
making the following change to a trapezoid in Figure 4.6:
$$\epsfbox{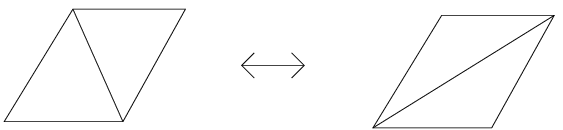}\leqno{\hbox{Figure $4.11$}}$$
Such a flop gives a resolution $W_2$ of $W'$. 
The point is that we can at the same time alter the fibration $f:X\rightarrow
B$ on the quintic to a fibration $f':X'\rightarrow B$, simply
by making a change in the discriminant locus pictured in Figure 4.3.
This change is represented by the following diagram:
$$\epsfbox{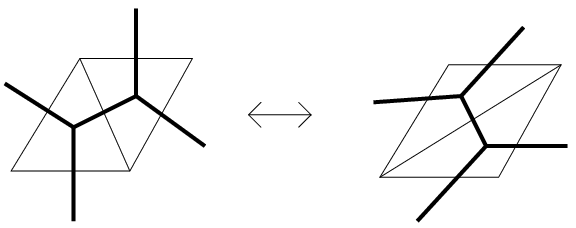}\leqno{\hbox{Figure $4.12$}}$$
It is clear that the new fibration $f':X'\rightarrow B$ so constructed
is still a $T^3$-fibration on the quintic. Nowhere in the details of the
proof that $X$ was diffeomorphic to the quintic did we make use of the detailed
structure of $\Delta$. Furthermore, by running through the proof that $\check
X$ is diffeomorphic to $W$, one finds that $\check X'$ is diffeomorphic
to $W_2$. Thus we see how a change in the fibration $f:X\rightarrow B$,
but {\it not} a change in $X$ itself, gives rise to two duals related by
a flop.

This is an interesting enough phenomenon that it is worth noting it
really is a local one. One can construct over a closed 3-ball $B$ two different
$T^3$-fibrations, $f_1:X_1\rightarrow B$ and $f_2:X_2\rightarrow B$
by applying the construction of Proposition 2.5 to $S_1,S_2\subseteq
\bar Y=T^2\times B$, where $S_1$ and $S_2$ are both homeomorphic
to $S^2-\{p_1,p_2,p_3,p_4\}$, and are fibred over the two discriminant
loci $\Delta_1,\Delta_2$ pictured
below, with classes in $T^2$ over each branch as labelled. Here
$\Delta_1\cap\partial B=\Delta_2\cap\partial B$, and $S_1$ and $S_2$
can be chosen so that $\partial\bar Y\cap S_1=\partial\bar Y\cap S_2$.
$$\epsfbox{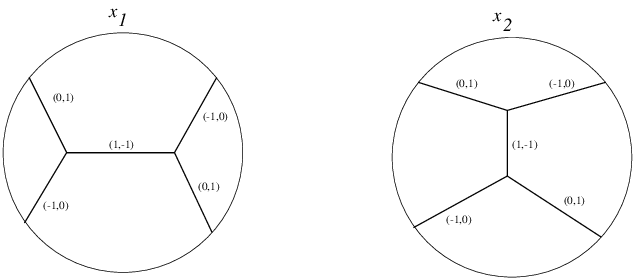}\leqno{\hbox{Figure $4.13$}}$$
There is in fact an isotopy  $I:\bar Y\times [0,1]
\rightarrow \bar Y\times [0,1]$ which is the identity on $\partial\bar Y$,
$I_0=id$ and $I_1(S_1)=S_2$.
This can be seen by first constructing
an isotopy  $I':S_1\times [0,1]\rightarrow \bar Y\times [0,1]$
with $I_0':S_1\rightarrow \bar Y$ the embedding of $S_1$ into
$\bar Y$, and $I_1':S_1=S_2\rightarrow\bar Y$ the embedding of $S_2$
into $\bar Y$.
One then applies
Corollary 2.4 of [11] to obtain $I$. In particular, $I_1$ induces a
homeomorphism
$\varphi:X_1\rightarrow X_2$
which is the identity on $\partial X_1=\partial X_2$. 

Now the duals of $f_1$ and $f_2$ can be constructed, giving
$\check f_1:\check X_1\rightarrow B$ and $\check f_2:\check X_2\rightarrow
B$. It is clear that $\check X_1$ and $\check X_2$ are related by a flop.
On the other hand, $\varphi:\partial X_1\rightarrow \partial X_2$ is
the identity, and this induces an isomorphism $\check\varphi:
\partial\check X_1\rightarrow\partial\check X_2$. However, it is possible
to show that this homeomorphism does not extend to a homeomorphism
$\check\varphi:\check X_1\rightarrow\check X_2$. We will not show this here, but
point out that we should expect this, for otherwise every
toric resolution of $W'$ would be homeomorphic, which we don't
expect to be the case. We also note that if we
naively tried to construct such a homeomorphism, we might try to proceed as 
follows. $\check X_1$ and $\check X_2$ are both expressed as $T^2$-fibrations
over $S^1\times B$, with these fibrations
failing to be $T^2$-bundles over $\{p\}\times \Delta_1$ and $\{p\}
\times\Delta_2\subseteq S^1\times B$. We might try to construct a homeomorphism
$S^1\times B$ taking $\{p\}\times \Delta_1$ to $\{p\}\times\Delta_2$, while
being the identity on $\partial(S^1\times B)$. This is clearly impossible,
because the two discriminant loci
$$\epsfbox{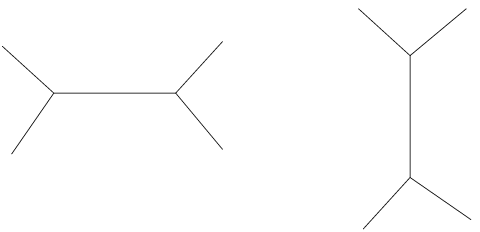}\leqno{\hbox{Figure $4.14$}}$$
are clearly not the same if the end-points are marked. On the other hand,
$X_1$ and $X_2$ were diffeomorphic because the two surfaces
$S_1$ and $S_2$ were the same:
$$\epsfbox{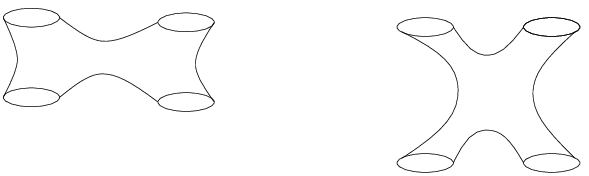}\leqno{\hbox{Figure $4.15$}}$$
This might remind the reader of one of the basic features underlying
string theory.

{\hd Bibliography}
\item{[1]} Aspinwall, P., Greene, B., and Morrison, D.,
``Calabi-Yau Moduli Space, Mirror Manifolds and Spacetime Topology
Change in String Theory,'' {\it Nuclear Phys. B}, {\bf 416}, 
(1994), 414--480.
\item{[2]} Batyrev, V., ``Dual Polyhedra and Mirror Symmetry
for Calabi-Yau Hypersurfaces in Toric Varieties,'' {\it J. Alg. Geom.},
{\bf 3}, (1994), 493--535.
\item{[3]} Bredon, G., {\it Sheaf Theory}, 2nd edition, Springer-Verlag,
1997.
\item{[4]} Chiang, T.-M., Klemm, A., Yau, S.-T. and Zaslow, E., 
``Local Mirror Symmetry: Calculations and Interpretations,'' preprint,
hep-th/9903053.
\item{[5]} Dais, D., ``Enumerative Combinatorics of Invariants of
Cetain Complex Threefolds with Trivial Canonical Bundle,'' {\it Bonner
Mathematische Schriften}, {\bf 279}, (1995).
\item{[6]} Gross, M., ``Special Lagrangian Fibrations I: Topology,'' 
in {\it Integrable Systems and Algebraic Geometry},
eds. M.-H. Saito, Y. Shimizu and
K. Ueno, World Scientific, 1998, 156--193.
\item{[7]} Gross, M., ``Special Lagrangian Fibrations II: Geometry,''
to appear in {\it Surveys in Differential Goemetry}.
\item{[8]} Gross, M., and Wilson, P.M.H., ``Mirror Symmetry via 3-tori for
a Class of Calabi-Yau Threefolds,'' 
{\it Math. Ann.}, {\bf 309}, (1997) 505--531.
\item{[9]} Harvey, R., and Lawson, H.B. Jr.,  ``Calibrated Geometries,'' {\it
 Acta
Math.} {\bf 148}, 47-157 (1982).
\item{[10]} Hirzebruch, F., {\it Topological Methods in Algebraic
Geometry,} Second corrected printing to Third edition, Springer-Verlag, 1978.
\item{[11]} Hudson, J.F.P., and Zeeman, E.C.,  ``On Combinatorial Isotopy,''
{\it Inst. Hautes \'Etudes Sci. Publ. Math.}, {\bf 19}, (1964) 69--94.
\item{[12]} Humphreys, J., {\it Linear Algebraic Groups,} Springer-Verlag,
1975.
\item{[13]} Lawson, H.B. Jr., and Michelsohn, M.-L., {\it Spin
Geometry}, Princeton University Press, 1989.
\item{[14]} Morrison, D., ``The Clemens-Schmid Exact Sequence
and Applications,'' in {\it Topics in Transcendental Algebraic
Geometry,} edited by P. Griffiths, Princeton University Press.
\item{[15]} Morrison, D., ``Mirror Symmetry and Rational
Curves on Quintic Threefolds: A Guide for Mathematicians,''
{\it J. Amer. Math. Soc.}, {\bf 6}, (1993), 223--247.
\item{[16]} Oda, T., {\it Convex Bodies and Algebraic Geometry}, Springer,
1988.
\item{[17]} Ruan, W.-D., ``Lagrangian Tori Fibration of Toric Calabi-Yau
Manifold I,'' preprint, math.DG/9904012.
\item{[18]} Ruan, W.-D., ``Lagrangian Tori Fibration of Toric Calabi-Yau
Manifold III: Symplectic Topological SYZ Mirror Construction for General
Quintics,'' preprint, math.DG/9909126.
\item{[19]} Strominger, A., Yau, S.-T., and Zaslow, E.,  ``Mirror Symmetry is
T-Duality,'' {\it Nucl. Phys.} {\bf B479}, (1996) 243--259.
\item{[20]} Wall, C.T.C., ``Classification Problems in Differential
Topology. V. On Certain $6$-manifolds,''  {\it Inv. Math.} {\bf 1},
(1966), 355--374.
\item{[21]} Zharkov, I., ``Torus Fibrations of Calabi-Yau Hypersurfaces in
Toric Varieties and Mirror Symmetry,'' preprint, alg-geom/9806091,
to appear in {\it Duke Math. J.}
\end

%% file: defs.tex
\font\hd=cmbx10 scaled\magstep1
\input cyracc.def
\newfam\cyrfam

\def\num{\global\advance\count10 by 1 \eqno(\the\count10)}

\def\Pfive{{\bf P}^5}
\def\Pfour{{\bf P}^4}
\def\Pthree{{\bf P}^3}

\def\Ptwo{{\bf P}^2}
\def\Pone{{\bf P}^1}
\def\P{{\bf P}}

\def\O{{\cal O}}

\def\T{{\cal T}}

\def\F{{\cal F}}

\def\N{{\cal N}}
\def\E{{\cal E}}

\def\HH{{\underline H}}

\def\EY{\E(1)|_Y}

\def\PEY{\P(\EY)}

\def\OP{\O_{\PEY}(1)}
\def\OP1{\O_{\Pone}}

\def\Pic{{\rm Pic}}

\def\coker{\mathop{\rm coker}}

\def\im{\mathop{\rm Im}}

\def\mod{\mathop{\rm mod}}
\def\ext{{\rm Ext}}

\def\hom{{\rm Hom}}

\def\rank{\mathop{\rm rank}}

\def\lhom{{\underline{\hom}}}

\def\boldz{{\bf Z}}

\def\dual#1{{#1}^{\scriptscriptstyle \vee}}

\def\exact#1#2#3{0\rightarrow#1\rightarrow#2\rightarrow#3\rightarrow0}

\def\mapright#1{\smash{
  \mathop{\longrightarrow}\limits^{#1}}}

\def\mapdown#1{\Big\downarrow
   \rlap{$\vcenter{\hbox{$\scriptstyle#1$}}$}}